\documentclass[12pt,epsfig,amsfonts]{amsart}
   \newtheorem{lemma}{Lemma}[section]
   \newtheorem{theorem}{Theorem}
   
   \newtheorem{prop}{Proposition}[section]
   \newtheorem{coro}{Corollary}[section]

\usepackage{epsfig}

\newtheorem{claim}{Claim}[section]
\newtheorem{sublemma}{Sublemma}[section]
\numberwithin{equation}{section}

\author{Qiudong Wang}
\address[Qiudong Wang]
{Department of Mathematics\\
    University of Arizona\\
    Tuscon, AZ 85721}
\email{dwang@math.arizona.edu}

\author{Ali Oksasoglu}
\address[Ali Oksasoglu]
{Honeywell Corporation\\
11100 N. Oracle Rd.\\
Tucson, AZ 85737}
\email{ali.oksasoglu@honeywell.com}

\title{On the dynamics of certain homoclinic tangles}

\begin{document}

\begin{abstract}
In this paper we study homoclinic tangles formed by transversal
intersections of the stable and the unstable manifold of a {\it
non-resonant, dissipative} homoclinic saddle point in periodically
perturbed second order equations. We prove that the dynamics of
these homoclinic tangles are that of {\it infinitely wrapped
horseshoe maps} (See Section \ref{s0}A). Using $\mu$ as a
parameter representing the magnitude of the perturbations, we
prove that (a) there exist infinitely many disjoint open intervals
of $\mu$, accumulating at $\mu = 0$, such that the entire
homoclinic tangle of the perturbed equation consists of one single
horseshoe of infinitely many symbols, (b) there are parameters in
between each of these parameter intervals, such that the
homoclinic tangle contains attracting periodic solutions, and (c)
there are also parameters in between where the homoclinic tangles
admit non-degenerate transversal homoclinic tangency of certain
dissipative hyperbolic periodic solutions. In particular, (c)
implies the existence of strange attractors with SRB measures for
a positive measure set of parameters.
\end{abstract}

\maketitle

\section{Introduction} \label{s0}

We start with an autonomous second-order ordinary differential
equation that contains a {\it non-resonant, dissipative} saddle
fixed point with a homoclinic solution. This autonomous equation
is then subjected to time periodic perturbations. In this paper we
study the dynamics of the homoclinic tangles formed by transversal
intersections of the stable and the unstable manifold of the
perturbed saddle point. These homoclinic tangles have been one of
the major inspirations for the dynamical systems theory and a long
standing puzzle in the studies of ordinary differential equations
in modern times.

\smallskip

\noindent {\bf A. Description of results.} \ Instead of focusing
on the picture of the globally induced time-T maps, by which H.
Poincar\'e observed an exceedingly complicated mess \cite{P} and
S. Smale constructed an embedded horseshoe map \cite{S}, we
compute the return maps induced by periodically perturbed
equations around the homoclinic solution in the extended phase
space. It has turned out that the return maps for these homoclinic
tangles are {\it infinitely wrapped horseshoe maps}, the geometric
structure of which is as follows. Take an annulus ${\mathcal A} =
S^1 \times I$. We represent points in $S^1$ and $I$ by using
variables $\theta$ and $z$ respectively. We call the direction of
$\theta$ the horizontal direction and the direction of $z$ the
vertical direction. To form an infinitely wrapped horseshoe map,
which we denote as ${\mathcal F}$, we first divide ${\mathcal A}$
into two vertical strips, which we denote as $V$ and $U$.
${\mathcal F}: V \to {\mathcal A}$ is defined on $V$ but not on
$U$. We compress $V$ in the vertical direction and stretch it in
the horizontal direction, making the image infinitely long towards
both ends. Then we fold it and wrap it around the annulus
${\mathcal A}$ infinitely many times. See Fig. 1.

We use a small parameter $\mu$ to represent the magnitude of the
time periodic perturbations. Denote the return maps obtained from
the periodically perturbed equations as ${\mathcal F}_{\mu}$ and
let\footnote{We caution that $V$ and $U$ depend also on $\mu$. So
to be completely rigorous we ought to write $V_{\mu}$ and
$U_{\mu}$ instead of $V$ and $U$. However, for ${\mathcal
F}_{\mu}$ derived from the periodically perturbed equations, $V$
and $U$ vary only slightly as $\mu$ varies, and we could
practically think them as being independent of $\mu$ at this
stage.}
$$
\Omega_{\mu} = \{ (\theta, z) \in V: \ {\mathcal
F}^n_{\mu}(\theta, z) \in V \ \ \ \forall n \geq 0 \}, \ \ \ \ \ \
\ \Lambda_{\mu} = \cap_{n\geq 0} {\mathcal
F}^n_{\mu}(\Omega_{\mu}).
$$
Then $\Omega_{\mu}$ represents all solutions that stay close to
the unperturbed homoclinic loop in forward times; $\Lambda_{\mu}$
is the set $\Omega_{\mu}$ is attracted to, representing all
solutions that stay close to the unperturbed homoclinic loop in
both the forward and the backward times. The geometrical and
dynamical structures of the homoclinic tangle represented by
${\mathcal F}_{\mu}$ are manifested in those of $\Omega_{\mu}$ and
$\Lambda_{\mu}$. $\Lambda_{\mu}$ obviously contain a horseshoe of
infinitely many symbols for all $\mu$. This horseshoe covers
Smale's horseshoe and all its variations. It is the one that
resides inside all homoclinic tangles.

\medskip

\begin{picture}(7, 5.5)
\put(3.7,0){ \psfig{figure=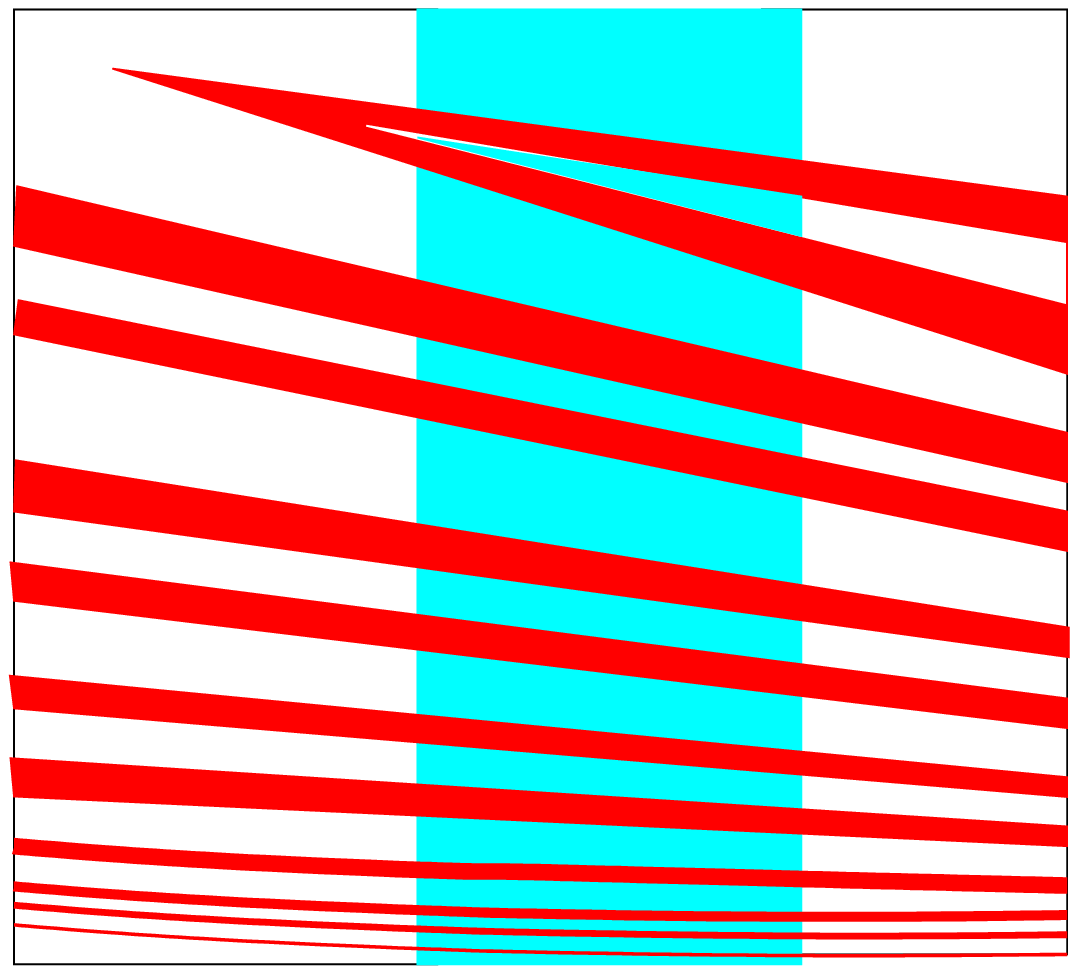,height = 5cm, width =
7cm} }
\end{picture}

\medskip

\centerline{Fig. 1 \ Infinitely wrapped horseshoe maps.}

\medskip

The structure of $\Omega_{\mu}$ and $\Lambda_{\mu}$ depend
sensitively on the location of the folded part of ${\mathcal
F}_{\mu}(V)$. If this part is deep inside of $U$, then the entire
homoclinic tangle is reduced to one horseshoe of infinitely many
symbols. If it is located inside of $V$, then the homoclinic
tangles are likely to have attracting periodic solutions or sinks
and observable chaos associated with non-degenerate transversal
homoclinic tangency. We prove that, as $\mu \to 0$, the folded
part of ${\mathcal F}_{\mu}(V)$ moves horizontally towards $\theta
= +\infty$ with a roughly constant speed with respect to $p = \ln
\mu$, crossing $V$ and $U$ infinitely many times along the way. It
then follows, under mild assumptions, that (a) there are
infinitely many disjoint open intervals of $\mu$, accumulating at
$\mu = 0$, such that the entire homoclinic tangle consists of one
single horseshoe of infinitely many symbols; (b) there are other
parameters in between these intervals, such that the homoclinic
tangle contains attracting periodic solutions; and (c) there are
also parameters in between where the homoclinic tangle admits
non-degenerate transversal homoclinic tangency. See Theorems
\ref{th1}-\ref{th3} in Sect. \ref{s3.2} for more details. In
particular, (c) implies the existence of strange attractors with
SRB measures for a positive measure set of parameters.

\noindent {\bf B. Method of study.} \ We use variables $(x, y)$ to
represent the phase space of the unperturbed equation and let $(x,
y) = (0, 0)$ be the saddle fixed point. Denote the homoclinic
solution for $(x, y) = (0, 0)$ as $\ell$. We construct a small
neighborhood of $\ell$ by taking the union of a small neighborhood
$U_{\varepsilon}$ of $(0, 0)$ and a small neighborhood $D$ around
$\ell$ out of $U_{\frac{1}{4} \varepsilon}$. See Fig. 2. Let
$\sigma^{\pm} \in U_{\varepsilon} \cap D$ be the two line segments
depicted in Fig. 2, both perpendicular to the homoclinic solution.
We use an angular variable $\theta \in S^1$ to represent the time.

\begin{picture}(7, 5.5)
\put(4,0){ \psfig{figure=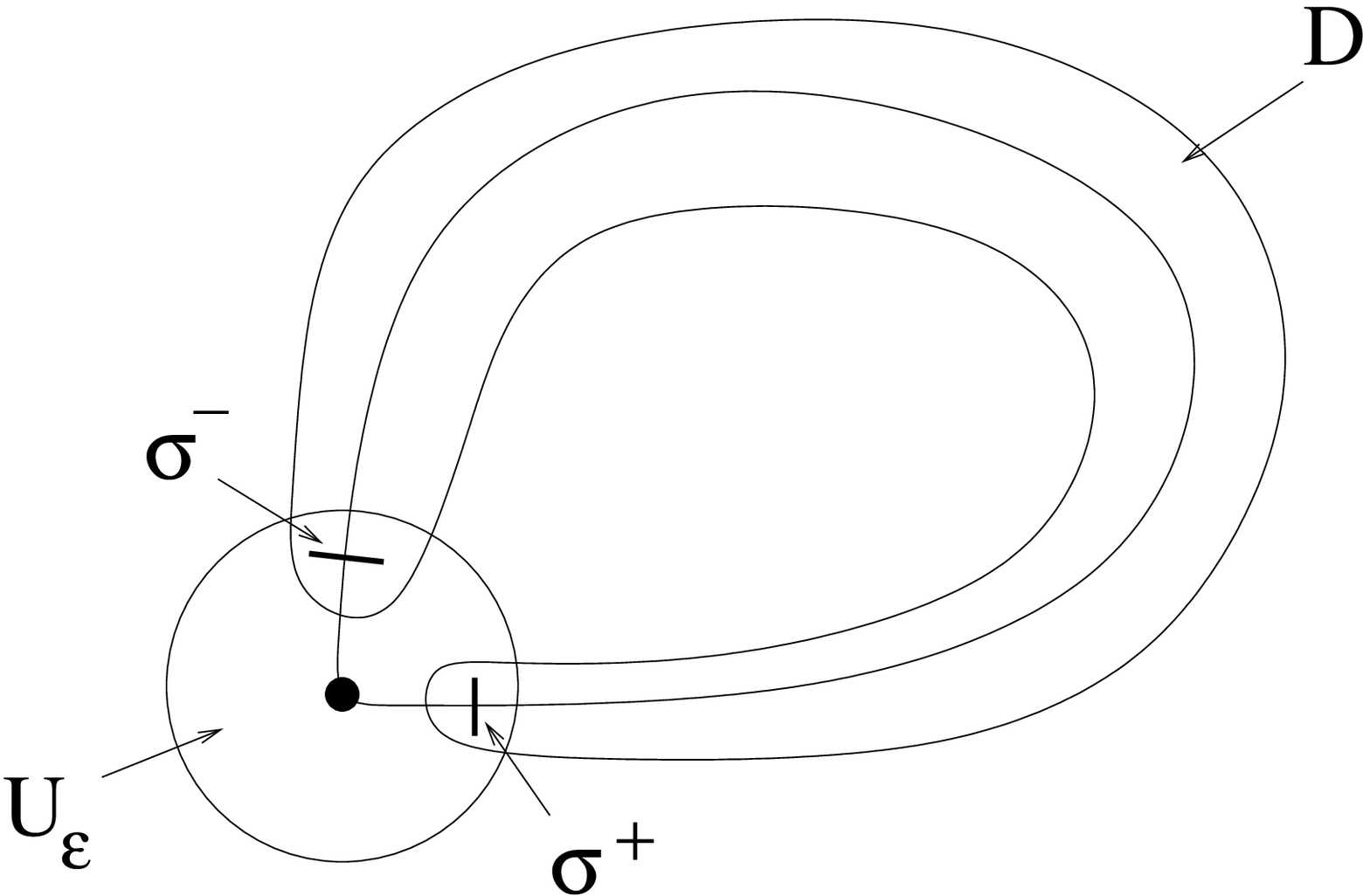,height = 5cm, width =
7cm} }
\end{picture}

\medskip

\centerline{Fig. 2 \ $U_{\varepsilon}$, $D$ and $\sigma^{\pm}$.}

\smallskip

In the extended phase space $(x, y, \theta)$ we denote
$$
{\mathcal U}_{\varepsilon} = U_{\varepsilon} \times S^1, \ \ \
{\bf D} = D \times S^1
$$
and let
$$
\Sigma^{\pm} = \sigma^{\pm} \times S^1.
$$

\begin{picture}(8, 6)
\put(3.5,0){ \psfig{figure=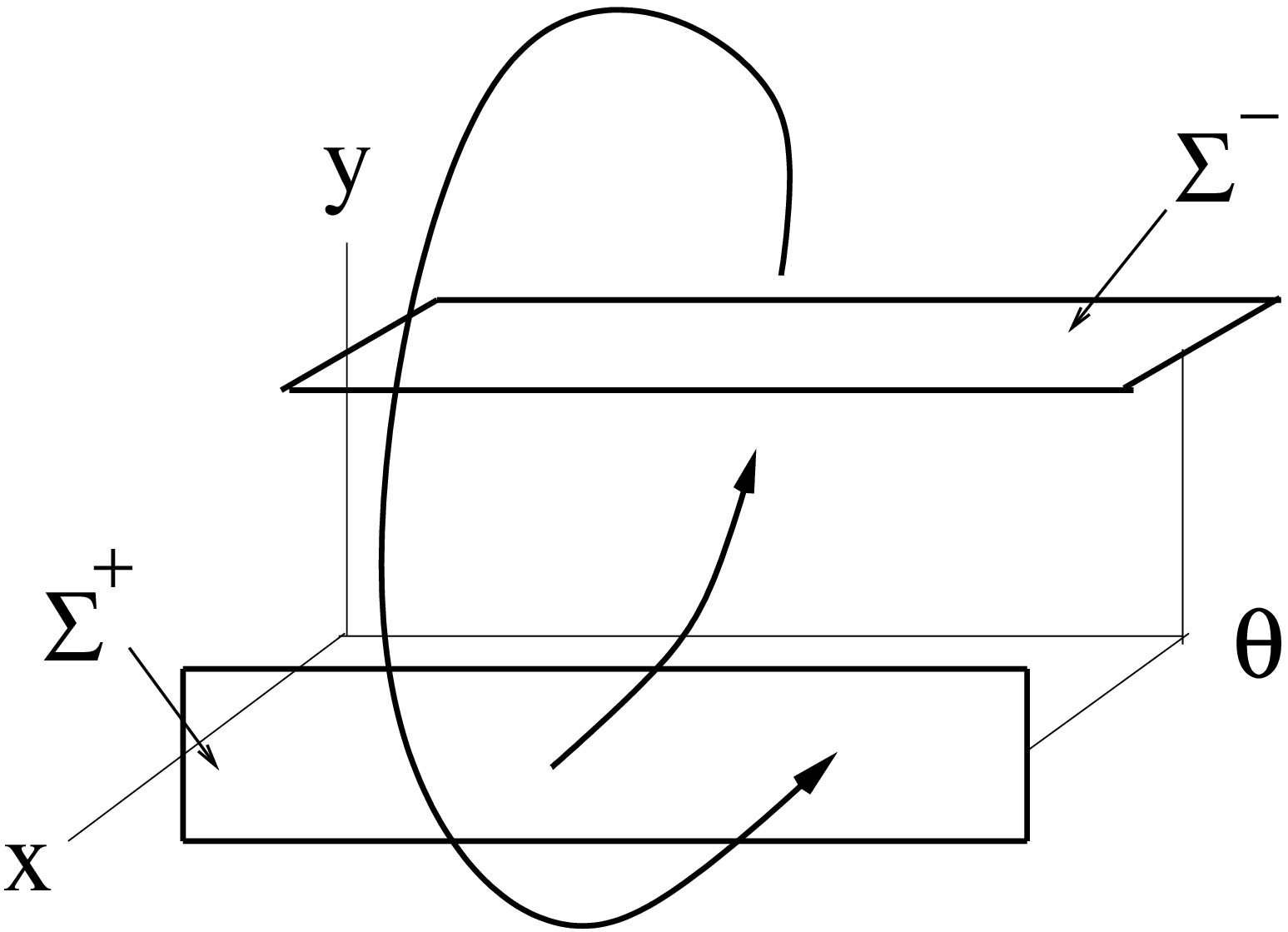,height = 5cm, width =
8cm} }
\end{picture}

\medskip

\centerline{Fig. 3 \ ${\mathcal N}$ and ${\mathcal M}$.}

\medskip

Let ${\mathcal N}: \Sigma^+ \to \Sigma^-$ be the maps induced by
the solutions on ${\mathcal U}_{\varepsilon}$ and ${\mathcal M}:
\Sigma^- \to \Sigma^+$ be the maps induced by the solutions on
${\bf D}$. See Fig. 3.  We first compute ${\mathcal M}$ and
${\mathcal N}$ separately, then compose ${\mathcal N}$ and
${\mathcal M}$ to obtain an explicit formula for the return map
${\mathcal N} \circ {\mathcal M}: \Sigma^- \to \Sigma^-$.

We follow the steps of \cite{WO} in deriving the return maps.
There are, however, two main differences between the classical
scenario of homoclinic tangles we now consider and the ones
studied in \cite{WO}. First, the return maps of this paper are
only partially defined on $\Sigma^-$. After following the entire
length of the homoclinic loop of the unperturbed equation, part of
$\Sigma^-$ (represented by $V$ in ${\mathcal A}$) would hit
$\Sigma^+$ on one side of the local stable manifold of the
perturbed saddle where they return to $\Sigma^-$; and the rest
(represented by $U$ in ${\mathcal A}$) would hit on the other side
where they sneak out. See Fig. 4. Second, analytic controls
represented by the $C^3$ estimates in \cite{WO} would deteriorate
as we approach to the transversal intersections of the stable and
the unstable manifold of the perturbed saddle, potentially
devastating the usefulness of the formulas obtained for the return
maps. Between the two, the second is essentially a technical issue
we need to overcome. The first is an intrinsic character of these
homoclinic tangles.

\begin{picture}(9, 5.5)
\put(4.5,0){ \psfig{figure=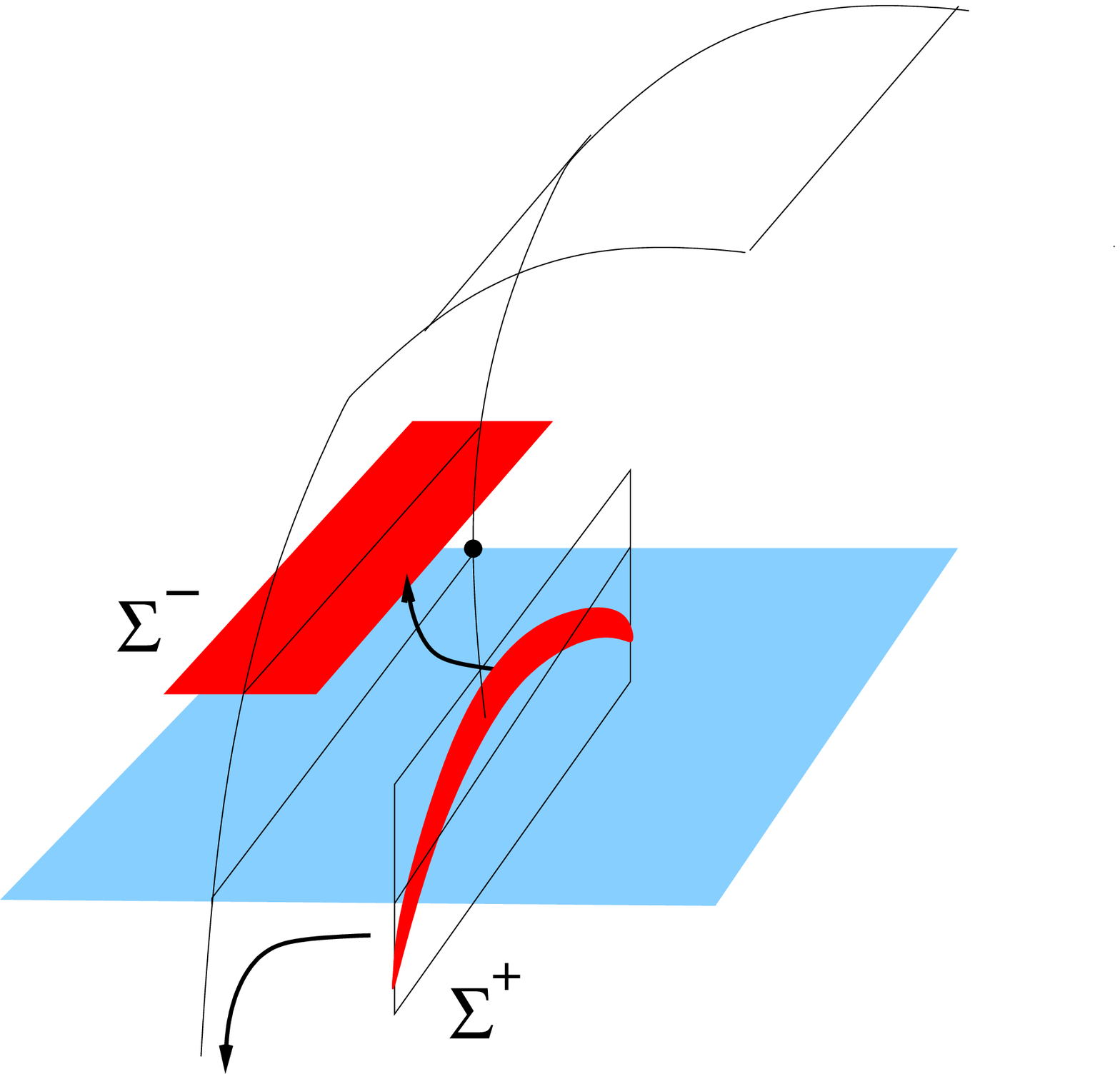,height = 5cm, width = 8cm}
}
\end{picture}

\medskip

\centerline{Fig. 4 \  Partial returns to $\Sigma^-$.}

\medskip

This paper is organized as follows. In Section \ref{s1} we
introduce the equations of study and a set of changes of variables
to transform the equations into certain canonical forms. In
Section \ref{s2} we compute the return maps. Proposition
\ref{prop1-s2.3c} in Section \ref{s2} is the main result of this
paper. Theorems about homoclinic tangles are then formulated and
proved in Section \ref{s3} assuming the forcing function is in the
form of $\sin \omega t$. Homoclinic tangles associated with
general forcing functions are studied in Section \ref{s4}.

\smallskip

\noindent {\bf C. Remarks on history.} \ Homoclinic tangles formed
by transversal intersections of the stable and the unstable
manifold of a periodically perturbed homoclinic saddle in systems
of ordinary differential equations were first observed by H.
Poincar\'e \cite{P} more than one hundred years ago. His
observation is regarded in general as the event that gave birth to
the modern theory of chaos and dynamical systems.

There exists a vast literature on periodically perturbed
differential equations (see for instance, the reference list of
\cite{GH}). We could put the related studies roughly into two
categories. The first category contains the ones that attempted to
understand the dynamics of the associated homoclinic tangles and
the second contains the ones that attempted to verify the
existence of these homoclinic tangles in concrete systems of
differential equations. Among the most influential in the first
category are the studies of Cartwright and Littlewood \cite{CL}
and Levinson \cite{L} on van der Pol's equation and the studies of
Sitnikov \cite{Sit} and Alekeseev [A] on Sitnikov's motions. These
studies led eventually to Smale's construction of his horseshoe
map \cite{S}. To the best of our knowledge, Smale's horseshoe is
essentially all that has been rigorously claimed for these
homoclinic tangles. In the second category, the most influential
is the development of Melnikov's method \cite{M}, purposed on
verifying the existence of homoclinic tangles in concrete systems
of differential equations.

Our idea of constructing return maps such as those derived in
\cite{WO} and in this paper is motivated by a work of Afraimovich
and Shilnikov, see also \cite{AH}. With the return maps obtained
in \cite{WO} and in this paper, we are able to apply many existing
theories on maps to the studies of homoclinic tangles of
differential equations. Among the theories directly applied are
the Newhouse theory \cite{N}, \cite{PT} on homoclinic tangency;
the theory of SRB measures \cite{Si}, \cite{R}, \cite{Bo}; the
theory of H\'enon-like attractors \cite{BC}, \cite{MV}, \cite{BY};
and the recent theory of rank one chaos \cite{WY1}-\cite{WY3}
based on the theory of Benedicks and Carleson on strongly
dissipative H\'enon maps \cite{BC}. These derived maps have led us
to many new results.

In Sect. \ref{s3.4} we present an overview on various dynamics
scenarios newly found around periodically perturbed homolcinic
solutions. We now know that, for the two main scenarios for the
time-T maps depicted in Fig. 5, the dynamics of the first is that
of an infinitely wrapped horseshoe maps and the second is that of
a rank one maps.

\begin{picture}(12, 4)
\put(2,0){ \psfig{figure=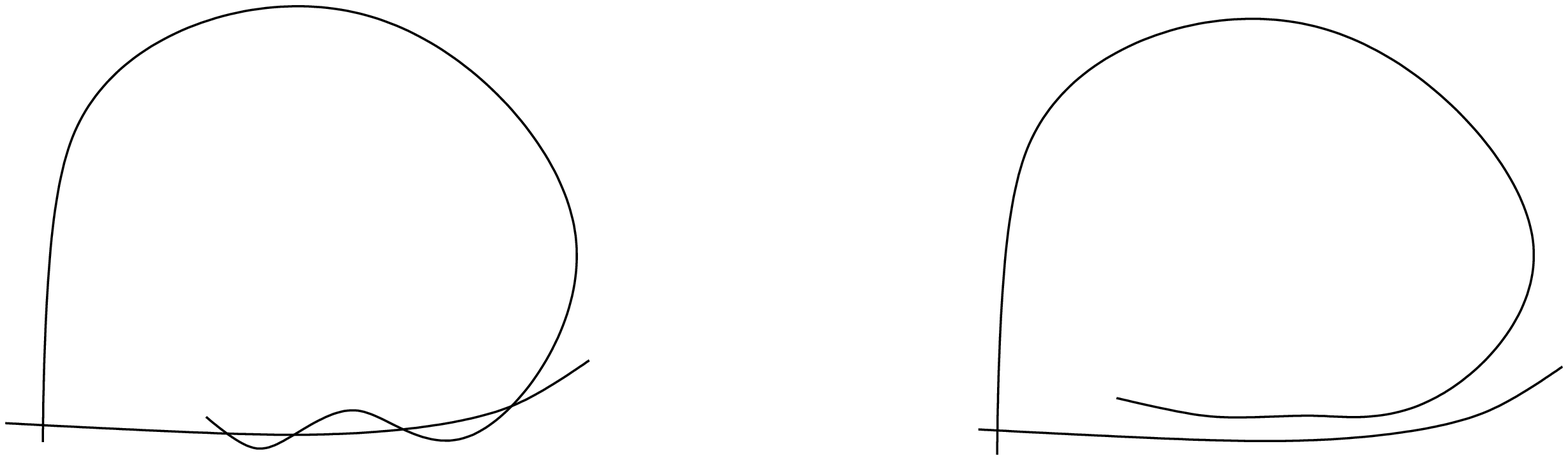,height = 3cm, width =
11cm} }
\end{picture}

\centerline{\ \ \ (a)  \hspace{2.2in} (b)}

\smallskip

\centerline{Fig. 5 \ The stable and the unstable manifold of a
perturbed saddle}

\medskip

\section{Equations and Canonical Forms}\label{s1}
In this section we first introduce the equations. We then
introduce coordinate changes to transform these equations into
canonical forms, which we will use in Section \ref{s2} to compute
the return maps.

\subsection{Equations of study} \label{s1.1} \ Let $(x, y) \in {\mathbb
R}^2$ be the phase variables and $t$ be the time. We start with an
autonomous system
\begin{equation}\label{f1-s1}
\frac{dx}{dt} = - \alpha x + f(x, y), \ \ \ \ \frac{d y}{dt} =
\beta y + g(x, y)
\end{equation}
where $f(x, y), g(x, y)$ are real-analytic functions defined on an
open domain ${\bf V} \subset {\mathbb R^2}$ satisfying $f(0, 0) =
g(0, 0) =\partial_x f(0, 0) = \partial_y f(0, 0) =
\partial_x g(0, 0) = \partial_y g(0, 0) = 0$. First we assume
that $(0, 0)$ is a {\it non-resonant, dissipative}  saddle point.
To be more precise we assume

\smallskip

{\bf (H1)} \ {\it (i) there exists $d_1, d_2> 0$ so that for all
$n, m \in {\mathbb Z}^+$,
$$
|n \alpha - m \beta| > d_1 (|n| + |m|)^{-d_2};
$$

\ \ \ \ \ \ \ \ \ (ii) $0 < \beta < \alpha$.}

\smallskip

(H1)(i) is a Diophantine non-resonance condition on $\alpha$ and
$\beta$. (H1)(ii) claims that the saddle point (0, 0) is
dissipative. Let us also assume that the positive $x$-side of the
local stable manifold and the positive $y$-side of the local
unstable manifold of $(0, 0)$ are included as part of a homoclinic
solution, which we denote as
$$
\ell = \{\ell(t)=(a(t), b(t)) \in {\bf V}, \  \ t \in {\mathbb
R}\}.
$$

Let ${\mathcal Q}(t): {\mathbb R} \to {\mathbb R}$ be a real
analytic function satisfying ${\mathcal Q}(t) = {\mathcal Q}(t+
2\pi)$. To the right of equation (\ref{f1-s1}) we add forcing
terms to form a non-autonomous system
\begin{equation}\label{f2-s1}
\begin{split}
\frac{dx}{dt} & = - \alpha x + f(x, y) + \mu A(x, y) (\rho + {\mathcal Q}(\omega t)), \\
\frac{dy}{dt} & = \beta y + g(x, y) + \mu B(x, y)(\rho + {\mathcal
Q}(\omega t))
\end{split}
\end{equation}
where $A(x, y), B(x, y)$ are also real analytic on ${\bf V}$
satisfying $A(0, 0) = B(0, 0) = 0$, $\partial_x A(0, 0) =
\partial_y A(0, 0) = \partial_x B(0, 0) = \partial_y B(0, 0) =0$.
We regard $\alpha, \beta$, $f(x, y), g(x, y)$, $A(x, y), B(x, y)$
and ${\mathcal Q}(t)$ as been fixed. $\omega, \rho, \mu$ are
forcing parameters. The range ${\mathbb P}$ for $(\omega, \rho,
\mu)$ is as follows: Let $R_{\omega}>>1$ be an arbitrarily picked
real number sufficiently large, $R_{\mu} >> R_{\rho} >>
R_{\omega}$. Define
\begin{equation}\label{f1a-s1.1}
{\mathbb P} = \{ (\omega, \rho, \mu) \in {\mathbb R}^3,  \ 0 <
\omega < R_{\omega}, \ R_{\rho}^{-1} < \rho < R_{\rho}, \ 0 < \mu
< R_{\mu}^{-1} \}.
\end{equation}
This study is exclusively on equation (\ref{f2-s1}) with
parameters inside of ${\mathbb P}$. From $R_{\mu}
>> R_{\rho} >> R_{\omega}$ we have, for all $(\omega, \rho, \mu)
\in {\mathbb P}$,
$$
\mu \rho, \ \ \mu \omega << 1.
$$
We can make the magnitude of $|\mu \rho|$ and $|\mu \omega|$ as
small as we desire by adjusting $R_{\mu}$.

We also fix the values of $\omega$ and $\rho$, leaving $\mu$ then
as the only parameter to vary. We study the solutions of equation
(\ref{f2-s1}) in the surroundings of the homoclinic loop $\ell$ in
the original phase space $(x, y)$, which we divide into a small
neighborhood $U_{\varepsilon}$ of $(0, 0)$ and a small
neighborhood $D$ around $\ell$ out of $U_{\frac{1}{4}
\varepsilon}$. See Fig. 2 in Section \ref{s0}. In Sect. \ref{s1.2}
we introduce a change of variables to linearize equation
(\ref{f2-s1}) on $U_{\varepsilon}$. In Sect. \ref{s1.3} we
introduce coordinate changes to transform (\ref{f2-s1}) on $D$
into a canonical form. In the rest of this paper $r>0$ is reserved
for an integer arbitrarily fixed, and we will control the
$C^r$-norm of the derived return maps in phase variables and
parameters.

\smallskip

\noindent {\bf Two small scales:} \  $\mu << \varepsilon<<1$
represent two small scales of different magnitude. $\varepsilon$
represents the size of a small neighborhood of $(x, y) = (0, 0)$,
where the linearizations of Sect. \ref{s1.2} is valid. Define
$$
U_{\varepsilon} = \{ (x, y): x^2 + y^2 < 4 \varepsilon^2 \} \ \ \
$$
and let $L^+, -L^-$ be the respective times at which the
homoclinic solution $\ell(t)$ enters $U_{\frac{1}{2}\varepsilon}$
in the positive and the negative directions. $L^+$ and $L^-$ are
related, both determined completely by $\varepsilon$ and
$\ell(t)$. The parameter $\mu << \varepsilon$ controls the
magnitude of the time-periodic perturbation. We make
$$
R_{\mu} >> \varepsilon^{-1} >> R_{\rho} >> R_{\omega}.
$$

\smallskip

\noindent {\bf Notation:} \ Quantities that are independent of the
phase variables, time and $\mu$ are regarded as constants and $K$
is used to denote a generic constant, the precise value of which
is allowed to change from line to line. On occasion, a specific
constant is used in different places. We use subscripts to denote
such constants as $K_0, K_1, \cdots$. We will also make
distinctions between constants depend on $\varepsilon$ and those
do not by making such dependencies explicit. A constant that
depends on $\varepsilon$ is written as $K(\varepsilon)$. A
constant written as $K$ is independent of $\varepsilon$.

\subsection{Linearization on $U_{\varepsilon}$}\label{s1.2}

From this point on all functions are regarded as functions in
phase variables, time $t$ and the parameter $\mu$. Let $X, Y$ be
such that
\begin{equation}\label{f1-s1.2}
\begin{split}
x & = X + P(X, Y) + \mu \tilde P(X, Y, \theta; \mu) \\
y & = Y + Q(X, Y) + \mu \tilde Q(X, Y, \theta; \mu)
\end{split}
\end{equation}
where $P, Q, \tilde P, \tilde Q$ as functions of $X$ and $Y$ are
real-analytic on $|(X, Y)| < 2 \varepsilon$, and the values of
these functions and their first derivatives with respect to $X$
and $Y$ at $(X, Y) = (0, 0)$ are all zero. As is explicitly
indicated in (\ref{f1-s1.2}), $P$ and $Q$ are independent of
$\theta$ and $\mu$. We also assume that
$$
\tilde P(X, Y, \theta+ 2 \pi; \mu)= \tilde P(X, Y, \theta; \mu), \
\ \tilde Q(X, Y, \theta+ 2 \pi; \mu)= \tilde Q(X, Y, \theta; \mu)
$$
are periodic of period $2 \pi$ in $\theta$ and they are also
real-analytic with respect to $\theta$ and $\mu$ for all $\theta
\in {\mathbb R}$ and $|\mu| < R_{\mu}^{-1}$. Substituting $\theta$
by $\omega t$ in (\ref{f1-s1.2}) defines a non-autonomous, near
identity coordinate transformation from $x, y$ to $(X, Y)$, which
we write explicitly as
\begin{equation}\label{f2-s1.2}
\begin{split}
x & = X + P(X, Y) + \mu \tilde P(X, Y, \omega t; \mu) \\
y & = Y + Q(X, Y) + \mu \tilde Q(X, Y, \omega t; \mu).
\end{split}
\end{equation}

We have
\begin{prop}\label{prop-2.1}
Assume that $\alpha$ and $\beta$ satisfy the Diophantine
non-resonance condition (H1)(i). Then there exists a small
neighborhood $U_{\varepsilon}$ of $(0, 0)$, the size of which are
completely determined by equation (\ref{f1-s1}) and $d_1, d_2$ in
(H1)(i), such that there exists an analytic coordinate
transformation in the form of (\ref{f2-s1.2}) that transforms
equation (\ref{f2-s1}) into
\begin{equation*}
\frac{d X}{dt} = - \alpha X, \ \ \ \ \frac{d Y}{dt} = \beta Y.
\end{equation*}
Moreover, the $C^r$-norms of $P, Q, \tilde P, \tilde Q$ as
functions of $X, Y, \theta, \mu$ are all uniformly bounded from
above by a constant $K$ that is independent of both $\varepsilon$
and $\mu$ on $(X, Y) \in U_{\varepsilon}$, $\theta \in {\mathbb
R}$ and $\mu \in (-R_{\mu}^{-1}, R_{\mu}^{-1})$.
\end{prop}
\noindent {\bf Proof:} \ This is a standard linearization result.
See for instance \cite{CLS} for a proof. \hfill $\square$

\subsection{A canonical form around homoclinic loop}\label{s1.3}
In this subsection we derive a standard form for equation
(\ref{f2-s1}) around the homoclinic loop of equation (\ref{f1-s1})
outside of ${\mathcal U}_{\frac{1}{4}\varepsilon}$. Let
$$
\ell = \{\ell(t)=(a(t), b(t)) \in {\mathbb R}^2, \  \ t \in
{\mathbb R}\}
$$
be the homoclinic solution of the unperturbed equation
(\ref{f1-s1}), and
$$
(u(t), v(t)) = \left|\frac{d}{dt} \ell(t)\right|^{-1} \frac{d}{dt}
\ell(t)
$$
be the unit tangent vector of $\ell$ at $\ell(t)$. Let us regard
$t$ in $\ell(t) =(a(t), b(t))$ not as time, but as a parameter
that parameterize the curve $\ell$ in $(x, y)$-space. We replace
$t$ by $s$ to write this homoclinic loop as $\ell(s) = (a(s),
b(s))$. We have
\begin{equation}\label{f1-s1.3}
\frac{d a(s)}{ds} =-\alpha a(s) +  f(a(s), b(s)), \ \ \ \ \frac{d
b(s)}{ds} = \beta b(s) + g(a(s), b(s)).
\end{equation}
By definition,
\begin{equation}\label{f2-s1.3}
\begin{split}
u(s) & = \frac{-\alpha a(s) + f(a(s), b(s))}{\sqrt{(-\alpha a(s) +
f(a(s), b(s)))^2 + (\beta b(s) + g(a(s),
b(s)))^2}}, \\
v(s) & = \frac{\beta b(s) + g(a(s), b(s))}{\sqrt{(-\alpha a(s) +
f(a(s), b(s)))^2 + (\beta b(s) +  g(a(s), b(s)))^2}}.
\end{split}
\end{equation}
Let
\begin{equation*}
{\bf e}(s) = (v(s), -u(s)).
\end{equation*}
We now introduce new variables $(s, z)$ such that
\begin{equation*}
(x, y) = \ell(s) + z {\bf e}(s).
\end{equation*}
This is to say that
\begin{equation}\label{f3-s1.3}
x = x(s, z) := a(s) + v(s) z, \ \  y =y(s, z):= b(s) - u(s) z.
\end{equation}

We derive the equations for (\ref{f2-s1}) in new variables $(s,
z)$ defined through (\ref{f3-s1.3}). Differentiating
(\ref{f3-s1.3}) we obtain
\begin{equation}\label{f4-s1.3}
\begin{split}
\frac{dx}{dt} & = (- \alpha a(s) + f(a(s), b(s)) + v'(s) z)
\frac{ds}{dt} + v(s)
\frac{dz}{dt} \\
\frac{dy}{dt} & = (\beta b(s) + g(a(s), b(s)) - u'(s) z)
\frac{ds}{dt} - u(s) \frac{dz}{dt}
\end{split}
\end{equation}
where $u'(s) = \frac{d u(s)}{ds}$, $v'(s) = \frac{d v(s)}{ds}$.
Let us denote
\begin{equation*}
\begin{split}
F(s, z) & = - \alpha (a(s) + z v(s)) + f(a(s) + z v(s), b(s) - z u(s)), \\
G(s, z) & = \beta (b(s) - z u(s)) + g(a(s) + z v(s), b(s) - z
u(s)), \\
{\mathbb A}(s, z) & = A(x(s, z), y(s, z)), \\
{\mathbb B}(s, z) & = B(x(s, z), y(s, z)).
\end{split}
\end{equation*}
By using equation (\ref{f2-s1}), we obtain from equation
(\ref{f4-s1.3}) the new equations for $s, z$ as
\begin{equation*}
\begin{split}
\frac{dz}{dt} & = v(s)F(s, z)- u(s) G(s, z) +\mu(v(s) {\mathbb
A}(s, z) - u(s) {\mathbb B}(s, z)) (\rho + {\mathcal Q}(\omega t))
\\ \frac{ds}{dt} &= \frac{v(s) G(s, z) + u(s) F(s, z) + \mu (v(s)
{\mathbb B}(s, z) + u(s) {\mathbb A}(s, z))(\rho + {\mathcal
Q}(\omega t))}{\sqrt{F(s, 0)^2 + G(s, 0)^2} + z(u(s) v'(s) - v(s)
u'(s))}.
\end{split}
\end{equation*}
We re-write these equations as
\begin{equation}\label{f5-s1.3}
\begin{split}
\frac{dz}{dt} & = E(s) z + z^2 w_2(s, z) + \mu(v(s) {\mathbb A}(s,
z) - u(s) {\mathbb B}(s, z))(\rho + {\mathcal Q}(\omega t)) \\
\frac{ds}{dt} & = 1 + z w_1(s, z, \omega t; \mu) + \frac{\mu (v(s)
{\mathbb B}(s, z) + u(s) {\mathbb A}(s, z))(\rho + {\mathcal
Q}(\omega t))}{\sqrt{F(s, 0)^2 + G(s, 0)^2}}
\end{split}
\end{equation}
where
\begin{equation}\label{f5a-s1.3}
\begin{split}
E(s) & = v^2(s) (-\alpha + \partial_x f(a(s), b(s))) + u^2(s)
(\beta + \partial_y g(a(s), b(s))) \\
& \ \ - u(s) v(s) (\partial_y f(a(s), b(s)) + \partial_x g(a(s),
b(s))).
\end{split}
\end{equation}
and $w_1(s, z, \theta + 2 \pi; \mu) = w_1(s, z, \theta; \mu)$ is
periodic in $\theta$ of period $2 \pi$.  Also in the rest of this
section we let $K_1(\varepsilon)$ be a given constant independent
of $\mu$, and regard equation (\ref{f5-s1.3}) as been defined on
$$
\{ s \in [-2 L^-, 2 L^+], \ \ \ |z| < K_1(\varepsilon) \mu; \ \ \
\mu \in (0, R_{\mu}^{-1}) \}.
$$
The $C^r$-norms of $w_1(s, z, \theta; \mu)$ and $w_2(s, z)$ are
bounded above by a constant $K(\varepsilon)$.

Finally we re-scale the variable $z$ by letting
\begin{equation}\label{f6-s1.3}
Z = \mu^{-1} z.
\end{equation}
We arrive at the following equations
\begin{equation}\label{f7-s1.3}
\begin{split}
\frac{dZ}{dt} & = E(s) Z + \mu \tilde w_2(s, Z, \omega t; \mu) +
{\mathbb H}(s) (\rho+ {\mathcal Q}(\omega t)) \\
\frac{ds}{dt} & = 1 + \mu \tilde w_1(s, Z, \omega t; \mu)
\end{split}
\end{equation}
where
\begin{equation}\label{f7a-s1.3}
{\mathbb H}(s) = v(s) A(a(s), b(s)) - u(s) B(a(s), b(s));
\end{equation}
and $(s, Z; \mu)$ are defined on
$$
{\mathbb D} = \{ (s, Z; \mu): \ s \in [-2L^-, 2L^+], \ |Z| \leq
K_1(\varepsilon), \ \mu \in (0, R_{\mu}^{-1}) \}.
$$
Note that here we assume that $R_{\mu}^{-1}$ is sufficiently small
so that
$$
\mu << \min_{s \in [-2L^-, 2L^+]} (F(s, 0)^2 + G(s, 0)^2).
$$
Again, the $C^r$ norms of $\tilde w_1, \tilde w_2$ with respect to
$s, Z; \mu$ and $t$ are uniformly bounded by a constant
$K(\varepsilon)$ for $(s, Z; \mu) \in {\mathbb D}$ and $t \in
{\mathbb R}$. Equation (\ref{f7-s1.3}) is the one we need.

\subsection{Technical estimates}\label{s1.4} Estimates presented in
this subsection are directly taken from [WO], which we include for
completeness.

\smallskip

\noindent {\bf Notation:} \ We are going to adopt the following
convention in comparing the magnitude of two functions $f(t)$ and
$g(t)$. We denote $f(t) \prec g(t)$ if there exists $K
> 0$ independent of $t$ so that $ |f(t)|< K |g(t)|$ as $t \to
\infty$ (or $-\infty$). We denote $f(t) \sim g(t)$ if in addition
we have $|f(t)| > K^{-1} |g(t)|$. We also denote $f(t) \approx
g(t)$ if
$$
\frac{f(t)}{g(t)} \to 1
$$
as $t \to \infty$ (or $-\infty$).

\smallskip

Recall that $\ell(t) = (a(t), b(t))$ is the homoclinic solution
for the hyperbolic fixed point $(0, 0)$ of equation (\ref{f1-s1}).
$(u(t), v(t))$ is the unit tangent vector of $\ell$ at $\ell(t)$
defined through (\ref{f2-s1.3}).
\begin{lemma} \label{lem1-s1.4}
As $t \to +\infty$,
\begin{equation*}
\begin{split}
& a(t) \sim e^{-\alpha t}, \ b(t) \prec e^{-2 \alpha t}, \ u(t)
\approx -1, \ v(t) \prec e^{-\alpha t}; \\
& a(-t) \prec e^{- 2\beta t}, \ b(-t) \sim e^{- \beta t}, \ u(-t)
\prec e^{-\beta t}, \ v(-t) \approx 1.
\end{split}
\end{equation*}
\end{lemma}
\noindent {\bf Proof:} We are simply re-stating the fact that
$\ell(t) \to (0, 0)$ with an exponential rate $-\alpha$ in the
positive time direction along the $x$-axis, and an exponential
rate $\beta$ in the negative time direction along the $y$-axis.
\hfill $\square$

\smallskip

Let $E(s)$ be as in (\ref{f5a-s1.3}).

\begin{lemma}\label{lem2-s1.4}
As $L^{\pm} \to +\infty$,

(i) $\int_{-L^-}^0 (E(s) + \alpha) ds \prec 1$, \ \ $\int_0^{L^+}
(E(s) - \beta) ds \prec 1$.

(ii) $\int_{-L^-}^0 E(s) ds \approx - \alpha L^-$, \ \
$\int_0^{L^+} E(s) ds \approx  \beta L^+$.
\end{lemma}
\noindent {\bf Proof:} \ (i) claims that the integrals are
convergent as $L^{\pm} \to \infty$. For the first integral, we
observe that by adding $\alpha$ to $E(t)$, we obtain $E(t) +
\alpha$ as a collection of terms, each of which decays
exponentially as $t \to -\infty$ according to Lemma
\ref{lem1-s1.4}. Similarly, taking $\beta$ away from $E(t)$, we
obtain $E(t) - \beta$ as a collection of terms, each of which
decays exponentially as $t \to \infty$.

For (ii) we write
\begin{eqnarray*}
\int_{-L^-}^0 E(s) ds & = & - \alpha L^{-} + \int_{-L^-}^0 (E(s) +
\alpha)
ds \\
\int_{0}^{L^+} E(s) ds & = & \beta L^+  + \int_{0}^{L^+} (E(s)-\beta) ds. \\
\end{eqnarray*}
(ii) now follows from (i). \hfill $\square$

\smallskip

We also have
\begin{lemma}\label{lem3-s1.4} As $\varepsilon \to 0$,
$\varepsilon \sim e^{-\alpha L^+} \sim e^{- \beta L^-}$.
\end{lemma}
\noindent {\bf Proof:} \ This follows directly from the definition
of $L^{\pm}$ and Lemma \ref{lem1-s1.4}. \hfill $\square$

\section{Derivation of Return Maps} \label{s2}
Let $\theta \in S^1$ be an angular variable for time. We re-write
equation (\ref{f2-s1}) as
\begin{equation}\label{f1-s2}
\begin{split}
\frac{dx}{dt} & = - \alpha x + f(x, y) + \mu A(x, y) (\rho + {\mathcal Q}(\theta)), \\
\frac{dy}{dt} & = \beta y + g(x, y) + \mu B(x, y)(\rho +
{\mathcal Q}(\theta)). \\
\frac{d \theta}{dt} & = \omega.
\end{split}
\end{equation}
$\Sigma^{\pm}$ are formally defined in Sect. \ref{s2.1} (See
Section \ref{s0}B). In Sect. \ref{s2.2}, we study coordinate
conversions between $(X, Y, \theta)$ and $(s, Z, \theta)$ on
$\Sigma^{\pm}$. ${\mathcal N}: \Sigma^+ \to \Sigma^-$, ${\mathcal
M}: \Sigma^- \to \Sigma^+$ and the return map ${\mathcal F} =
{\mathcal N} \circ {\mathcal M}: \Sigma^- \to \Sigma^-$ are
computed in Sect. \ref{s2.3}.

\subsection{Poincar\'e sections $\Sigma^{\pm}$}\label{s2.1}

We start with equation (\ref{f1-s2}) on ${\mathcal
U}_{\varepsilon}$. We have obtained in Sect. \ref{s1.2} a change
of variables on ${\mathcal U}_{\varepsilon}$ in the form of
(\ref{f1-s1.2}), that is,
\begin{equation}\label{f1-s2.1}
\begin{split}
x & = X + P(X, Y) + \mu \tilde P(X, Y, \theta; \mu) \\
y & = Y + Q(X, Y) + \mu \tilde Q(X, Y, \theta; \mu)
\end{split}
\end{equation}
that transforms equation (\ref{f1-s2}) to the linear equation
\begin{equation}\label{f2-s2.1}
\frac{d X}{dt} = -\alpha X, \ \ \ \ \frac{dY}{dt} = \beta Y, \ \ \
\ \frac{d \theta}{dt} = \omega
\end{equation}
on ${\mathcal U}_{\varepsilon}$. We define $\Sigma^{\pm}$ inside
of ${\mathcal U}_{\varepsilon} \cap {\bf D}$ by letting
$$
\Sigma^- = \{ (X, Y, \theta): \ Y = \varepsilon, \ |X| < \mu, \
\theta \in S^1 \},
$$
and
$$
\Sigma^+ = \{ (X, Y, \theta): \ X = \varepsilon, \ |Y| <
K_1(\varepsilon) \mu, \ \theta \in S^1 \}.
$$
$K_1(\varepsilon)$ will be precisely defined in Sect. \ref{s2.3}.
Observe that in [WO], $\Sigma^{\pm}$ are defined in slightly
different terms. The current definition is designed to avoid the
long and deteriorating derivative estimates of [WO].

\smallskip

We turn to the canonical form for equation (\ref{f1-s2}) on ${\bf
D}$. Let
\begin{equation}\label{f1-s2.1b}
x = a(s) + \mu v(s) Z, \ \  y = b(s) - \mu u(s) Z.
\end{equation}
Then according to Sect. \ref{s1.3}, equation (\ref{f1-s2}) on
${\bf D}$ is written in $(s, Z, \theta)$ as
\begin{equation}\label{f2-s2.1b}
\begin{split}
\frac{dZ}{dt} & = E(s) Z + \mu \tilde w_2(s, Z, \theta; \mu) +
{\mathbb H}(s) (\rho+ {\mathcal Q}(\theta)) \\
\frac{ds}{dt} & = 1 + \mu \tilde w_1(s, Z, \theta; \mu) \\
\frac{d \theta}{dt} & = \omega
\end{split}
\end{equation}
where
\begin{equation}\label{f3-s2.1b}
\begin{split}
E(s) & = v^2(s) (-\alpha + \partial_x f(a(s), b(s))) + u^2(s)
(\beta + \partial_y g(a(s), b(s))) \\
& \ \ - u(s) v(s) (\partial_y f(a(s), b(s)) + \partial_x g(a(s),
b(s)));
\end{split}
\end{equation}
\begin{equation}\label{f4-s2.1b}
{\mathbb H}(s) = v(s) A(a(s), b(s)) - u(s) B(a(s), b(s));
\end{equation}
and the $C^r$-norms of $\tilde w_1, \tilde w_2$ are bounded from
above by a $K(\varepsilon)$ on ${\bf D} \times (0, R_{\mu}^{-1})$
where
$$
{\bf D} = \{ (s, Z, \theta): \ s \in [-2L^-, 2L^+], \ |Z| \leq
K_1(\varepsilon), \ \theta \in S^1 \}.
$$

\medskip

Let $q \in \Sigma^+$ or $\Sigma^-$. We represent $q$ by using the
$(X, Y, \theta)$-coordinates, for which we have $X = \varepsilon$
on $\Sigma^+$ and $Y = \varepsilon$ on $\Sigma^-$. We can also use
$(s, Z, \theta)$-coordinate to represent the same $q$, for which
the defining equations for $\Sigma^{\pm}$ are not as direct. To
compute the return maps, we need to first attend two issues that
are technical in nature. First, we need to derive the defining
equations for $\Sigma^{\pm}$ for $(s, Z, \theta)$. Second, we need
to be able to change coordinates from $(X, Y, \theta)$ to $(s, Z,
\theta)$ and vice versa on $\Sigma^{\pm}$. We start with some
preparations in notation.

\smallskip

\noindent {\bf Notation:} \ The intended formula for the return
maps would inevitably contain terms that are explicit and terms
that are implicit. Implicit terms are usually ``error" terms, and
the usefulness of a derived formula would depend completely on how
well the error terms are controlled. In this paper we aim on
$C^r$-control on all error terms. The derivations of the return
maps involve a composition of maps and multiple coordinate
changes. To facilitate our presentation, from this point on we
adopt specific conventions for indicating controls on magnitude.
For a given constant, we write ${\mathcal O}(1)$, ${\mathcal
O}(\varepsilon)$ or ${\mathcal O}(\mu)$ to indicate that the
magnitude of the constant is bounded by $K$, $K \varepsilon$ or
$K(\varepsilon) \mu$, respectively. For a function of a set $V$ of
variables on a specific domain, we write ${\mathcal O}_V(1),
{\mathcal O}_V(\varepsilon)$ or ${\mathcal O}_V(\mu)$ to indicate
that the $C^r$-norm of the function on the specified domain is
bounded by $K, K \varepsilon$ or $K(\varepsilon) \mu$,
respectively. We chose to specify the domain in the surrounding
text rather than explicitly involving it in the notation. For
example, ${\mathcal O}_{Z, \theta}(\mu)$ represents a function of
$Z, \theta$, the $C^r$-norm of which is bounded above by
$K(\varepsilon) \mu$.

\smallskip

\noindent {\bf The new parameter $p$:} \ For the formulas obtained
to be most useful, it is also desirable that we have control on
the derivatives with respect to the forcing parameter $\mu$.
Taking derivative with respect $\mu$, however, are problematic
because such action takes $\mu$ out of the needed places. To
resolve this potentially damaging problem we introduce a new
parameter $p = \ln \mu$ and regard $p$, not $\mu$, as our
bottom-line parameter. In another word, we regard $\mu$ as a
shorthand for $e^p$, and all functions written in $\mu$ as
functions in $p$. Observe that $\mu \in (0, \mu_0)$ corresponds to
$p \in (-\infty, \ln \mu_0)$. This is a {\it very important
conceptual point} because by regarding a function $F(\mu)$ of
$\mu$ as a function of $p$, we have
$$
\partial_p F(\mu) = \mu \partial_{\mu} F(\mu).
$$
So regarding $F(\mu)$ as a function of $p$ would give back to us
that much needed factor $\mu$ in derivative estimates.

\subsection{Conversion of coordinates on
$\Sigma^{\pm}$}\label{s2.2}

We start with the defining equations for $\Sigma^+$ in $(s, Z,
\theta)$. Results for $\Sigma^-$ are similar.
\begin{lemma}\label{lem1-s2.2}
We have for $(s, Z, \theta) \in \Sigma^+$
$$
s = L^+ + {\mathcal O}_{Z, \theta, p}(\mu).
$$
\end{lemma}
\noindent {\bf Proof:} \ From (\ref{f1-s2.1}) and
(\ref{f1-s2.1b}), we have on $\Sigma^+$
\begin{equation}\label{f1-s2.2}
\begin{split}
a(s) + v(s) z & = \varepsilon + P(\varepsilon, Y) + \mu \tilde
P(\varepsilon, Y, \theta; \mu) \\
b(s) - u(s) z & = Y + Q(\varepsilon, Y) + \mu \tilde
Q(\varepsilon, Y, \theta; \mu).
\end{split}
\end{equation}
By definition
\begin{equation}\label{f2-s2.2}
\begin{split}
a(L^+) & = \varepsilon + P(\varepsilon, 0) \\
b(L^+) & = Q(\varepsilon, 0).
\end{split}
\end{equation}
Let
\begin{equation}\label{f3-s2.2}
\begin{split}
W_1 & = a(s) - a(L^+) + v(s) z - \mu \tilde P(\varepsilon, 0,
\theta; \mu), \\
W_2 & = b(s) - b(L^+) - u(s) z - \mu \tilde Q(\varepsilon, 0,
\theta; \mu).
\end{split}
\end{equation}
We have from (\ref{f1-s2.2}) and (\ref{f2-s2.2}),
\begin{equation*}
\begin{split}
W_1 & = P(\varepsilon, Y) - P(\varepsilon, 0) + \mu (\tilde
P(\varepsilon, Y, \theta; \mu) - \tilde P(\varepsilon, 0, \theta; \mu)) \\
W_2 & = Y + Q(\varepsilon, Y) - Q(\varepsilon, 0) + \mu ((\tilde
Q(\varepsilon, Y, \theta; \mu) - \tilde Q(\varepsilon, 0, \theta;
\mu))
\end{split}
\end{equation*}
which we re-write as
\begin{equation}\label{f4-s2.2}
\begin{split}
W_1 & = ({\mathcal O}(\varepsilon) + \mu {\mathcal O}_{\theta,
p}(1)) Y + {\mathcal O}_{Y, \theta, p}(1) Y^2\\
W_2 & = (1 + {\mathcal O}(\varepsilon) + \mu {\mathcal O}_{\theta,
p}(1))Y + {\mathcal O}_{Y, \theta, p}(1) Y^2.
\end{split}
\end{equation}
We first obtain
\begin{equation}\label{f5-s2.2}
Y = (1 + {\mathcal O}(\varepsilon) + \mu {\mathcal O}_{\theta,
p}(1))W_2 + {\mathcal O}_{W_2, \theta, p}(1) W_2^2
\end{equation}
by inverting the second line in (\ref{f4-s2.2}). We then
substitute into the first line in (\ref{f4-s2.2}) to obtain
\begin{equation*}
\begin{split}
W_1 & = ({\mathcal O}(\varepsilon) + \mu {\mathcal O}_{\theta,
p}(1))((1 + {\mathcal O}(\varepsilon) + \mu {\mathcal O}_{\theta,
p}(1)) W_2  + {\mathcal O}_{W_2, \theta, p}(1)
W_2^2) \\
& \ \ + {\mathcal O}_{Y, \theta, p}(1)  ((1 + {\mathcal
O}(\varepsilon) + \mu {\mathcal O}_{\theta, p}(1))W_2 +
{\mathcal O}_{W_2, \theta, p}(1) W_2^2)^2 \\
& = ({\mathcal O}(\varepsilon) + \mu {\mathcal O}_{\theta,
p}(1))W_2  + {\mathcal O}_{W_2, \theta, p}(1) W_2^2.
\end{split}
\end{equation*}
Consequently,
\begin{equation}\label{f6-s2.2}
F(s, Z, \theta, \mu) := W_1 - ({\mathcal O}(\varepsilon) + \mu
{\mathcal O}_{\theta, p}(1))W_2  + {\mathcal O}_{W_2, \theta,
p}(1) W_2^2 = 0,
\end{equation}
where $W_1, W_2$ as function of $s, Z, \theta$ and $\mu$ are
defined by (\ref{f3-s2.2}). To re-write $W_1, W_2$ we let
\begin{equation}
\xi = s - L^+
\end{equation}
and expand $a(s)$ in $\xi$ as
\begin{equation*}
a(s) = a(L^+) + a'(L^+) \xi + \sum_{i=2}^{\infty} a_i(L^+) \xi^i.
\end{equation*}
Expansions for $b(s), u(s)$ and $v(s)$ are similar. We have
\begin{equation}\label{f7-s2.2}
\begin{split}
W_1 & = a'(L^+) \xi + \sum_{i=2}^{\infty} a_i(L^+) \xi^i + v(L^+)
z + (v'(L^+)\xi + \sum_{i=2}^{\infty} v_i(L^+) \xi^i) z \\
& \ \ - \mu \tilde P(\varepsilon, 0, \theta; \mu)
\\
W_2 & = b'(L^+) \xi + \sum_{i=2}^{\infty} b_i(L^+) \xi^i - u(L^+)
z - (u'(L^+)\xi + \sum_{i=2}^{\infty} u_i(L^+) \xi^i)z \\
& \ \ -\mu \tilde Q(\varepsilon, 0, \theta; \mu).
\end{split}
\end{equation}
We now put (\ref{f7-s2.2}) for $W_1, W_2$ back into equation
(\ref{f6-s2.2}) and replace $z$ by $\mu Z$. We obtain
$$
(a'(L^+) - {\mathcal O}(\varepsilon) b'(L^+) + h(\theta, p, \xi)
\xi) \xi = {\mathcal O}_{Z, \theta, p}(\mu)
$$
where the $C^r$ norm of $h(\theta, p, \xi)$ is bounded from above
by $K(\varepsilon)$. From Lemma \ref{lem1-s1.4}, $a'(L^+) \approx
-\alpha \varepsilon, \ b'(L^+) = {\mathcal O}(\varepsilon^2)$. We
finally obtain
$$
s = L^+ + {\mathcal O}_{Z, \theta, p}(\mu)
$$
by solving $\xi$. This proves Lemma \ref{lem1-s2.2}. \hfill
$\square$

\smallskip

Lemma \ref{lem1-s2.2} is not precise enough. We need the following
refinement.
\begin{lemma}\label{lem2-s2.2}
We have on $\Sigma^+$,
$$
s-L^+ = - \frac{v(L^+) + {\mathcal O}(\varepsilon) u(L^+)}{a'(L^+)
- {\mathcal O}(\varepsilon) b'(L^+)} z + \frac{\mu}{a'(L^+) -
{\mathcal O}(\varepsilon) b'(L^+)} {\mathcal O}_{\theta, p}(1) +
{\mathcal O}_{Z, \theta, p}(\mu^2).
$$
\end{lemma}
\noindent {\bf Proof:} \ It suffices for us to drop all terms that
is ${\mathcal O}_{Z, \theta, p}(\mu^2)$ in equation
(\ref{f6-s2.2}) to solve for $\xi$. From Lemma \ref{lem1-s2.2} we
conclude that all terms in $\xi, z$ of degree higher than one are
${\mathcal O}_{Z, \theta, p}(\mu^2)$. With these terms all
dropped, (\ref{f6-s2.2}) becomes
\begin{equation}
(a'(L^+) - {\mathcal O}(\varepsilon) b'(L^+))\xi + (v(L^+) +
{\mathcal O}(\varepsilon) u(L^+)) z = \mu {\mathcal O}_{\theta,
p}(1),
\end{equation}
from which the estimates of Lemma \ref{lem2-s2.2} on $\Sigma^+$
follows. \hfill $\square$

\smallskip

From this point on we let
$$
{\mathbb X} = \mu^{-1} X, \ \ {\mathbb Y} = \mu^{-1} Y.
$$
\begin{lemma}\label{lem3-s2.2}
On $\Sigma^+$ we have
$$
{\mathbb Y} = (1 + {\mathcal O}(\varepsilon))Z + {\mathcal
O}_{\theta, p}(1) + {\mathcal O}_{Z, \theta, p}(\mu).
$$
\end{lemma}
\noindent {\bf Proof:} \ We have
\begin{equation}
\begin{split}
Y & = (1 + {\mathcal O}(\varepsilon))(b'(L^+) \xi - u(L^+) z - \mu
\tilde Q(\varepsilon, 0, \theta; \mu)) + {\mathcal O}_{Z, \theta,
p}(\mu^2) \\
& =  (1 + {\mathcal O}(\varepsilon))\left(-\left(u(L^+) + b'(L^+)
\frac{v(L^+) + {\mathcal O}(\varepsilon) u(L^+)}{a'(L^+) -
{\mathcal O}(\varepsilon) b'(L^+)}\right) z \right. \\
& \ \ \ \left. + \frac{\mu b'( L^+)}{a'(L^+) - {\mathcal
O}(\varepsilon) b'(L^+)} {\mathcal O}_{\theta, p}(1) - \mu \tilde
Q(\varepsilon, 0, \theta; \mu)\right) + {\mathcal O}_{Z, \theta,
p}(\mu^2) \\
& = (1 + {\mathcal O}(\varepsilon)) z + \mu {\mathcal O}_{\theta,
p}(1) + {\mathcal O}_{Z, \theta, p}(\mu^2),
\end{split}
\end{equation}
where the first equality follows from using (\ref{f5-s2.2}),
(\ref{f7-s2.2}) and Lemma \ref{lem1-s2.2}; the second equality
from using Lemma \ref{lem2-s2.2}. To obtain the third equality we
use $u(L^+) = -1 + {\mathcal O}(\varepsilon)$, $a'(L^+) \approx
-\alpha \varepsilon, \ b'(L^+) = {\mathcal O}(\varepsilon^2)$.
\hfill $\square$

\smallskip

Along similar lines we can also prove
\begin{lemma}\label{lem4-s2.2}
On $\Sigma^-$, we have

(i) $s = -L^- + {\mathcal O}_{Z, \theta, p}(\mu)$; and

(ii) $Z = (1 + {\mathcal O}(\varepsilon)) {\mathbb X} + {\mathcal
O}_{\theta, p}(1) + {\mathcal O}_{{\mathbb X}, \theta, p}(\mu)$.
\end{lemma}
\noindent {\bf Proof:} \ Left to the reader as an exercise. \hfill
$\square$

\subsection{The return map ${\mathcal F} = {\mathcal N} \circ
{\mathcal M}$}\label{s2.3}

First we compute ${\mathcal N}: \Sigma^+ \to \Sigma^-$ and
${\mathcal M}: \Sigma^- \to \Sigma^+$ separately. We then compose
${\mathcal N}$ and ${\mathcal M}$ by using Lemmas \ref{lem3-s2.2}
and \ref{lem4-s2.2}.

\medskip

\noindent {\bf A. The induced map ${\mathcal N}: \Sigma^+ \to
\Sigma^-$.} \ For $({\mathbb X}, {\mathbb Y}, \theta) \in
\Sigma^+$ we have ${\mathbb X} = \varepsilon \mu^{-1}$ by
definition. Similarly, for $({\mathbb X}, {\mathbb Y}, \theta) \in
\Sigma^-$ we have ${\mathbb Y}= \varepsilon \mu^{-1}$.  Denote a
point on $\Sigma^+$ by using $({\mathbb Y}, \theta)$ and a point
on $\Sigma^-$ by using $({\mathbb X}, \theta)$, and let
$$
({\mathbb X}_1, \theta_1) = {\mathcal N}({\mathbb Y}, \theta)
$$
for $({\mathbb Y}, \theta)\in \Sigma^+$.
\begin{prop}\label{prop1-s2.3a}
We have for $({\mathbb Y}, \theta) \in \Sigma^+$,
\begin{equation}\label{f1-s2.3a}
\begin{split}
{\mathbb X}_1 & =  (\mu \varepsilon^{-1})^{\frac{\alpha}{\beta}-1}
{\mathbb Y}^{\frac{\alpha}{\beta}} \\
\theta_1 & = \theta + \frac{\omega} {\beta} \ln (\varepsilon
\mu^{-1}) - \frac{\omega}{\beta} \ln {\mathbb Y}.
\end{split}
\end{equation}
\end{prop}
\noindent {\bf Proof:} \ Let $T$ be the time it takes for the
solution of (\ref{f2-s2.1}) from $(\varepsilon, Y, \theta) \in
\Sigma^+$ to get to $(X_1, \varepsilon, \theta_1) \in \Sigma^-$.
We have
\begin{equation*}
X_1 = \varepsilon e^{-\alpha T}, \ \ \ \ \varepsilon = Y e^{\beta
T}, \ \ \ \ \theta_1 = \theta + \omega T,
\end{equation*}
from which (\ref{f1-s2.3a}) follows. \hfill $\square$

\medskip

\noindent {\bf B. The induced map ${\mathcal M}: \Sigma^- \to
\Sigma^+$.} \ Let ${\mathbb H}(s)$ be as in (\ref{f4-s2.1b}). In
what follows, we write
\begin{equation}\label{f1a-s2.3b}
\begin{split}
A_L & = \int_{-L^-}^{L^+} {\mathbb H}(s) e^{- \int_{0}^{s} E(\tau)
d\tau} ds \\
\phi_L(\theta) & = \int_{-L^-}^{L^+} {\mathbb H}(s) {\mathcal
Q}(\theta + \omega s + \omega L^-) e^{- \int_{0}^{s} E(\tau)
d\tau} ds
\end{split}
\end{equation}
We also write
\begin{equation} \label{f2a-s2.3b}
P_L  = e^{\int_{-L^-}^{L^+} E(s) ds}, \ \ \ P_L^+ =
e^{\int_{0}^{L^+} E(s) ds}.
\end{equation}
Note that for $P_L$ we integrate from $s = -L^-$ to $s = L^+$,
while for $P_L^+$ the integration starts from $s = 0$. First we
have \
\begin{lemma}\label{lem1-s2.3b}
$$
P_L \sim \varepsilon^{\frac{\alpha}{\beta} - \frac{\beta}{\alpha}}
<< 1, \ \ \ \ \ \ \ P_L^+ \sim \varepsilon^{-\frac{\beta}{\alpha}}
>> 1.
$$
\end{lemma}
\noindent {\bf Proof:} \ Both estimates follows directly from
Lemmas \ref{lem2-s1.4} and \ref{lem3-s1.4}. \hfill $\square$

\smallskip

For $q = (s^-, Z, \theta) \in \Sigma^-$, the value of $s^-$ is
uniquely determined by that of $(Z, \theta)$ through Lemma
\ref{lem4-s2.2}(i). So it is allowed for us to use $(Z, \theta)$
to represent $q$. Let $(s(t), Z(t), \theta(t))$ be the solution of
equation (\ref{f2-s2.1b}) initiated at $(s^-, Z, \theta)$, and
$\tilde t$ be the time $(s(\tilde t), Z(\tilde t), \theta(\tilde
t))$ hit $\Sigma^+$. By definition ${\mathcal M}(q) = (s(\tilde
t), Z(\tilde t), \theta(\tilde t))$. In what follows we write
$$
s^+ = s(\tilde t), \ \ \hat Z = Z(\tilde t), \ \  \hat \theta =
\theta(\tilde t).
$$
\begin{prop}\label{prop1-s2.3b}
Denote $(\hat Z, \hat \theta) = {\mathcal M}(Z, \theta)$. We have
\begin{equation}\label{f1-s2.3b}
\begin{split}
\hat Z & =  P_L^+ (\rho A_L + \phi_L(\theta))+ P_L Z + {\mathcal
O}_{Z, \theta, p}(\mu) \\
\hat \theta & = \theta + \omega (L^+ + L^-) + {\mathcal O}_{Z,
\theta, p}(\mu) .
\end{split}
\end{equation}
\end{prop}
\noindent {\bf Proof:} \ Let us re-write equation (\ref{f2-s2.1b})
as
\begin{equation}\label{f2-s2.3b}
\begin{split}
\frac{dZ}{ds} & = E(s) Z  +
{\mathbb H}(s) (\rho+ {\mathcal Q}(\theta)) + {\mathcal O}_{s, Z, \theta, p}(\mu)\\
\frac{d \theta}{ds} & = \omega + {\mathcal O}_{s, Z, \theta,
p}(\mu)
\end{split}
\end{equation}
on ${\bf D} \times (0, R_{\mu}^{-1})$ where
$$
{\bf D} = \{(s, Z, \theta): \ \ s \in [-2 L^-, 2 L^+], \ |Z| <
K_1(\varepsilon), \ \theta \in S^1 \}.
$$
Dropping all error terms in (\ref{f2-s2.3b}) we have
\begin{equation}\label{f3-s2.3b}
\begin{split}
\frac{dZ}{ds} & = E(s) Z  +
{\mathbb H}(s) (\rho+ {\mathcal Q}(\theta))  \\
\frac{d \theta}{ds} & = \omega.
\end{split}
\end{equation}
We estimate the solution of equation (\ref{f2-s2.3b}) initiated at
$(Z, \theta)$ from $s = s^-$ to $s = s^+$ by the solution of
equation (\ref{f3-s2.3b}) initiated at the same $(Z, \theta)$ from
$s = -L^-$ to $s = L^+$. By the smooth dependencies of solutions
with respect to equations and initial conditions, the error of
such estimates, according to Lemma \ref{lem1-s2.2} and Lemma
\ref{lem4-s2.2}(i), is
$$
{\mathcal O}_{Z, \theta, p}(\mu) + {\mathcal O}_{\hat Z, \hat
\theta, p}(\mu)
$$
provided that both solutions stay inside of ${\bf D}$. By solving
(\ref{f3-s2.3b}), we obtain
\begin{equation}\label{f4-s2.3b}
\begin{split}
\hat Z & = P_L (Z + \Phi_L(\theta))+ {\mathcal O}_{Z, \theta,
p}(\mu) + {\mathcal O}_{\hat Z, \hat \theta, p}(\mu) \\
\hat \theta & = \theta + \omega (L^+ + L^-) + {\mathcal O}_{Z,
\theta, p}(\mu) + {\mathcal O}_{\hat Z, \hat \theta, p}(\mu)
\end{split}
\end{equation}
where $P_L$ is as in (\ref{f2a-s2.3b}) and
\begin{equation}\label{f5-s2.3b}
\Phi_L(\theta) = \int_{-L^-}^{L^+} {\mathbb H}(s)(\rho + {\mathcal
Q}(\theta + \omega L^- + \omega \tau))\cdot e^{-\int_{-L^-}^\tau
E(\hat \tau) d\hat \tau} d\tau.
\end{equation}
From (\ref{f4-s2.3b}) we have
\begin{equation}\label{f6-s2.3b}
\begin{split}
\hat Z & = P_L (Z + \Phi_L(\theta))+ {\mathcal O}_{Z, \theta,
p}(\mu) \\
\hat \theta & = \theta + \omega(L^+ + L^-) + {\mathcal O}_{Z,
\theta, p}(\mu).
\end{split}
\end{equation}

Let
\begin{equation}\label{f7-s2.3b}
K_1(\varepsilon) = \max_{\theta \in S^1, \ s \in [-2L^-, 2L^+],
\mu \in (-R_{\mu}^{-1}, R_{\mu}^{-1})} P_s (2 + |\Phi_s(\theta)|)
\end{equation}
where $P_s$ and $\Phi_s$ are obtained by replacing $L^+$ with $s$
in $P_L$ and $\Phi_L$. $K_1(\varepsilon)$ is the one we use for
${\bf D}$ and $\Sigma^+$. Solutions of (\ref{f2-s2.3b}) initiated
on $\Sigma^-$ will stay inside of ${\bf D}$ before hitting
$\Sigma^+$. To finish, we observe that
\begin{equation*}
\begin{split}
P_L \Phi_L(\theta)  & =P_L^+ \cdot \int_{-L^-}^{L^+} {\mathbb
H}(s) (\rho + {\mathcal Q}(\theta+ \omega L^- + \omega s)) e^{-
\int_{0}^{s} E(\tau)
d\tau} ds \\
& = P_L^+ (\rho A_L + \phi_L(\theta)).
\end{split}
\end{equation*}
This finishes the proof of Proposition \ref{prop1-s2.3b}. \hfill
$\square$

\medskip

\noindent {\bf C. The return map ${\mathcal F} = {\mathcal N}
\circ {\mathcal M}$} \ We are now ready to compute the return map
${\mathcal F} = {\mathcal N} \circ {\mathcal M}: \Sigma^- \to
\Sigma^-$. We use $({\mathbb X}, \theta)$ to represent a point on
$\Sigma^-$ and denote $(\tilde {\mathbb X}, \tilde \theta) =
{\mathcal F}({\mathbb X}, \theta)$.
\begin{prop}\label{prop1-s2.3c}
The map ${\mathcal F} ={\mathcal N} \circ {\mathcal M}: \Sigma^-
\to \Sigma^-$ is given by
\begin{equation}\label{f1-s2.3c}
\begin{split}
\tilde {\mathbb X} & =  (\mu
\varepsilon^{-1})^{\frac{\alpha}{\beta}-1} [(1 + {\mathcal
O}(\varepsilon))P_L^+ {\mathbb F}({\mathbb X},
\theta)]^{\frac{\alpha}{\beta}} \\
\tilde \theta & = \theta + \omega (L^+  + L^-) +
\frac{\omega}{\beta} \ln \mu^{-1} \varepsilon (1 + {\mathcal
O}(\varepsilon))P_L^+ - \frac{\omega}{\beta} \ln {\mathbb
F}({\mathbb X}, \theta) + {\mathcal O}_{{\mathbb X}, \theta,
p}(\mu)
\end{split}
\end{equation}
where
\begin{equation}\label{f2-s2.3c}
\begin{split}
{\mathbb F}({\mathbb X}, \theta) & = (\rho A_L + \phi_L(\theta))+
P_L (P_L^+)^{-1} (1+ {\mathcal O}(\varepsilon)){\mathbb X} \\
&  \ \ \ \ + (P_L^+)^{-1}(1+P_L) {\mathcal O}_{\theta, p}(1) +
{\mathcal O}_{{\mathbb X}, \theta, p}(\mu),
\end{split}
\end{equation}
and $P_L, P_L^+$ and $\phi_L(\theta)$ are as in (\ref{f1a-s2.3b})
and (\ref{f2a-s2.3b}).
\end{prop}
\noindent {\bf Proof:} \ By using Proposition \ref{prop1-s2.3b}
and Lemma \ref{lem4-s2.2}, we have
\begin{equation*}
\begin{split}
\hat Z & = P_L(1+ {\mathcal O}(\varepsilon)){\mathbb X} + P_L^+
(\rho A_L + \phi_L(\theta)) + P_L {\mathcal O}_{\theta, p}(1) +
{\mathcal O}_{{\mathbb X}, \theta, p}(\mu) \\
\hat \theta & = \theta + \omega (L^+ + L^-) + {\mathcal
O}_{{\mathbb X}, \theta, p}(\mu).
\end{split}
\end{equation*}

Let $\hat {\mathbb Y}$ be the ${\mathbb Y}$-coordinate for $(\hat
Z, \hat \theta)$, we have from Lemma \ref{lem3-s2.2},
$$
\hat {\mathbb Y} = (1 + {\mathcal O}(\varepsilon))P_L^+ {\mathbb
F}({\mathbb X}, \theta)
$$
where ${\mathbb F}({\mathbb X}, \theta)$ is as in
(\ref{f2-s2.3c}). We then obtain (\ref{f1-s2.3c}) by using
(\ref{f1-s2.3a}).  \hfill $\square$

\smallskip

We remark that ${\mathcal F} = {\mathcal N} \circ {\mathcal M}$ is
only defined on the part of $\Sigma^-$ where
$$
{\mathbb F}({\mathbb X}, \theta) > 0,
$$
and the set in $\Sigma^-$ defined by ${\mathbb F} = 0$ is on the
stable manifold of the saddle $(x, y) = (0, 0)$. Proposition
\ref{prop1-s2.3c} is the main result of this paper.

\section{Dynamics of Homoclinic Tangles: ${\mathcal Q}(t)= \sin t$} \label{s3}

In this section we let ${\mathcal Q}(t) = \sin t$ in equation
(\ref{f2-s1}). In Sect. \ref{s3.1} we derive the return maps. In
Sect. \ref{s3.2}, we prove that these return maps are {\it
infinitely wrapped horseshoe maps} (See Section \ref{s0}A). In
particular, we prove that (a) there exist infinitely many disjoint
open intervals of $\mu$, accumulating at $\mu = 0$, such that the
entire homoclinic tangle is one single horseshoe represented by a
full shift of infinitely many symbols (Theorem \ref{th1}); (b)
there are parameters in between each of these intervals, such that
the homoclinic tangle contains attracting periodic solutions
(Theorem \ref{th2}); and (c) there are also parameters in between
where the homoclinic tangle admits non-degenerate transversal
homoclinic tangency (Theorem \ref{th3}). The existence of
H\'enon-like attractors, following directly from Theorem \ref{th3}
and \cite{MV}, is stated in Corollary \ref{coro-add-s3}. In Sect.
\ref{s3.3} we study the associated homoclinic tangles by
numerically iterating the derived return maps. Finally in Sect.
\ref{s3.4}, we summarize various dynamics scenarios in the
surroundings of periodically perturbed homoclinic solutions newly
found through the return maps of Proposition \ref{prop1-s2.3c}.

\subsection{The return maps for homoclinic tangle}\label{s3.1}
Let $Q(t) = \sin t$ in equation (\ref{f2-s1}). Let ${\mathbb
F}({\mathbb X}, \theta)$ be as in Proposition \ref{prop1-s2.3c}.
The stable and the unstable manifold of $(x, y) = (0, 0)$ of
equation (\ref{f2-s1}) intersect if and only if there exists
$\theta$ such that ${\mathbb F}(0, \theta) = 0$. In [WO], the
authors excluded the possibility of these intersections by
restricting to a specific range of forcing parameters. We now
allow ${\mathbb F}(0, \theta) =0$.

Let
\begin{equation}\label{f1-s3.1}
\begin{split}
A & = \int_{-\infty}^{\infty}{\mathbb H}(s) e^{- \int_{0}^{s}
E(\tau) d\tau} ds  \\
C(\omega) & = \int_{-\infty}^{\infty} {\mathbb H}(s)
\cos (\omega s) e^{- \int_{0}^{s} E(\tau) d\tau} ds \\
S(\omega) & = \int_{-\infty}^{\infty}{\mathbb H}(s) \sin (\omega
s) e^{- \int_{0}^{s} E(\tau) d\tau} ds.
\end{split}
\end{equation}
Recall that $\ell(s) = (a(s), b(s)), s \in {\mathbb R}$ is the
homoclinic solution of equation (\ref{f1-s1}) and $(u(s), v(s))$
is the unit tangent vector of $\ell(s)$. Also recall that
${\mathbb H}(s)$ is as in (\ref{f7a-s1.3}) and $E(s)$ is as in
(\ref{f5a-s1.3}). Using the conclusions of Sect. \ref{s2.2}, it is
easy to verify that $A, C$ and $S$ are all well-defined. In the
rest of this section we assume that

\medskip

{\bf (H2)} \ (i) $A \neq 0$; and (ii) $C^2(\omega) + S^2(\omega)
\neq 0$.

\medskip

For a given equation (\ref{f1-s1}) satisfying (H1), (H2)(i) holds
for majority of $A(x, y)$ and $B(x, y)$. (H2)(ii) requires that,
as a function of $s$, the Fourier spectrum of the function
$$
R(s) = {\mathbb H}(s) e^{-\int_0^s E(\tau) d \tau}
$$
is not identically zero. We know that $R(s)$ decays exponentially
as a function of $s$, and it follows that the Fourier transform
$\hat R(\omega)$ is analytic in a strip contain the real
$\omega$-axis. Consequently, $\hat R(\omega) = 0$ for at most a
discrete set of values of $\omega$ because $R(s)$ is not
identically zero. Note that $\hat R(\omega) = C(\omega)+ i
S(\omega)$.

\smallskip

\noindent {\bf Specifications on parameters:} \ The parameters
$\omega, \rho, \varepsilon, \mu$ are specified as follows. First
we fix (arbitrarily) an $\omega$ such that (H2)(ii) holds. Then we
fix a value of $\rho$ such that\footnote{Let us assume $A>0$ here
to maintain a positive range for $\rho$. If $A<0$ we need to
switch $\rho$ to $-\rho$.}
$$
3 <\frac{1}{\rho A} \sqrt{C^2(\omega) + S^2(\omega)} < 9.
$$
Numbers $3$ and $9$ here have no specific meaning and can be
replaced by any other two numbers larger than 1. We let
$\varepsilon$ be small enough for a variety of reasons: one is to
validate the derivations of the previous sections and another is
to make
\begin{equation}\label{f-add-s3.1a}
2 <\frac{1}{\rho A_L} \sqrt{C^2_L(\omega) + S^2_L(\omega)} < 10
\end{equation}
where $A_L, C_L, S_L$ are obtained by replacing the integral
bounds $\pm \infty$ with $\pm L^{\pm}$ respectively in $A, C, S$.
$\mu$ ($ << \varepsilon$) is the only parameter we allow to vary.

\medskip

\noindent {\bf The return maps:} \ In the rest of this section we
use $z$ for ${\mathbb X}$, ${\mathcal A}$ for $\Sigma^-$. So we
write
$$
{\mathcal A} = \{ (\theta, z): \ \ \theta \in {\mathbb R}/(2 \pi
{\mathbb Z}), \ \ |z| < 1 \}.
$$
We regard $\omega, \rho, \varepsilon$ as been fixed. Let
$(\theta_1, z_1) = {\mathcal F}(\theta, z)$ for $(\theta, z) \in
{\mathcal A}$ where ${\mathcal F}$ is from Proposition
\ref{prop1-s2.3c}. We have
\begin{equation}\label{f1-s3.1a}
\begin{split}
\theta_1 & = \theta + {\bf a} - \frac{\omega}{\beta} \ln {\mathbb
F}(\theta, z, \mu) \\
z_1 & =  {\bf b} [{\mathbb F}(\theta, z,
\mu)]^{\frac{\alpha}{\beta}}
\end{split}
\end{equation}
where
\begin{equation} \label{f2-s3.1a}
\begin{split}
{\bf a} & = \frac{\omega}{\beta} \ln \mu^{-1} + \omega (L^+  + L^-) +
\frac{\omega}{\beta} \ln(\varepsilon (1 +
{\mathcal O}(\varepsilon))P_L^+ A_L \rho) \\
{\bf b} & = (\mu \varepsilon^{-1})^{\frac{\alpha}{\beta}-1} [(1 +
{\mathcal O}(\varepsilon))P_L^+ A_L \rho]^{\frac{\alpha}{\beta}}  \\
\end{split}
\end{equation}
and
\begin{equation} \label{f3-s3.1a}
{\mathbb F}(\theta, z, \mu) = 1 + {\bf c} \sin \theta + {\bf k} z
+ {\mathbb E}(\theta, \mu) + {\mathcal O}_{\theta, z, p}(\mu),
\end{equation}
in which
\begin{equation} \label{f4-s3.1a}
\begin{split}
{\bf c} & = (A_L \rho)^{-1} \sqrt{C_L^2 + S_L^2} \\
{\bf k} & = (A_L \rho)^{-1} P_L (P_L^+)^{-1} (1+ {\mathcal
O}(\varepsilon))
\end{split}
\end{equation}
and
\begin{equation} \label{f5-s3.1a}
{\mathbb E}(\theta, \mu)  = (A_L \rho)^{-1} (P_L^+)^{-1}(1+P_L)
{\mathcal O}_{\theta, p}(1).
\end{equation}
Note that in getting (\ref{f1-s3.1a}) we have changed $\theta +
\omega L^- + c_0$ to $\theta$ where $c_0$ is such that $\tan c_0 =
C_L^{-1} S_L$. ${\bf a}, {\bf b}, {\bf c}, {\bf k}$ and ${\mathbb
E}(\theta, \mu)$ are as follows:

\medskip

(i) ${\bf b} \to 0$ as $\mu \to 0$. We can think ${\mathcal F}$ as
an unfolding of the 1D maps
$$
f(\theta) = \theta + {\bf a} - \frac{\omega}{\beta} \ln (1 + {\bf
c} \sin \theta + {\mathbb E}(\theta, 0)).
$$

\smallskip

(ii) ${\bf a} \to +\infty$ as $\mu \to 0$. ${\bf a}$ is a large
number. But since it appears in the angular component we can
module it by $2 \pi$. With $\omega, \rho$ and $\varepsilon$ been
fixed, ${\bf a}$ is essentially $\omega \beta^{-1} \ln \mu^{-1}$.
Varying $\mu$ from a small $\mu_0>0$ to zero is to run ${\bf a}$
over $({\bf a}_0, +\infty)$ for some ${\bf a}_0 \sim \omega
\beta^{-1} \ln \mu_0^{-1}$.

\smallskip

(iii) By (\ref{f-add-s3.1a}), ${\bf c} \in [2, 10]$ is a constant
independent of $\mu$. Consequently, there exists an interval for
$\theta$ so that $1 + {\bf c} \sin \theta \leq 0$, and the stable
and the unstable manifold of the perturbed saddle of equation
(\ref{f2-s1}) do intersect. Also observe that from Lemma
\ref{lem1-s2.3b} we have
$$
{\mathbb E}(\theta, \mu) \sim \varepsilon^{\beta \alpha^{-1}}
{\mathcal O}_{\theta, p}(1).
$$
When $\varepsilon$ is sufficiently small, ${\mathbb E}(\theta,
\mu)$ is a $C^r$-small perturbation to $1 + {\bf c} \sin \theta$.

\smallskip

(iv) ${\bf k}$ is a small number independent of $\mu$. In fact, $
{\bf k} \sim \varepsilon^{\alpha \beta^{-1}}$ from Lemma
\ref{lem1-s2.3b}. {\bf k} is, however, much larger than $\mu$ and
it follows that the first derivative of ${\mathbb F}(\theta, z,
\mu)$ with respect to $z$ is $\approx {\bf k}$. This implies that
the unfolding from $f(\theta)$ in (i) to ${\mathcal F}$ is
non-degenerate in $z$-direction, and is controlled completely by
the linear term ${\bf k} z$.

\smallskip

(v) It is important that ${\mathbb E}(\theta, \mu)$ is {\it
independent} of $z$. Otherwise we would have trouble in
controlling what happens in $z$-direction. See (iv) above.

\smallskip

\noindent {\bf New notation on parameter:} \ In the rest of this
section we put ${\bf a}$ in the place of $p$, regarding it as the
bottom line parameter. Both $\mu$ and $p$ are regarded as
functions of ${\bf a}$. Since we have fixed $\omega, \rho$ and
$\varepsilon$, ${\bf c}$ and ${\bf k}$ are fixed constants
independent of ${\bf a}$. We denote the return maps as ${\mathcal
F}_{\bf a}$ to emphasize that ${\bf a}$ is the parameter. ${\bf
b}$ is a function of ${\bf a}$. Because ${\bf a}$ and $p$ are
linearly related, we have ${\mathcal O}_{\theta, p}(1) = {\mathcal
O}_{\theta, {\bf a}}(1)$ in (\ref{f5-s3.1a}) and ${\mathcal
O}_{\theta, z, p}(\mu) = {\mathcal O}_{\theta, z, {\bf a}}(\mu)$
in (\ref{f3-s3.1a}).

\subsection{Homoclinic tangles as an infinitely wrapped
horseshoe map}\label{s3.2} For $q = (\theta, z) \in {\mathcal A}$,
let ${\bf v} = (u, v)$ be a tangent vector of ${\mathcal A}$ at
$q$ and let $s({\bf v}) = v u^{-1}$. $s({\bf v})$ is the slope of
${\bf v}$. We say that ${\bf v}$ is {\it horizontal} if $|s({\bf
v})| < \frac{1}{100}$ and ${\bf v}$ is {\it vertical} if $|s({\bf
v})| > 100$. A curve in ${\mathcal A}$ is a {\it horizontal curve}
if all its tangent vectors are horizontal and it is a {\it
vertical curve} if all its tangent vectors are vertical. A
vertical curve is {\it fully extended} if it reaches both
boundaries of ${\mathcal A}$ in $z$-direction. A region in
${\mathcal A}$ that is bounded by two non-intersecting, fully
extended vertical curves is a {\it vertical strip}. For a given
vertical strip $V$, {\it a horizontal strip} in $V$ is a region
bounded by two non-intersecting horizontal curves traversing $V$
in $\theta$-direction.

Observe that
\begin{equation}\label{f1-s3.1b}
{\mathbb F}(\theta, z, \mu) = {\bf k} z + 1 + {\bf c} \sin \theta
+ {\mathbb E}(\theta, \mu) + {\mathcal O}_{\theta, z, {\bf
a}}(\mu) =0
\end{equation}
defines two fully extended vertical curves that divide ${\mathcal
A}$ into two vertical strips, which we denote as $V$ and $U$. Let
${\mathbb F}
> 0$ on $V$ and ${\mathbb F} < 0$ on U. ${\mathcal F}_{\bf a}$ is
well-defined on $V$ but not on $U$. $U$ is the window through
which the solutions of equation (\ref{f2-s1}) sneak out.

Let
\begin{equation}\label{f2-s3.1b}
\Omega_{\bf a}  = \{ (\theta, z) \in V: \ {\mathcal F}^n_{\bf
a}(\theta, z) \in V, \ \forall n \geq 0 \}, \ \ \ \ \ \ \ \
\Lambda_{\bf a} = \cap_{n\geq 0} {\mathcal F}^n_{\bf
a}(\Omega_{\bf a}).
\end{equation}
$\Omega_{\bf a}$ represents all solutions of equation
(\ref{f2-s1}) that stay close to $\ell$ in forward times;
$\Lambda_{\bf a}$ is the set $\Omega_{\bf a}$ is attracted to,
representing all solutions that stay close to $\ell$ in both the
forward and the backward times. $\Omega_{\bf a}$ and $\Lambda_{\bf
a}$ together represent the homoclinic tangles, the structure of
which we now unravel through ${\mathcal F}_{\bf a}$.

For a fixed $z \in [-1, 1]$, let
$$
I_z = \{\theta \in (-\frac{1}{2} \pi, \frac{3}{2} \pi]: \ (\theta,
z) \in V \}.
$$
$I_z$ is an interval in $(-\frac{1}{2} \pi, \frac{3}{2} \pi)$,
which we denote as $(\theta_l(z), \theta_r(z))$. Let $h_z = \{
(\theta, z): \theta \in I_z \}$. ${\mathcal F}_{\bf a}(h_z)$ is a
1D curve in ${\mathcal A}$ parameterized in $\theta$, which we
denote as $(z_1(\theta), \theta_1(\theta))$. By definition
\begin{equation}\label{f3-s3.1b}
\theta_1(\theta)  = \theta + {\bf a} - \frac{\omega}{\beta} \ln
{\mathbb F}(\theta, z, \mu).
\end{equation}
We have
\begin{lemma}\label{lem1-s3.2} Assume that $\omega \beta^{-1} > 100$.

(a) $\lim_{\theta \to \theta_r(z)^-} (\theta_1, z_1)  =
\lim_{\theta \to \theta_l(z)^+} (\theta_1, z_1) = (+\infty, 0)$.

(b) For every fixed $z \in [-1, 1]$, there exists a unique value
of $\theta$, which we denote as $\theta_c(z)$, such that
$$
\frac{d \theta_1}{d \theta}(\theta_c(z)) = 0.
$$

(c) Let
\begin{equation}\label{f4-s3.1b}
V_f = \cup_{z \in [-1, 1]} \{ (\theta, z) \in V: \  \left| \frac{d
\theta_1}{d \theta} \right|< 2\}.
\end{equation}
Then $V_f$ is a vertical strip, the horizontal size of which is $<
10 \omega^{-1} \beta$.
\end{lemma}
\noindent {\bf Proof:} \ Observe, from (\ref{f1-s3.1b}), that
$\theta_l \in (-\frac{1}{2} \pi, 0)$ where $\cos \theta > 0$, and
$\theta_r \in (\pi, \frac{3}{2} \pi)$ where $\cos \theta < 0$. (a)
follows directly from the fact that, as $\theta \to \theta_l^+,
\theta_r^-$, ${\mathbb F} \to 0$. To prove (b) we first observe
that, because ${\mathbb F} \to 0$ as $ \theta \to \theta_r^-$,
$$
|1 + {\bf c} \sin \theta_r^-| < K \varepsilon^{\beta \alpha^{-1}}
<< 1,
$$
and it follows that
$$
\frac{\partial {\mathbb F}}{\partial \theta}(\theta_r(z)^-, z)
\approx {\bf c} \cos \theta_r < -1.
$$
Consequently,
$$
\lim_{\theta \to \theta_r^-} \frac{d \theta_1}{d \theta} =
\lim_{\theta \to \theta_r^-} \left(1 - \omega \beta^{-1}
\frac{1}{\mathbb F} \frac{\partial {\mathbb F}}{\partial
\theta}\right) = +\infty.
$$
Similarly, we have
$$
\lim_{\theta \to \theta_l^+} \frac{d \theta_1}{d \theta} =
\lim_{\theta \to \theta_l^+} \left(1 - \omega \beta^{-1}
\frac{1}{\mathbb F} \frac{\partial {\mathbb F}}{\partial \theta}
\right)= - \infty.
$$
Therefore there exists at least one $\theta_c(z)$ satisfying
$\frac{d \theta_1}{d \theta} = 0$. For the uniqueness we observe
that
$$
\frac{d^2 \theta_1}{d \theta^2} = - \frac{\omega
\beta^{-1}}{{\mathbb F}^2}\left(\frac{\partial^2 {\mathbb
F}}{\partial \theta^2}{\mathbb F} - \left(\frac{\partial {\mathbb
F}}{\partial \theta} \right)^2 \right) \approx \frac{\omega
\beta^{-1}}{{\mathbb F}^2}({\bf c}^2 + {\bf c} \sin \theta) > 0
$$
for all $\theta$. Recall that ${\bf c} > 2$.

To prove (c) we observe that the boundary of $V_f$ is defined by
$$
\left|1 - \omega \beta^{-1} \frac{1}{\mathbb F} \frac{\partial
{\mathbb F}}{\partial \theta} \right| = 2,
$$
from which we obtain
$$
|\cos \theta| \leq \frac{9}{2} \omega^{-1} \beta + K
\varepsilon^{\beta \alpha^{-1}}.
$$
(c) follows directly from this estimate. \hfill $\square$

\smallskip

We are now ready to formally state and prove the first of our
theorems.

\begin{theorem}[{\bf Horseshoe of infinitely many symbols}]\label{th1}
Let ${\mathcal Q}(t) = \sin t$ and assume (H1) and (H2) for
equation (\ref{f2-s1}). Let the parameters $\omega, \rho,
\varepsilon$ been specified as in Sect. \ref{s3.1}. If in addition
$\omega \beta^{-1} > 100$, then there exists a sequence of $\mu$,
which we denote as
$$
1 >> \mu_1^{(r)} > \mu_1^{(l)} > \cdots > \mu_n^{(r)} >
\mu_n^{(l)} > \cdots > 0
$$
such that for all $\mu \in [\mu^{(l)}_n, \mu^{(r)}_n]$, ${\mathcal
F}_{\bf a}$ on
$$
\Lambda = \{ (\theta, z) \in V: \ \ {\mathcal F}^i_{\bf a}
(\theta, z) \in V, \ \forall i \in {\mathbb Z} \}
$$
conjugates to a full shift of countably many symbols.
\end{theorem}
\noindent {\bf Proof:} \ For different values of $\mu$, the
corresponding vertical curves in ${\mathcal A}$ defined by
(\ref{f1-s3.1b}) are ${\mathcal O}(\mu)$ close. So $V$ and $U$ are
almost stationary as ${\bf a}$ varies from ${\bf a}_0$ to
$+\infty$. On the other hand, it follows from (\ref{f1-s3.1a})
that, by varying ${\bf a}$ from ${\bf a}_0$ to $+\infty$, we move
${\mathcal F}_{\bf a}(V)$ horizontally towards $\theta = +\infty$.
Denote ${\mathcal F} = {\mathcal F}_{\bf a}$ and let $V_f$ be the
vertical strip defined through (\ref{f4-s3.1b}). The horizontal
size of ${\mathcal F}(V_f)$ is smaller than $20 \beta \omega^{-1}$
from Lemma \ref{lem1-s3.2} assuming $\omega \beta^{-1}
> 100$, which is in turn smaller than the horizontal size of $U$.
Therefore ${\mathcal F}(V_f)$ traverses ${\mathcal A}$ infinitely
many times in horizontal direction as we vary ${\bf a}$ from ${\bf
a}_0$ to $+\infty$ and  there are infinitely many sub-intervals of
${\bf a}$, such that ${\mathcal F}(V_f) \subset U$. For these
parameter values ${\mathcal F}(V) \cap V$ consists of countably
many horizontal strips in $V$ (see Fig. 1 in Section \ref{s0}A),
to each of which we assign a positive integer according naturally
to the order in which these strips are stacked in the downward
$z$-direction.

For $q \in {\mathcal A}$, let ${\bf v}$ be a tangent vector at
$q$. Let ${\mathcal C}_h(q)$ be the collection of all ${\bf v}$
satisfying $|s({\bf v})|< \frac{1}{100}$, and ${\mathcal C}_v(q)$
be the collection of all ${\bf v}$ satisfying $|s({\bf v})| >
100$. To prove that $\Lambda$ conjugates to a full shift of all
positive integers, it suffices to verify that we have, assuming
${\mathcal F}(V_f) \subset U$,
\begin{itemize}
\item[(i)] $D{\mathcal F}({\mathcal C}_h(q)) \subset {\mathcal
C}_h({\mathcal F}(q))$ on ${\mathcal F}^{-1}({\mathcal F}(V) \cap
V)$, and \item[(ii)] $D{\mathcal F}^{-1}({\mathcal C}_v(q))
\subset {\mathcal C}_v({\mathcal F}(q))$ on ${\mathcal F}(V)\cap
V$.
\end{itemize}

To prove (i) we first compute $D{\mathcal F}$ by using
(\ref{f1-s3.1a}). Let $(\theta_1, z_1) = {\mathcal F}(\theta, z)$,
we have
\begin{equation}\label{f1-add-s3.2}
D{\mathcal F} = \left( \begin{array}{cc} \frac{\partial
\theta_1}{\partial \theta} & \frac{\partial \theta_1}{\partial z}
\\
\frac{\partial z_1}{\partial \theta} & \frac{\partial
z_1}{\partial z}
\end{array} \right) = \left(
\begin{array}{cc} 1
- \omega \beta^{-1} \frac{1}{{\mathbb F}} \frac{\partial {\mathbb
F}}{\partial \theta} &  \omega \beta^{-1} \frac{1}{{\mathbb F}}
\frac{\partial {\mathbb F}}{\partial z} \\
\alpha \beta^{-1}  {\bf b} {\mathbb F}^{\alpha \beta^{-1}-1}
\frac{\partial {\mathbb F}}{\partial \theta} & \alpha \beta^{-1}
{\bf b} {\mathbb F}^{\alpha \beta^{-1} -1} \frac{\partial {\mathbb
F}}{\partial z}
\end{array} \right)
\end{equation}
where ${\mathbb F} = {\mathbb F}(\theta, z, \mu)$ is as in
(\ref{f3-s3.1a}) and
\begin{equation*}
\begin{split}
\frac{\partial {\mathbb F}}{\partial \theta} & = {\bf c} \cos
\theta + \varepsilon^{\beta \alpha^{-1}} {\mathcal O}_{\theta,
{\bf a}}(1)
+ {\mathcal O}_{\theta, z, {\bf a}}(\mu)\\
\frac{\partial {\mathbb F}}{\partial z} & = {\bf k} + {\mathcal
O}_{\theta, z, {\bf a}}(\mu).
\end{split}
\end{equation*}
Let ${\bf v}$ be such that $|s({\bf v})| < \frac{1}{100}$, we have
from (\ref{f1-add-s3.2})
\begin{equation}\label{f2-add-s3.2}
|s(D{\mathcal F}({\bf v}))| = \left|\frac{\alpha \beta^{-1}  {\bf
b} {\mathbb F}^{\alpha \beta^{-1}-1} \frac{\partial {\mathbb
F}}{\partial \theta}  +\alpha \beta^{-1} {\bf b} {\mathbb
F}^{\alpha \beta^{-1}-1} \frac{\partial {\mathbb F}}{\partial z}
s({\bf v}) }{(1- \omega \beta^{-1} \frac{1}{{\mathbb F}}
\frac{\partial {\mathbb F}}{\partial \theta}) +  \omega \beta^{-1}
\frac{1}{{\mathbb F}} \frac{\partial {\mathbb F}}{\partial z}
s({\bf v})}\right|.
\end{equation}
We have two cases to consider.

\smallskip

{\it Case 1: ${\mathbb F} \geq \sqrt{\bf k}$.} \ In this case we
have
\begin{equation*}
\omega \beta^{-1} \frac{1}{{\mathbb F}} \frac{\partial {\mathbb
F}}{\partial z} < \omega \beta^{-1} \sqrt{\bf k} << 1.
\end{equation*}
From $(\theta, z) \in {\mathcal F}^{-1}({\mathcal F}(V) \cap V)$
and ${\mathcal F}(V_f) \subset U$, it follows that $(\theta, z)
\not \in V_f$ therefore
\begin{equation*}
\left|\frac{\partial \theta_1}{\partial \theta}\right| = \left|1-
\omega \beta^{-1} \frac{1}{{\mathbb F}} \frac{\partial {\mathbb
F}}{\partial \theta}\right|
> 2.
\end{equation*}
These two estimates together implies that the denominator for
$|s(D{\mathcal F}({\bf v}))|$ in (\ref{f2-add-s3.2}) is $> 1$, and
it follows that $|s(D{\mathcal F}({\bf v}))| < \frac{1}{100}$.

\smallskip

{\it Case 2: ${\mathbb F} < \sqrt{\bf k}$.} \ In this case
$$
|1 +
{\bf c} \sin \theta| < K \varepsilon^{\frac{\beta}{\alpha}} +
\sqrt{\bf k},
$$
from which we have
\begin{equation} \label{f3-add-s3.2}
|{\bf c} \cos \theta|
> 1.
\end{equation}
It then follows that the denominator for $|s(D{\mathcal F}({\bf
v}))|$ in (\ref{f2-add-s3.2}) is $
> \frac{1}{2} \frac{1}{\sqrt{\bf k}}$, which implies
$|s(D{\mathcal F}({\bf v}))| < \frac{1}{100}$. This finishes our
proof for (i).

\smallskip

To prove (ii) we let ${\bf v}$ be such that $|s({\bf v})|
> 100$. From (\ref{f1-add-s3.2}),
\begin{equation}\label{f4-add-s3.2}
D{\mathcal F}^{-1} = \frac{1}{\alpha \beta^{-1} {\bf b} {\mathbb
F}^{\alpha \beta^{-1} -1} \frac{\partial {\mathbb F}}{\partial z}}
\left(
\begin{array}{cc} \alpha \beta^{-1}
{\bf b} {\mathbb F}^{\alpha \beta^{-1} -1} \frac{\partial {\mathbb
F}}{\partial z} &  - \omega \beta^{-1} \frac{1}{{\mathbb F}}
\frac{\partial {\mathbb F}}{\partial z} \\
- \alpha \beta^{-1}  {\bf b} {\mathbb F}^{\alpha \beta^{-1}-1}
\frac{\partial {\mathbb F}}{\partial \theta} & 1 - \omega
\beta^{-1} \frac{1}{{\mathbb F}} \frac{\partial {\mathbb
F}}{\partial \theta}
\end{array} \right),
\end{equation}
and we have
\begin{equation*}
|s(D{\mathcal F}^{-1}({\bf v}))| =  \left|\frac{- \alpha
\beta^{-1} {\bf b} {\mathbb F}^{\alpha \beta^{-1}-1}
\frac{\partial {\mathbb F}}{\partial \theta}s^{-1}({\bf v}) + (1 -
\omega \beta^{-1} \frac{1}{{\mathbb F}} \frac{\partial {\mathbb
F}}{\partial \theta})}{\alpha \beta^{-1} {\bf b} {\mathbb
F}^{\alpha \beta^{-1} -1} \frac{\partial {\mathbb F}}{\partial z}
s^{-1}({\bf v}) - \omega \beta^{-1} \frac{1}{{\mathbb F}}
\frac{\partial {\mathbb F}}{\partial z}} \right|
\end{equation*}
We again divide into the cases of ${\mathbb F} > \sqrt{k}$ and
${\mathbb F} < \sqrt{k}$. For the case of ${\mathbb F} >
\sqrt{k}$, the magnitude of the denominator $<< 1$ and that of the
numerator is $> 1$ again because
$$
\left|1 - \omega \beta^{-1} \frac{1}{{\mathbb F}} \frac{\partial
{\mathbb F}}{\partial \theta}\right| > 2
$$
from the assumption that $(\theta, z) \not \in  V_f$. For the case
of ${\mathbb F} < \sqrt{\bf k}$, we re-write $|s(D{\mathcal
F}^{-1}({\bf v}))|$ as
\begin{equation*}
|s(D{\mathcal F}^{-1}({\bf v}))|=  \left|\frac{- \alpha \beta^{-1}
{\bf b} {\mathbb F}^{\alpha \beta^{-1}} \frac{\partial {\mathbb
F}}{\partial \theta}s^{-1}({\bf v}) + ({\mathbb F} - \omega
\beta^{-1} \frac{\partial {\mathbb F}}{\partial \theta})}{\alpha
\beta^{-1} {\bf b} {\mathbb F}^{\alpha \beta^{-1}} \frac{\partial
{\mathbb F}}{\partial z} s^{-1}({\bf v}) - \omega \beta^{-1}
\frac{\partial {\mathbb F}}{\partial z}} \right|.
\end{equation*}
The denominator is again $<< 1$ and the dominating term in the
numerator is
$$
\left|\omega \beta^{-1} \frac{\partial {\mathbb F}}{\partial
\theta} \right| > 1.
$$
The last estimate is from
$$
\left|\frac{\partial {\mathbb F}}{\partial \theta} \right| \approx
|{\bf c} \cos \theta| > 1
$$
again by (\ref{f3-add-s3.2}). This proves (ii). \hfill $\square$

\smallskip

We refer the reader to Chapter III.1 of \cite{Mo} for a detailed
discussion on horseshoes of infinitely many symbols.

\medskip

{\bf Remarks:} \ 1. For the parameters of Theorem \ref{th1}, the
entire homoclinic tangle consists of a single horseshoe of
infinitely many symbols.

2. Due to the expansions associated with singularities of the
logarithmal function in (\ref{f1-s3.1a}), ${\mathcal F}_{\bf a}$
induces a horseshoe of infinitely many symbols on $V \setminus
V_f$ {\bf for all $|\mu| < \mu_0$}. This horseshoe covers Smale's
horseshoe and all its variations. It is the one that resides
inside all homoclinic tangles.

3. $\Lambda$ is much more complicated when ${\mathcal F}_{\bf
a}(V_f)$ intersects $V$. As ${\mathcal F}_{\bf a}(V_f)$ traverses
$V$, we encounter complicated dynamical patterns caused by our
allowing the images of the unstable manifold of the horseshoe in
$V \setminus V_f$ (see remark 2) to come back to traverse the
stable manifold of the same horseshoe. We will prove, momentarily,
that there are parameters that admit periodic sinks and there are
also others that admit non-degenerate transversal homoclinic
tangency. We unfortunately do not have a bifurcation diagram for
${\mathcal F}_{\bf a}$. However, we know from (\ref{f1-s3.1a})
that the same diagram are repeated infinitely many times as $\mu
\to 0$.

4. We caution that, though the horseshoe of Theorem \ref{th1}
represents all solutions of the perturbed equation that stay
forever inside of a small neighborhood of the homoclinic loop
$\ell$, solutions sneaked out through $U$ might find a way to come
back to ${\mathcal A}$, creating more complicated structures. One
particular mechanism for such coming back is for the unperturbed
equation to have two homoclinic solutions. See Fig. 6(a). In this
case, part of $U$ would come back to ${\mathcal A}$ following the
other homoclinic loop. On the other hand, it is easy to obtain
examples for which the solutions sneaked out of $U$ would never
come back. In this case the entire homoclinic tangle for the
perturbed equation is in fact reduced to the horseshoe of Theorem
\ref{th1}: all it takes for this to happen is for us to send the
other branch of the local unstable manifold of $(0, 0)$ to a sink.
See Fig. 6(b).

\begin{picture}(10, 5.5)
\put(1.5,0){ \psfig{figure=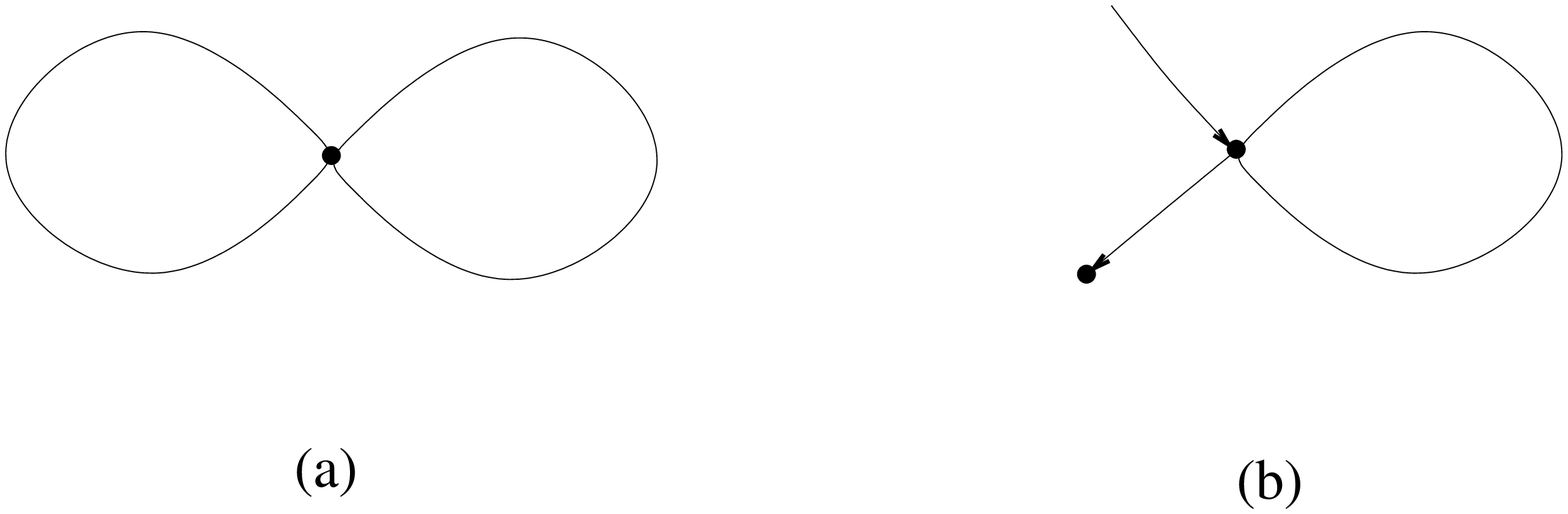,height = 4.5cm,
width = 10cm} }
\end{picture}

\medskip

\centerline{Fig. 6 \ (a) $U$ gets back to ${\mathcal A}$, and (b)
all points in $U$ approach a sink.}

\medskip

Our next Theorem is about the existence of periodic sinks. We
remark that these periodic sinks are {\bf not} Newhouse sinks
associated with homoclinic tangency.

\begin{theorem}[{\bf Periodic sinks}]\label{th2}
Let the assumptions be identical to that of Theorem \ref{th1}.
Then there exists an open set of $\mu$ inside each of the
intervals $[\mu_n^{(r)}, \mu_{n+1}^{(l)}]$, such that the
corresponding homoclinic tangle admits a periodic sink.
\end{theorem}
\noindent {\bf Proof:} \ Let $\theta_c(z)$ be as in Lemma
\ref{lem1-s3.2}(b). To make the dependency on $\mu$ explicit we
write it as $\theta_c(z, \mu)$. Let ${\bf a}_n$ be the value of
${\bf a}$ at $\mu = \mu_n^{(r)}$ and $[{\bf a}_n] = {\bf a}_n -
{\bf a}_n  mod  (2 \pi)$. Observe that there exists a $\hat \mu
\in [\mu_n^{(r)}, \mu_{n+1}^{(l)}]$ so that $\theta_1(\theta_c) =
\theta_c + [{\bf a}_n]$ where $\theta_c = \theta_c(0, \hat \mu)$.
This is because when $\mu$ traverses $[\mu_n^{(r)},
\mu_{n+1}^{(l)}]$, $\theta_1(\theta_c)$ traverses the interval
$(\theta_l + [{\bf a}_n], \theta_r + [{\bf a}_n])$. Let $\hat {\bf
a}$ be the value of ${\bf a}$ for $\hat \mu$. To solve for a fixed
point we let
\begin{equation}\label{f4-add-s3.2}
\begin{split}
\theta + [{\bf a}_n] & = \theta + \hat {\bf a} - \omega \beta^{-1} \ln {\mathbb F} \\
z & = {\bf b} {\mathbb F}^{\alpha \beta^{-1}}
\end{split}
\end{equation}
to obtain
\begin{equation}\label{f5-add-s3.2}
\begin{split}
{\mathbb F} & = e^{\omega^{-1} \beta (\hat
{\bf a}-[{\bf a}_n])}, \\
z & = {\bf b} e^{\omega^{-1} \alpha (\hat {\bf a}- [{\bf a}_n])}.
\end{split}
\end{equation}
From the first line we have
\begin{equation}\label{f6-add-s3.2}
 1 + {\bf c} \sin \theta + {\mathbb E}(\theta, \hat \mu) +
{\mathcal O}_{\theta, z, {\bf a}}(\hat \mu) = e^{\omega^{-1} \beta
(\hat {\bf a}-[{\bf a}_n])}.
\end{equation}
To solve (\ref{f6-add-s3.2}) for $\theta$, first we observe that
$\theta_c = \theta_c(0, \hat \mu)$ is a solution of
(\ref{f6-add-s3.2}) for $z = 0$. We then observe that
$$
|\cos \theta_c| > K^{-1}.
$$
This estimate follows from the fact that $\theta_c$ is defined by
$$
1 - \omega \beta^{-1} \frac{1}{\mathbb F} \frac{\partial {\mathbb
F}}{\partial \theta} = 0
$$
and ${\mathbb F} = e^{\omega^{-1} \beta (\hat {\bf a} - [{\bf
a}_n])}$ from (\ref{f5-add-s3.2}). Applying the inverse value
theorem to (\ref{f6-add-s3.2}) we obtain a solution $\hat \theta$
satisfying
$$
|\hat \theta - \theta_c| < K \hat \mu.
$$
In summary we have obtained a fixed point $(\hat \theta, \hat z)$
satisfying
$$
\hat \theta \approx \theta_c; \ \ \ \hat z = {\bf b}
e^{\omega^{-1} \alpha (\hat {\bf a}-[{\bf a}_n])}.
$$

To prove that $(\hat \theta, \hat z)$ is an attracting fixed
point, we compute the eigenvalues. The eigen-equation for
$D{\mathcal F}$ is
$$
\lambda^2 - Tr(D{\mathcal F}) \lambda + \det(D{\mathcal F}) = 0.
$$
From (\ref{f1-add-s3.2}) we have
\begin{equation}\label{f7-add-s3.2}
\begin{split}
Tr(D{\mathcal F}) & = \frac{\partial \theta_1}{\partial \theta} +
\alpha \beta^{-1} {\bf b} {\mathbb F}^{\alpha \beta^{-1} -1}
\frac{\partial {\mathbb
F}}{\partial z} << 1  \\
\det(D{\mathcal F}) & = \alpha \beta^{-1} {\bf b} {\mathbb
F}^{\alpha \beta^{-1} -1} \frac{\partial {\mathbb F}}{\partial z}
<< 1
\end{split}
\end{equation}
where for the first inequality we use
$$
\frac{\partial \theta_1}{\partial \theta} < K (|\hat \theta -
\theta_c| + |\hat z|)
$$
at $(\hat \theta, \hat z)$ with
$$
K = \max_{\theta \in (\theta_c, \hat \theta), z \in [0, \hat z]}
\left( \left|\frac{\partial^2 \theta_1}{\partial \theta^2} \right|
+ \left|\frac{\partial^2 \theta_1}{\partial \theta
\partial z}\right| \right).
$$
Note that, on the domain the maximum is taken, ${\mathbb F}>
\frac{1}{2}$. The rest of (\ref{f7-add-s3.2}) are obvious. It
follows  from (\ref{f7-add-s3.2}) that both eigenvalues of
$D{\mathcal F}$ are close to $0$. \hfill $\square$

\medskip

Our next Theorem is about the existence of non-degenerate
transversal homoclinic tangency.
\begin{theorem}[{\bf Homoclinic tangency}]\label{th3}
Let the assumptions be identical to that of Theorem \ref{th1}.
Then for every $n >0$ given, there exists $\hat \mu \in
[\mu^{(r)}_n, \mu^{(l)}_{n+1}]$, the corresponding value for ${\bf
a}$ we denote as $\hat {\bf a}$, such that

(i) ${\mathcal F}_{\hat {\bf a}}$ has a saddle fixed point, which
we denote as $q(\hat {\bf a})$, so that $W^u(q(\hat {\bf a})) \cap
W^s(q(\hat {\bf a}))$ contains a point of non-degenerate tangency.

(ii) Let $q({\bf a})$ be the continuous extension of $q(\hat {\bf
a})$ for ${\bf a}$ sufficiently close to $\hat {\bf a}$. Then as
${\bf a}$ passes through $\hat {\bf a}$, $W^u(q({\bf a}))$ crosses
$W^s(q({\bf a}))$ at the tangential intersection point of (i) with
a relative speed $> \frac{1}{2}$ with respect to ${\bf a}$  in
$\theta$-direction.
\end{theorem}
\noindent {\bf Proof:} \ Our plan of proof is as follows. We know
that ${\mathcal F}_{\bf a}$ induces a horseshoe of infinitely many
symbols in $V \setminus V_f$, creating many saddle fixed points.
Pick one and denote it as $q$. We prove that $q$ is continuously
extended over the $\mu$ interval $[\mu^{(r)}_n, \mu^{(l)}_{n+1}]$,
which we denote as $q({\bf a})$. Let $W^u(q({\bf a}))$ be the
unstable and $W^s(q({\bf a}))$ be the stable manifold of $q({\bf
a})$. We prove that $W^u(q({\bf a})) \cap V_f$ has a horizontal
segment traversing $V_f$, which we denote as $\ell^u({\bf a})$. We
also prove that $W^s(q({\bf a}))$ has a vertical segment fully
extended in $V$, which we denote as $\ell^s({\bf a})$. Observe
that ${\mathcal F}_{\bf a}(\ell^u({\bf a}))$ has a sharp quadratic
turn, and as $\mu$ varies from $\mu_n^{(r)}$ to $\mu_{n+1}^{(l)}$,
it moves from one side of $V$ to the other, transversally crossing
$\ell^s({\bf a})$. See Fig. 7.

\begin{picture}(12, 4.5)
\put(1,0){ \psfig{figure=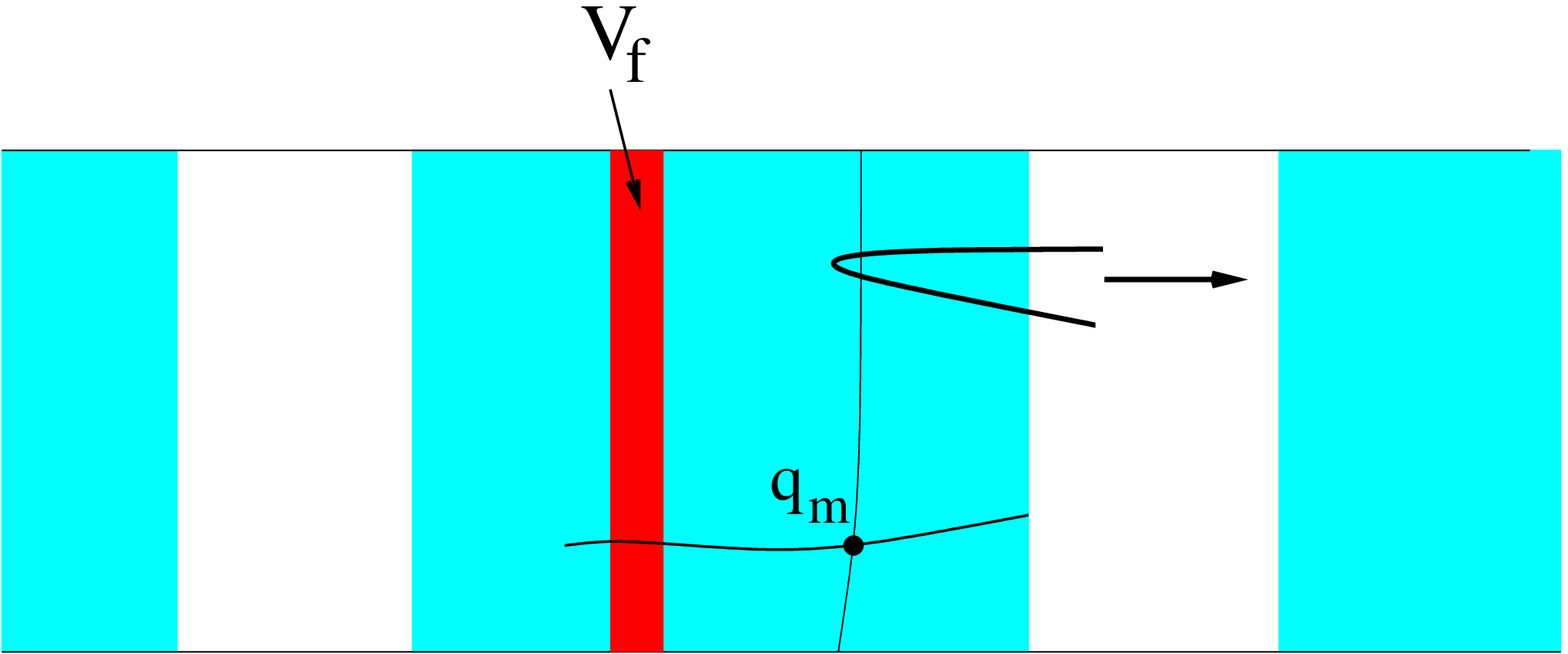,height = 4cm, width =
12cm} }
\end{picture}

\vskip .2in

\centerline{Fig. 7 \ Transversal homoclinic tangency.}

\vskip .2in

Detailed proof for Theorem \ref{th3} is long and include some
tedious computations. A complete proof is included in the
Appendices. \hfill $\square$

\medskip

The following is a direct consequence of Theorem \ref{th3}.
\footnote{We thank Marcelo Viana for assuring us that, with
Theorem \ref{th3},  \cite{MV} directly applies.}

\begin{coro}\label{coro-add-s3}
Let the assumptions be identical to that of Theorem \ref{th1}.
Then inside of every parameter interval $[\mu^{(r)}_n,
\mu^{(l)}_{n+1}]$, there is a set of parameters of positive
Lebesgue measure, such that the homoclinic tangle associated with
these parameters admits strange attractors with SRB measures.
\end{coro}
\noindent  {\bf Proof:} \ This follows from Theorem \ref{th3}
applying \cite{MV} and \cite{BY}, both are based on [BC], to
${\mathcal F}_{\bf a}$. \hfill $\square$

\subsection{Homoclinic tangles and observable chaos}\label{s3.3}
Let ${\mathcal F}_{\bf a}$ be as in (\ref{f1-s3.1a}) and $\Omega$,
$\Lambda$ be as in (\ref{f2-s3.1b}). $\Omega$ represents all
solutions that stay close to $\ell$ in forward times, and
$\Lambda$ represents all solutions that stay close to $\ell$ in
both the forward and the backward times. In this subsection we
study numerically the structures of $\Omega$ and $\Lambda$.

We start with a concept of {\it observability} in numerical
simulations. We say that a homoclinic tangle is {\it observable}
in phase space if $\Omega$ has positive Lebesgue measure.
Otherwise we say that this homoclinic tangle is {\it not
observable}. We only expect observable homoclinic tangles to show
up in numerical simulations. For maps with parameters, there is
also an issue of observability in parameter space: a
sub-collection of maps is observable only if it is from a
parameter set of positive Lebesgue measure. See \cite{WOk} for
more detailed discussions on observable dynamical scenarios in
numerical simulations.

To numerically study homoclinic tangles through ${\mathcal F}_{\bf
a}$, we drop the error terms in (\ref{f1-s3.1a}) and re-write
${\bf k} z$ as $z$. We obtain from
(\ref{f1-s3.1a})-(\ref{f4-s3.1a}) a family of 2D maps in the form
of
\begin{equation}\label{f1-s3.1d}
\begin{split}
\theta_1 & = \theta + a - d
\ln (1 + c \sin \theta + z) \\
z_1 & =  b [1 + c \sin \theta + z]^{\gamma}.
\end{split}
\end{equation}
where $a, b, c, d, \gamma$ are parameters. $a \in S^1$, $b << 1$,
$c > 1$, $d \in {\mathbb R}$ and $\gamma > 1$. From Theorems
\ref{th1}-\ref{th3} and Corollary \ref{coro-add-s3}, we would
expect at least three dynamical scenarios that are observable in
parameter space. They are as follows.

\medskip

(i) For parameters of Theorem \ref{th1}, $\Lambda$ for ${\mathcal
F}_{\bf a}$ is a uniformly hyperbolic invariant set, and $\Omega$
is the stable manifold for $\Lambda$ inside of $\Sigma^-$. Both
$\Omega$ and $\Lambda$ are Lebesgue measure zero sets. The
corresponding homoclinic tangle for these parameters is therefore
not observable in phase space.

\smallskip

(ii) For parameters of Theorem \ref{th2}, $\Omega$ contains an
open neighborhood of a periodic sink so the associated homoclinic
tangle is observable. Plots of individual orbits from $\Omega$
would lead us to periodic sinks in $\Lambda$.

\smallskip

(iii) For parameters of Corollary \ref{coro-add-s3}, we expect
strange attractors with SRB measures to show up as an observable
phenomenon in numerical simulations.

\medskip

In Figs 8-10, we plot $\Omega$ and $\Lambda$ for the maps defined
in (\ref{f1-s3.1d}) with various choices of parameters that
reflect the scenarios (i)-(iii) above respectively. Fig. 8 is for
scenario (i), with $a=0.2, b = 0.005, c=3, d = 2$ and $\gamma =
\sqrt{2}$. Fig. 8(a) is a plot of all points in $V$, the orbits of
which remain inside of $V$ after 3 iterations, Fig. 8(b) is for
after 6 iterations. Nothing is left in $V$ after 15 iterations.

\begin{picture}(5, 5)
\put(1,0){ \psfig{figure=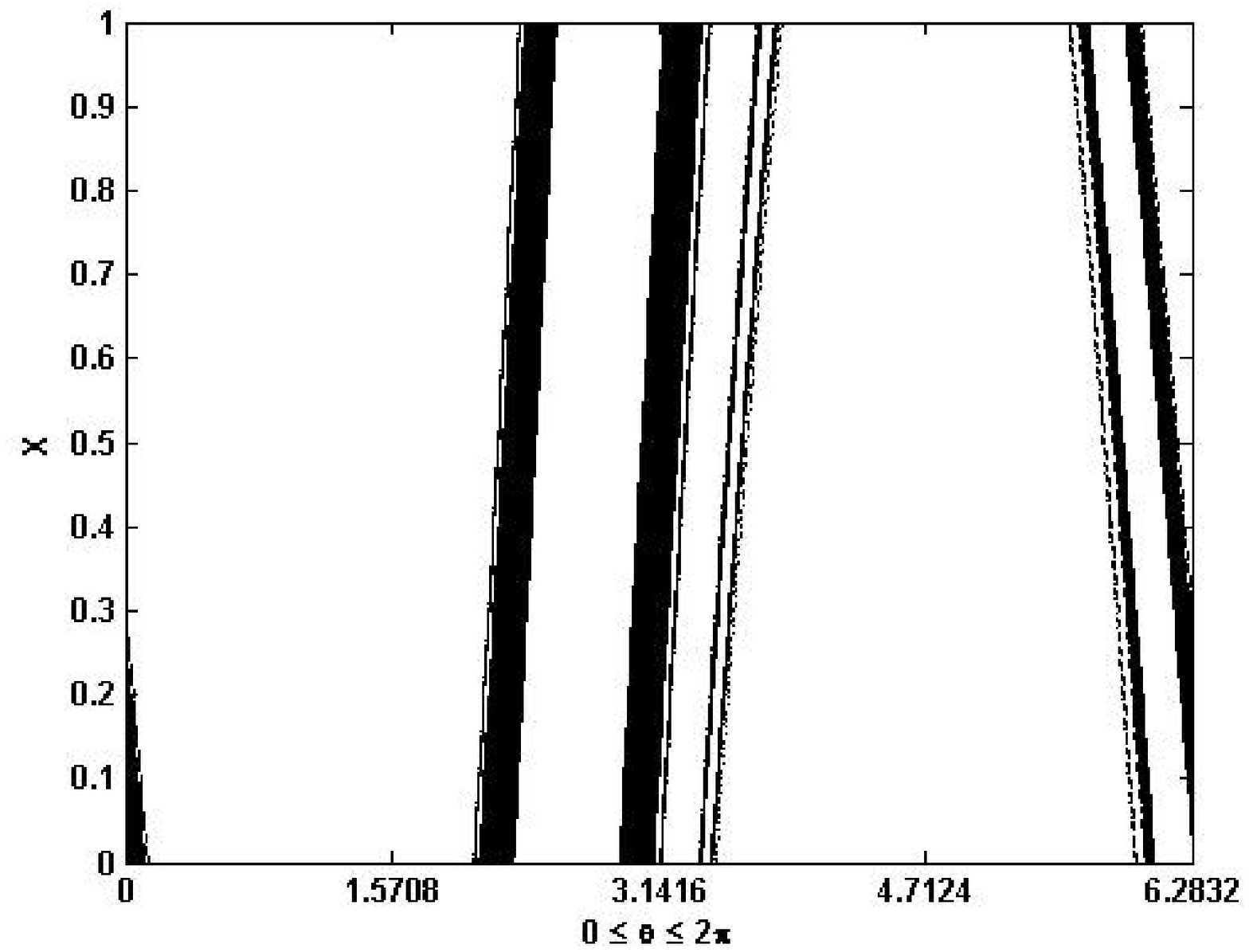,height = 5cm, width =
5cm} }
\end{picture}
\begin{picture}(9, 6)
\put(3, 0){ \psfig{figure=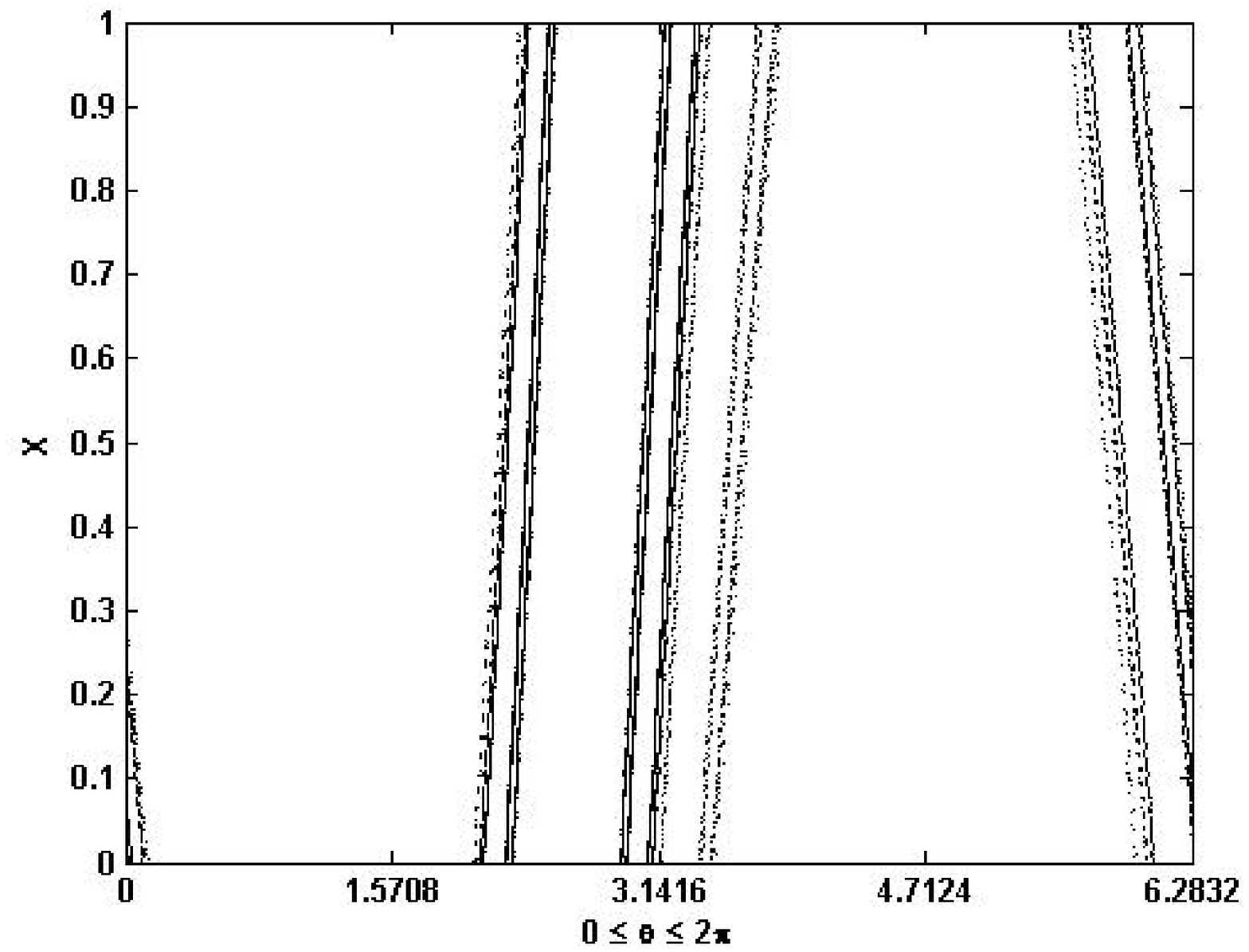,height = 5cm, width =
5cm} }
\end{picture}

\hspace{3.5cm} (a) \hspace{6cm} (b)

\smallskip

\centerline{Fig. 8 \ Homoclinic tangles with no sinks nor
observable chaos.}

\centerline{($a = 0.2, b = 0.005, c = 3, d = 2$ and $\gamma =
\sqrt{2}$)}

\medskip

Fig. 9 is for scenario (ii) with  $a = 2$. The values for $b$,
$c$, $d$, $\gamma$ are kept the same as in Fig. 8. In this case
$\Lambda$ contains a periodic sink with a relatively large basin.
Fig. 9(a) is for $\Omega$. All orbits initiated from $\Omega$
quickly converge to an attracting periodic orbit. In Fig. 9(b) we
depict $\theta_k$ v.s. $k$ for one orbit from $\Omega$. Picture
for $z_k$ v.s. $k$ is similar. Only one periodic sink shows up for
$\Lambda$ in numerical simulations. The horseshoe associated with
the singularity of the logarithmal function (Smale's horseshoe),
though exists inside of $\Lambda$, does not show up because the
set it attracts is a set of zero Lesbegue measure.

\begin{picture}(5, 5)
\put(1,0){ \psfig{figure=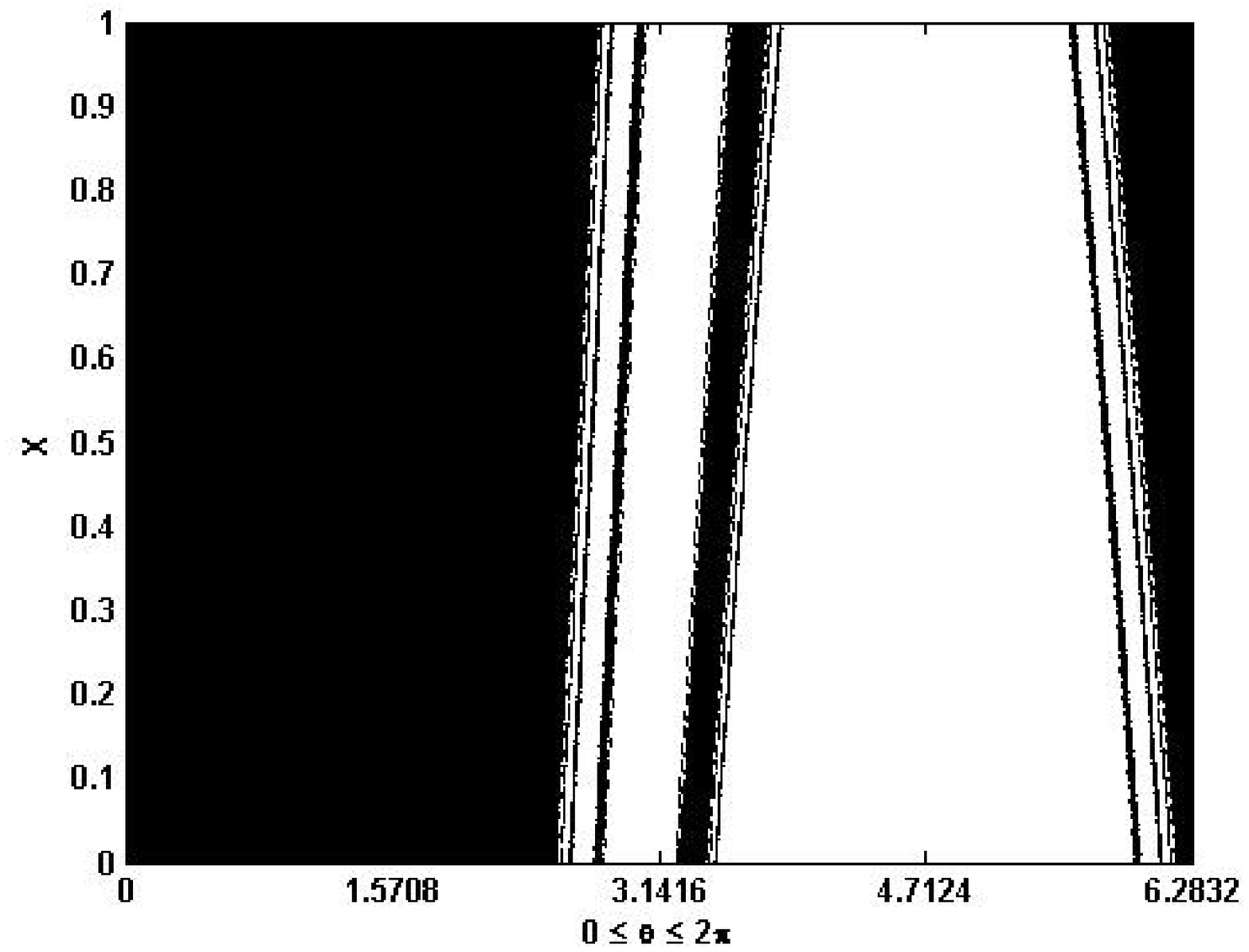,height = 5cm, width =
5cm} }
\end{picture}
\begin{picture}(9, 6)
\put(3, 0){ \psfig{figure=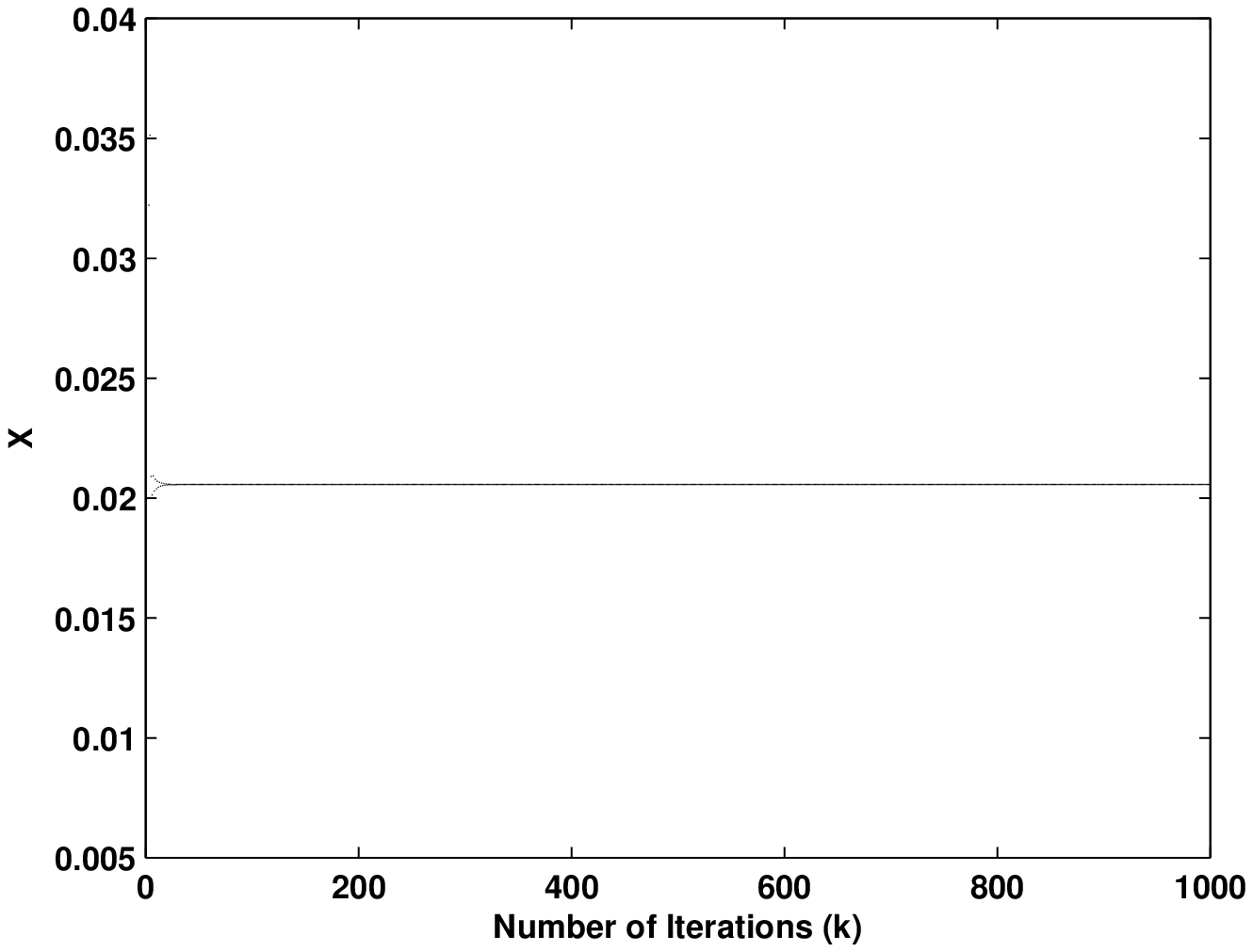,height = 5cm, width =
5cm} }
\end{picture}

\hspace{3.5cm} (a) \hspace{6cm} (b)

\smallskip

\centerline{Fig. 9 \ Homoclinic tangles with an attracting
periodic solution.}

\centerline{($a = 2, b = 0.005, c = 3, d = 2$ and $\gamma =
\sqrt{2}$)}

\medskip

Fig. 10 is for scenario (iii) with $a = 1.5$. The values for $b,
c, d$ and $\gamma$  are kept the same as before. $\Omega$ is
depicted in Fig. 10(a). In Fig. 10(b) we depicted again $\theta_k$
v.s. $k$ for one orbit. As $k$ moves forward, $\theta_k$ jumps
randomly in a fixed range. These pictures represent a strange
attractor with an SRB measure associated to transversal homoclinic
tangency of a saddle periodic orbit of relatively high period. All
orbits from $\Omega$ in fact offer the same picture.

\begin{picture}(5, 5)
\put(1,0){ \psfig{figure=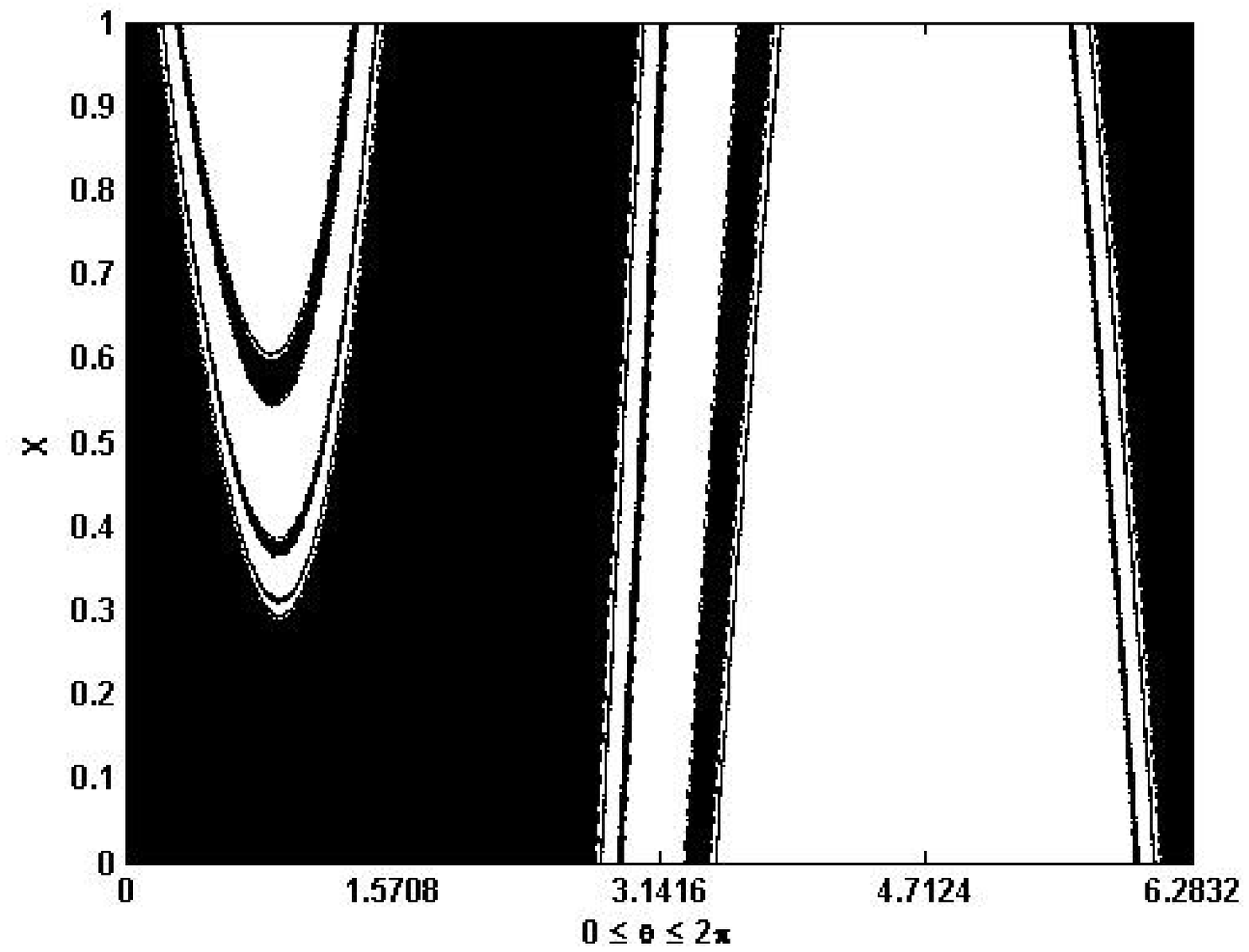,height = 5cm, width =
5cm} }
\end{picture}
\begin{picture}(9, 6)
\put(3, 0){ \psfig{figure=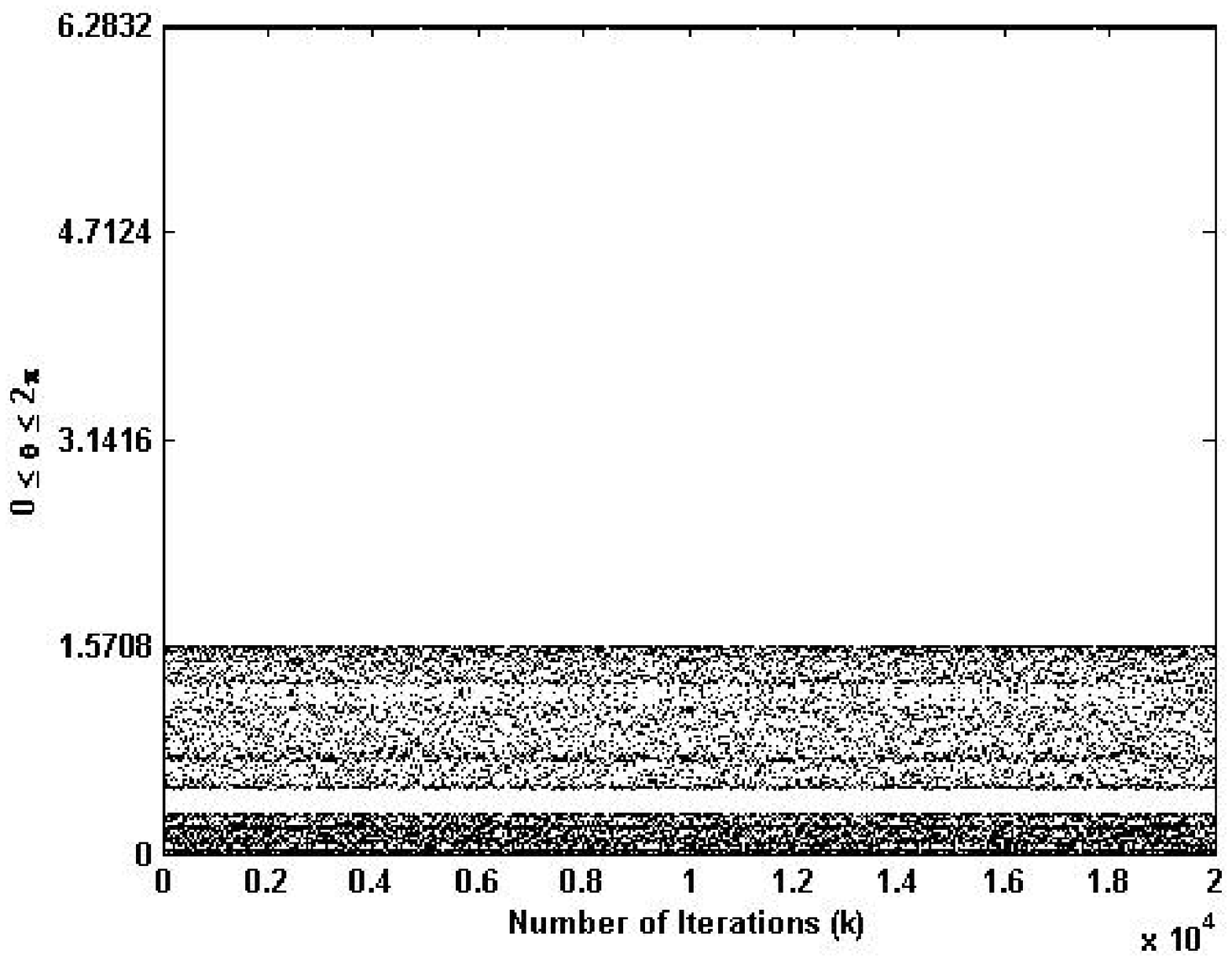,height = 5cm, width =
5cm} }
\end{picture}

\hspace{3.5cm} (a) \hspace{6cm} (b)

\smallskip

\centerline{Fig. 10 \ Tangles with an observable chaos.}

\centerline{($a = 1.5, b = 0.005, c = 3, d = 2$ and $\gamma =
\sqrt{2}$)}

\medskip

We also performed systematic search over all combinations of
parameters with $b$ reasonably small. We persistently run into one
of the three scenarios above. In the case of Fig. 8, however,
sometimes it takes much longer for all points to be completely
iterated out of $V$. This is particularly the case when $d$ is
small, and is more or less expected: as the overall strength of
expansions around the horseshoe of Theorem \ref{th1} gets weaker,
the points in $V$ tends to linger longer inside of $V$ before been
pushed out into $U$.

\subsection{Dynamical scenarios for periodically
perturbed homoclinic solutions} \label{s3.4}

In this paragraph we summarize all that have been obtain so far
for equation (\ref{f2-s1}) in \cite{WO}, \cite{LW} and in this
paper through the return map of Proposition \ref{prop1-s2.3c}.
Again, we let ${\mathcal Q}(t) = \sin t$ and assume (H1) and (H2).
The forcing parameters are inside of
$$
{\mathbb P} = \{ (\omega, \rho, \mu): \ \omega \in (0,
R_{\omega}), \ \rho \in (R_{\rho}^{-1}, R_{\rho}), \ \mu \in (0,
R_{\mu}^{-1}) \}
$$
where $R_{\mu} >> R_{\rho} >> R_{\omega} >> 1$. Let $W^s$ be the
stable and $W^u$ be the unstable manifold of $(x, y) = (0, 0)$ in
the extended phase space. Various dynamics scenarios for different
parts of ${\mathbb P}$ are illustrated in Fig. 11. Since our
purpose is to provide an overview, only descriptive statements are
presented. Rigorous formulations and their proofs are either
directly included in \cite{WO}, \cite{LW} and in this paper, or
obtained by reasonable modifications of existing text.

\smallskip

1. There is a surface $S^*$ in ${\mathbb P}$ (See Fig. 11), such
that for all parameters under $S^*$, $W^u  \cap W^s \neq
\emptyset$ and we have homoclinic tangles for equation
(\ref{f2-s1}). The dynamics of these homoclinic tangles are
studied in Sect. \ref{s3.2} of this paper. In particular, there
are open sets of parameters, such that the entire homnoclinic
tangle is one uniformly hyperbolic horseshoe. There are also
parameters for periodic sinks, and parameters for non-degenerate,
transversal homoclinic tangency.

\begin{picture}(9, 7.5)
\put(2,0){ \psfig{figure=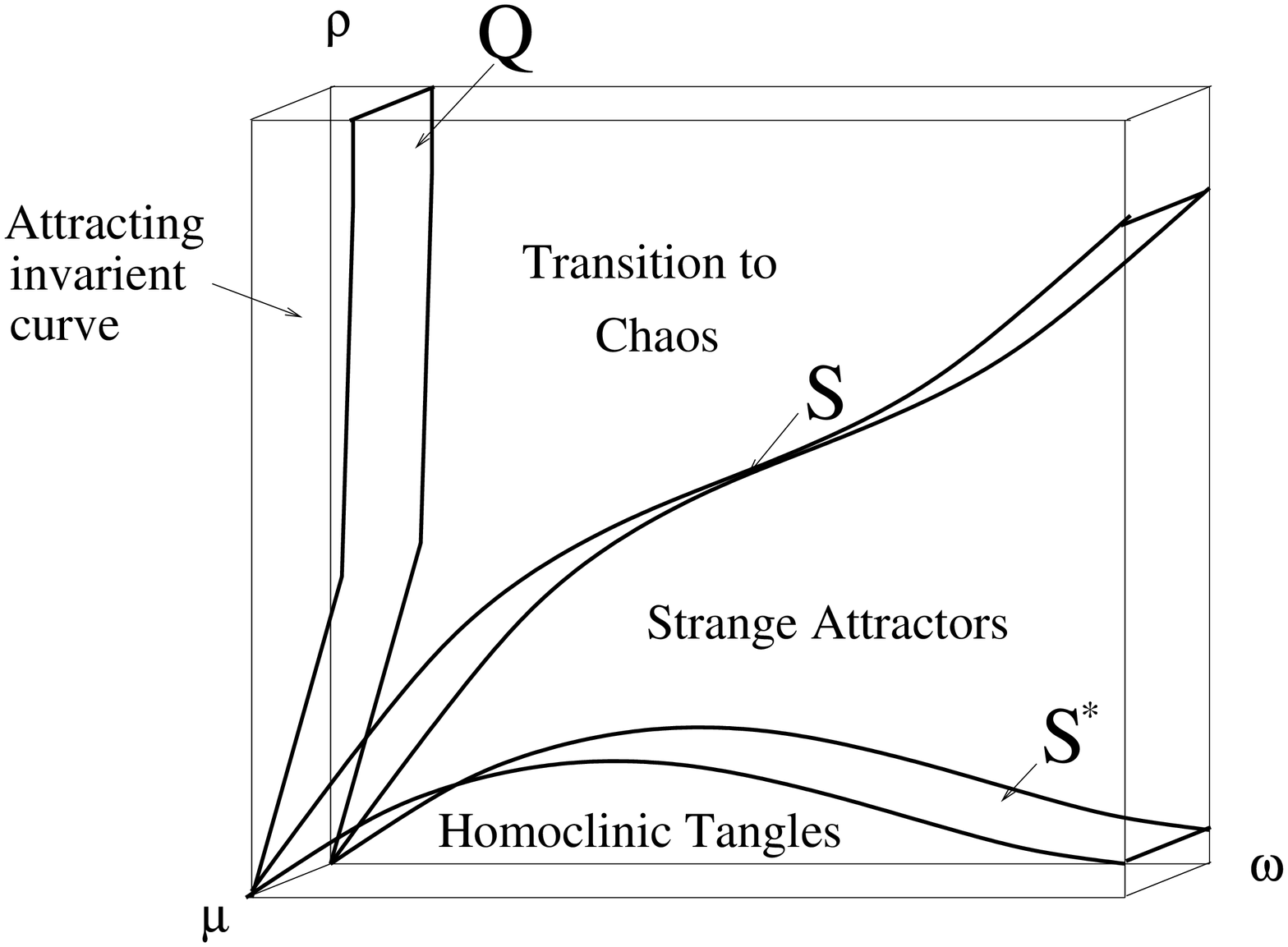,height = 7.5cm, width =
8cm} }
\end{picture}

\vskip .1in

\centerline{Fig. 11 \ Dynamical scenarios in parameter space.}

\vskip .1in

2. For parameters over $S^*$, $W^s \cap W^u = \emptyset$. The
return maps are again defined through
(\ref{f1-s3.1a})-(\ref{f5-s3.1a}), but for these parameters
${\mathbb F} > 0$ on $\Sigma^-$ so ${\mathcal F}_{\bf a}$ are
well-defined on $\Sigma^-$. These maps have been studied
systematically in \cite{WY4}. We know that
\begin{itemize}
\item[(a)] There is a surface $S$ above $S^*$ for which the
following holds. For all parameters in between $S$ and $S^*$,
${\mathcal F}_{\bf a}$ admit global attractors in $\Sigma^-$ that
are {\em chaotic} in the sense that they all contain a horseshoe
(See \cite{LW}). If the forcing frequency $\omega$ is reasonably
large, then there is a positive measure set of forcing parameters
such that the strange attractors are rank one attractors of
\cite{WY1} and \cite{WY2} with SRB measures (See \cite{WO}).

\item[(b)] There is a surface $Q$ (See Fig. 11) for which the
following holds. For any give set of parameters on the left of
$Q$, $\Sigma^-$ is attracted globally to a simple closed curve, on
which the induced map conjugates to a circle diffeomorphisms. If
the rotation number of this circle diffeomorphism is rational,
then there are saddles and attracting periodic solutions. If the
rotation number is irrational, then the solutions are
quasi-periodic. We also know as a fact that there are positive
measure sets of parameters such that the corresponding rotation
numbers are irrational and the corresponding solutions are
quasi-periodic.

\item[(c)] What happens between $Q$ and $S$ are as follows: as we
move from the left to the right in the $\omega$-direction, larger
forcing frequency first deforms, then breaks the attracting
invariant curve, inducing sinks and saddles. The unstable
manifolds of these induced saddles would eventually fold in
$\theta$-direction, intersecting the stable manifolds to create
strange attractors and rank one chaos. As we go down in the
$\rho$-direction, $W^u$ and $W^s$ are pulled gradually closer.
Reflected in the return maps of Proposition \ref{prop1-s2.3c} is
the growing relevance of the expansions associated with the
singularity of the logarithmical function of (\ref{f1-s3.1a}).
\end{itemize}

\section{Homoclinic tangles for general forcing functions}
\label{s4}

In this section we let ${\mathcal Q}(t)$ be an arbitrary periodic
function of period $2 \pi$.  We explain how the different choices
of the forcing function ${\mathcal Q}(t)$ affect the dynamics of
the associated homoclinic tangles.

Let $A$, $S(\omega)$ and $C(\omega)$ be the same as before (See
(\ref{f1-s3.1})). $A$, $S(\omega)$ and $C(\omega)$ are independent
of ${\mathcal Q}(t)$. We assume (H1) and (H2)(i) for equation
(\ref{f2-s1}) and replace (H2)(ii) by (H3) below.

\medskip

\noindent {\bf (H3)} \ There exists a constant $\xi
> 0$ so that
$$
\sqrt{S^2(\omega) + C^2(\omega)} \sim e^{-\xi |\omega|}
$$
as $| \omega| \to + \infty$.

\medskip

(H3) is stronger than (H2)(ii). It requires that the magnitude of
the Fourier transformation $\hat R(\omega)$ of the function
$$
R(s) = {\mathbb H}(s) e^{-\int_0^s E(\tau) d \tau}
$$
decays exponentially as $|\omega| \to \infty$. We note that there
is no lack of known systems satisfying (H3).

Expanding ${\mathcal Q}(t)$ in Fourier series we write
\begin{equation}\label{f1-s4.2}
{\mathcal Q}(t) = \sum_{n=1}^{\infty} (c_n \cos{nt} + s_n
\sin{nt}).
\end{equation}
If the mean value of a ${\mathcal Q}(t)$ is not zero, we give it
to $\rho$. So there is no loss of generality in assuming
(\ref{f1-s4.2}). Let us assume in addition that
\begin{equation}\label{f2-s4.2}
c_1^2 + s_1^2 \neq 0.
\end{equation}

Let ${\mathcal Q}(t)$ be as in (\ref{f1-s4.2}) satisfying
(\ref{f2-s4.2}) and assume (H1), (H2)(i) and (H3) for equation
(\ref{f2-s1}). Parameters $\omega, \rho, \varepsilon, \mu$ are
specified as follows. First we fix (arbitrarily) an $\omega \in
[\frac{1}{100} R_{\omega}, R_{\omega}]$. We then fix a value of
$\rho$ such that
$$
3 <\frac{1}{\rho A} \cdot \sqrt{c_1^2 + s_1^2} \cdot
\sqrt{C^2(\omega) + S^2(\omega)} < 9.
$$
After that we fix $\varepsilon$ sufficiently small so that,
$$
|\phi(\theta) - \phi_L(\theta)| << \rho A
$$
where $\phi_L(\theta)$ is from (\ref{f1a-s2.3b}) and
$\phi(\theta)$ is obtained by replacing $-L^-, L^+$ with $-\infty,
+\infty$ respectively in $\phi_L(\theta)$. We also make
$\varepsilon$ sufficiently small so that
$$
2 <\frac{1}{\rho A} \cdot \sqrt{c_1^2 + s_1^2} \cdot
\sqrt{C^2_L(\omega) + S^2_L(\omega)} < 10.
$$
$\mu$ ($<< \varepsilon$) is the only parameter we allow to vary.
In what follows ${\mathcal F}_p: \Sigma^- \to \Sigma^-$ is the
return maps of Proposition \ref{prop1-s2.3c} induced by equation
(\ref{f2-s1}). Recall that $p = \ln \mu$.
\begin{theorem}\label{th4}
Let ${\mathcal Q}(t)$ be as in (\ref{f1-s4.2}) satisfying
(\ref{f2-s4.2}) and assume (H1), (H2)(i) and (H3) for equation
(\ref{f2-s1}). Let the parameters $\omega, \rho, \varepsilon$ be
specified as in the above. Then

(a) there are infinitely many disjoint open intervals of $\mu$,
accumulating at $\mu =0$, so that the corresponding homoclinic
tangles of equation (\ref{f2-s1}) are reduced to one single
horseshoe of infinitely many symbols;

(b) in between each of these parameter intervals, there are values
of $\mu$ so that the homoclinic tangles of equation (\ref{f2-s1})
contain stable periodic solutions; and

(c) there are also parameters in between where the homoclinic
tangles admit non-degenerate transversal homoclinic tangency.
\end{theorem}
\noindent {\bf Proof:} \ We argue that our previous proofs for
Theorems \ref{th1}-\ref{th3} remain valid for the current setups.
By assuming (H3), (\ref{f2-s4.2}) and $\omega \in [\frac{1}{000}
R_{\omega}, R_{\omega}]$, which is $>>1$, we make the first order
term for ${\mathcal Q}(t)$ dominate in $\phi_L(\theta)$. Let us
recall that
$$
\phi_L(\theta) = \int_{-L^-}^{L^+} {\mathbb H}(s) {\mathcal
Q}(\theta + \omega s + \omega L^-) e^{- \int_{0}^{s} E(\tau)
d\tau} ds
$$
is a critical element of the return map ${\mathcal F}_p$ in
Proposition \ref{prop1-s2.3c}. By definition
\begin{equation*}
\begin{split}
\phi(\theta) & = \int_{-\infty}^{\infty} {\mathbb H}(s) {\mathcal
Q}(\theta + \omega s + \omega L^-) e^{- \int_{0}^{s} E(\tau)
d\tau} ds \\
& = \sum_{n=1}^{\infty} \sqrt{c_n^2 + s_n^2}\cdot \sqrt{C^2(n
\omega) + S^2(n \omega)} \cdot \sin(n \theta - n \omega L^- -
\theta_n)
\end{split}
\end{equation*}
where $\theta_n$ are constants completely determined by $c_n, s_n,
C(n \omega), S(n \omega)$. We re-write ${\mathcal F}_p$ of
Proposition \ref{prop1-s2.3c} following the steps of Sect.
\ref{s3.1}, using $z$ for ${\mathbb X}$ and denoting $(\theta_1,
z_1) = {\mathcal F}_p(\theta, z)$ for $(\theta, z) \in \Sigma^-$.
From Proposition \ref{prop1-s2.3c} we have
\begin{equation}\label{f3-s4.2}
\begin{split}
\theta_1 & = \theta + {\bf a} - \frac{\omega}{\beta} \ln {\mathbb
F}(\theta, z, \mu) \\
z_1 & =  {\bf b} [{\mathbb F}(\theta, z,
\mu)]^{\frac{\alpha}{\beta}}
\end{split}
\end{equation}
where
\begin{equation} \label{f4-s4.2}
\begin{split}
{\bf a} & = \frac{\omega}{\beta} \ln \mu^{-1} + \omega (L^+  +
L^-) + \frac{\omega}{\beta} \ln(\varepsilon (1 +
{\mathcal O}(\varepsilon))P_L^+ A_L \rho) \\
{\bf b} & = (\mu \varepsilon^{-1})^{\frac{\alpha}{\beta}-1} [(1 +
{\mathcal O}(\varepsilon))P_L^+ A_L \rho]^{\frac{\alpha}{\beta}}; \\
\end{split}
\end{equation}
\begin{equation} \label{f5-s4.2}
{\mathbb F}(\theta, z, \mu) = 1 + {\bf c} (\sin \theta +
\Phi(\theta)) + {\bf k} z + {\mathbb E}(\theta, \mu) + {\mathcal
O}_{\theta, z, p}(\mu)
\end{equation}
with
\begin{equation} \label{f6-s4.2}
\begin{split}
{\bf c} & = (A_L \rho)^{-1} \cdot \sqrt{c_1^2 + s_1^2} \cdot
\sqrt{C^2(\omega) + S^2(\omega)} \\
{\bf k} & = (A_L \rho)^{-1} (P_L P_L^+)^{-1} (1+ {\mathcal
O}(\varepsilon))
\end{split}
\end{equation}
and
\begin{equation} \label{f7-s4.2}
\begin{split}
{\mathbb E}(\theta, \mu) & = (A_L \rho)^{-1} (P_L^+)^{-1}(1+P_L)
{\mathcal O}_{\theta, p}(1) + (A_L \rho)^{-1} (\phi(\theta) -
\phi_L(\theta)) \\
\Phi(\theta) & = \sum_{n=2}^{\infty} \sqrt{\frac{c_n^2 +
s_n^2}{c_1^2 + s_1^2}}\cdot \sqrt{\frac{C^2(n \omega) + S^2(n
\omega)}{C^2(\omega) + S^2(\omega)}} \cdot \sin(n \theta - n
\omega L^- - \theta_n).
\end{split}
\end{equation}

By (H3) and the assumption that $\omega > \frac{1}{100}
R_{\omega}>>1$, $\Phi(\theta)$ is an added error term, toward
which our previous proofs of Theorems \ref{th1}-\ref{th3} are
indifferent. \hfill $\square$

\smallskip

From (\ref{f3-s4.2})-(\ref{f7-s4.2}) for ${\mathcal F}_p$ we see
that (\ref{f1-s3.1d}) is a prototype of return maps for all
${\mathcal Q}(\omega t)$ provided that $\omega >>1$. If the
forcing frequency is lower, then $\Phi(\theta)$ in (\ref{f7-s4.2})
remains an important part of $\phi_L(\theta)$ in ${\mathcal F}_p$.
It is then possible to have a number of disjoint vertical strips
for $V$, and an equal number of vertical strips for $U$. The
images of each of the $V$-components again wrap around $\Sigma^-$
infinitely many times in $\theta$ direction. We might, however,
have more turns for ${\mathcal F}(V)$, as shown in Fig. 12.

\begin{picture}(7, 5)
\put(3.7,0){ \psfig{figure=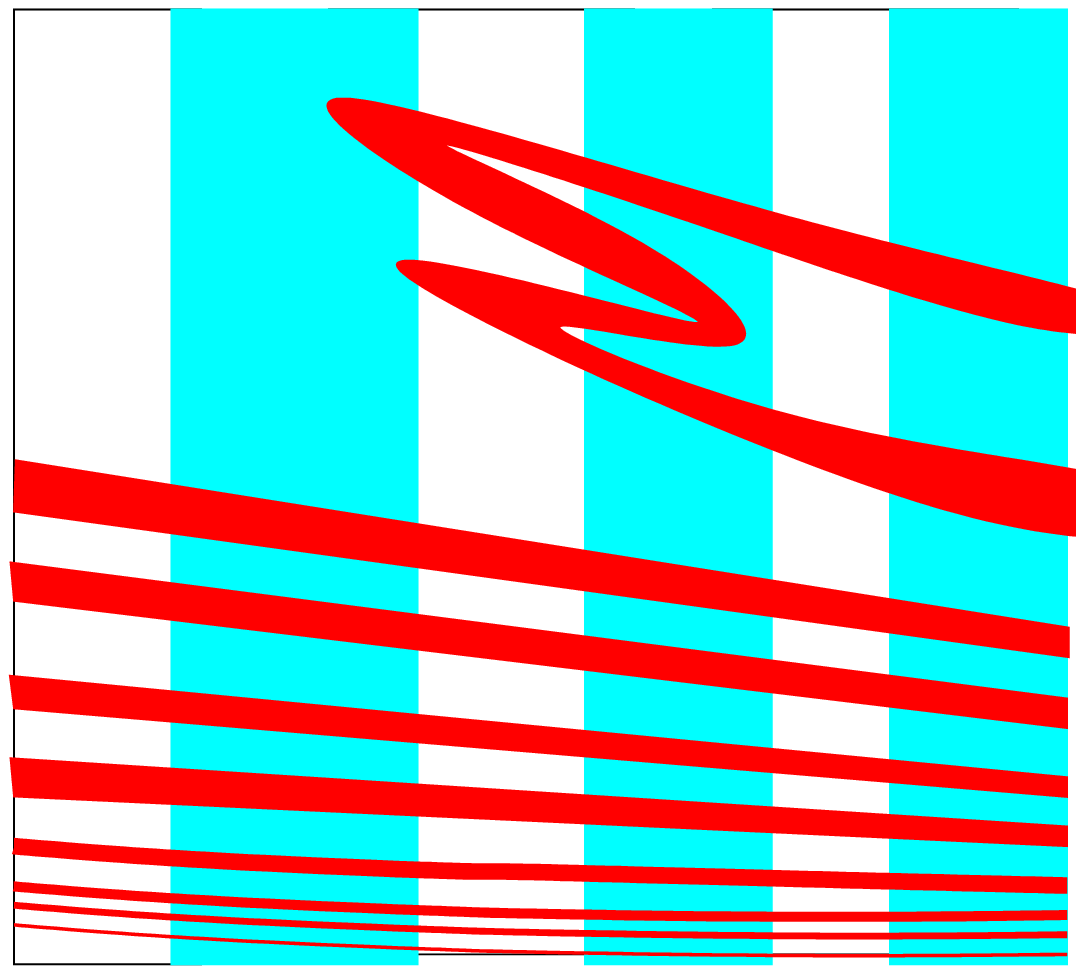,height = 4.5cm, width
= 7cm} }
\end{picture}

\vskip .1in

\centerline{Fig. 12 \ Infinitely wrapped horseshoe maps for
${\mathcal Q}(t)$ in general.}

\vskip .2in

We finish by presenting two more numerical pictures for $\Omega$
and $\Lambda$. These are for the maps assuming the form of
(\ref{f1-s3.1d}), but with $\sin \theta$ replaced by $\sin \theta
+ \sin 3 \theta$. Fig. 13 is for the case of an attracting
periodic sink, with $a = 1$, $b = 0.005, c = 1, d = 2$. $\Omega$
is depicted in Fig. 13(a). $\theta_k$ v.s. $k$ for one orbit from
$\Omega$ is depicted in Fig. 13(b). All orbits in $\Omega$ is
attracted to a periodic sink.

Fig. 14 is for a strange attractor with an SRB measure. The values
for $b$, $c, d$ and $\gamma$ are kept the same as in Fig. 13, but
$a$ is changed to $0.5$. $\Omega$ is depicted in Fig. 14(a), and
$\theta_k$ v.s. $k$ for one orbit from $\Omega$ is depicted in Fig
14(b). This orbit is attracted to an SRB measure.

\begin{picture}(5, 5.5)
\put(1,0){ \psfig{figure=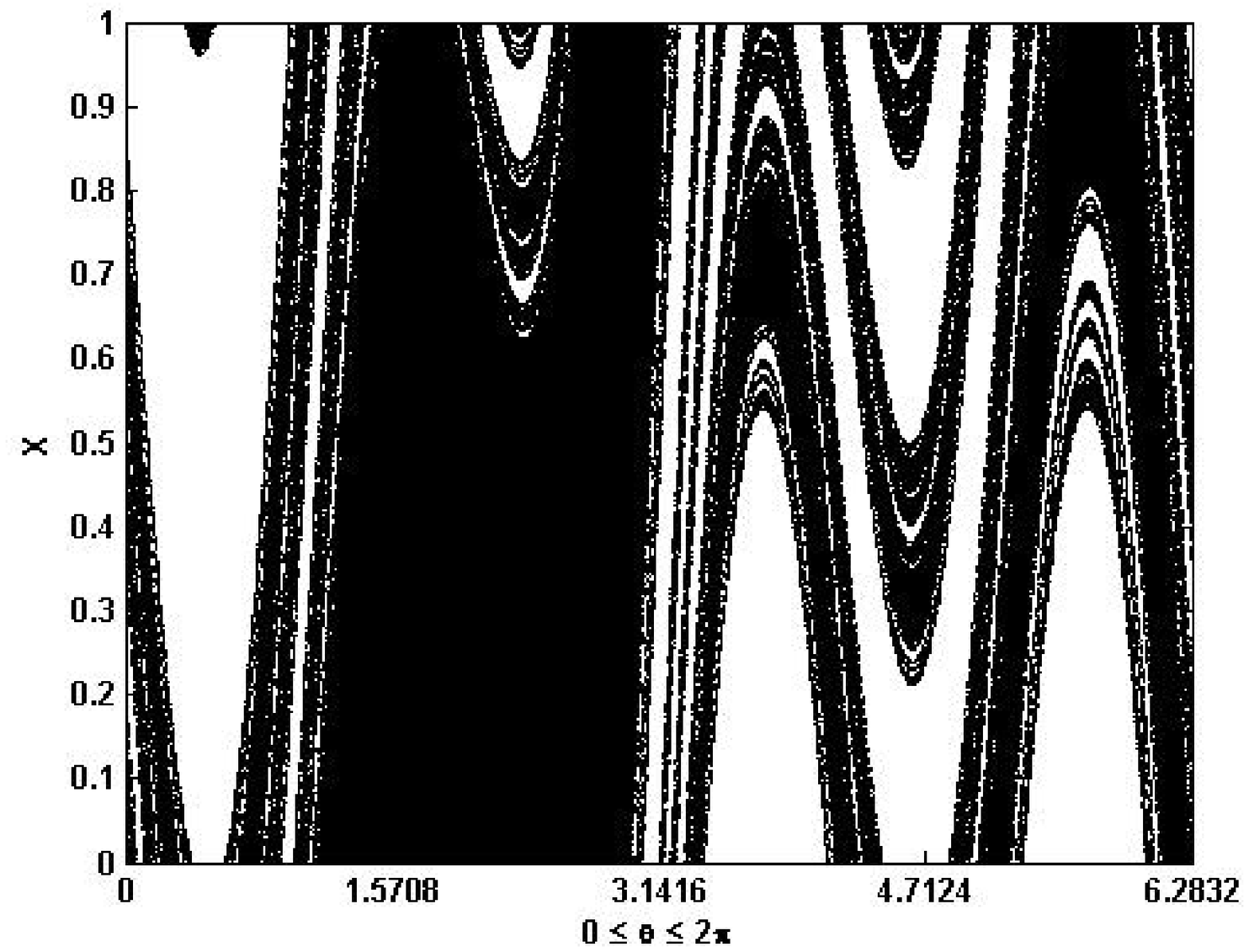,height = 5cm, width =
5cm} }
\end{picture}
\begin{picture}(9, 5.6)
\put(3, 0){ \psfig{figure=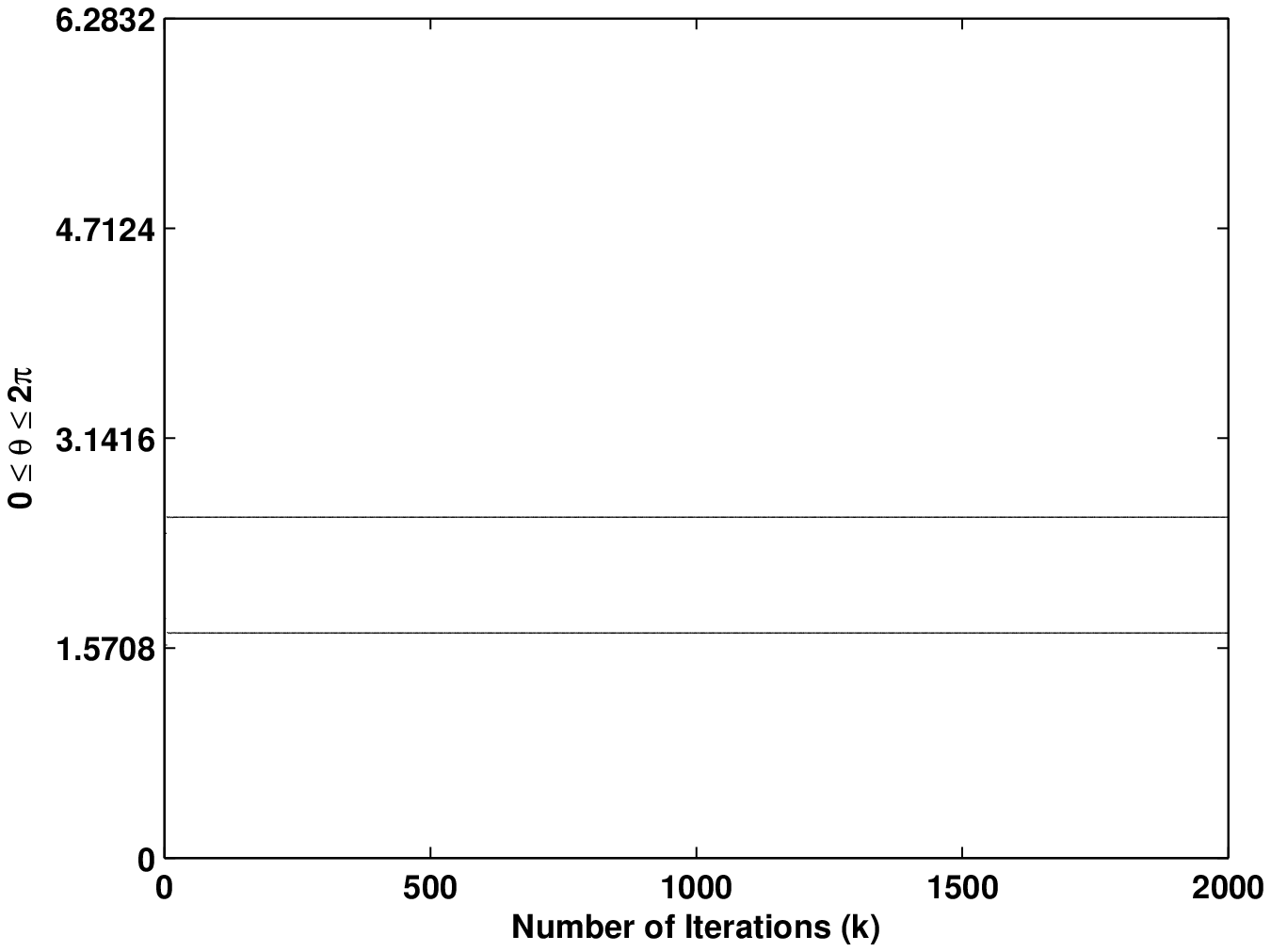,height = 5cm, width =
5cm} }
\end{picture}

\hspace{3.5cm} (a) \hspace{6cm} (b)

\vskip .1in

\centerline{Fig. 13 \ Homoclinic tangle with a periodic sink.}

\centerline{($a = 1, b = 0.005, c = 1, d = 2$ and $\gamma =
\sqrt{2}$)}

\begin{picture}(5, 5.5)
\put(1,0){ \psfig{figure=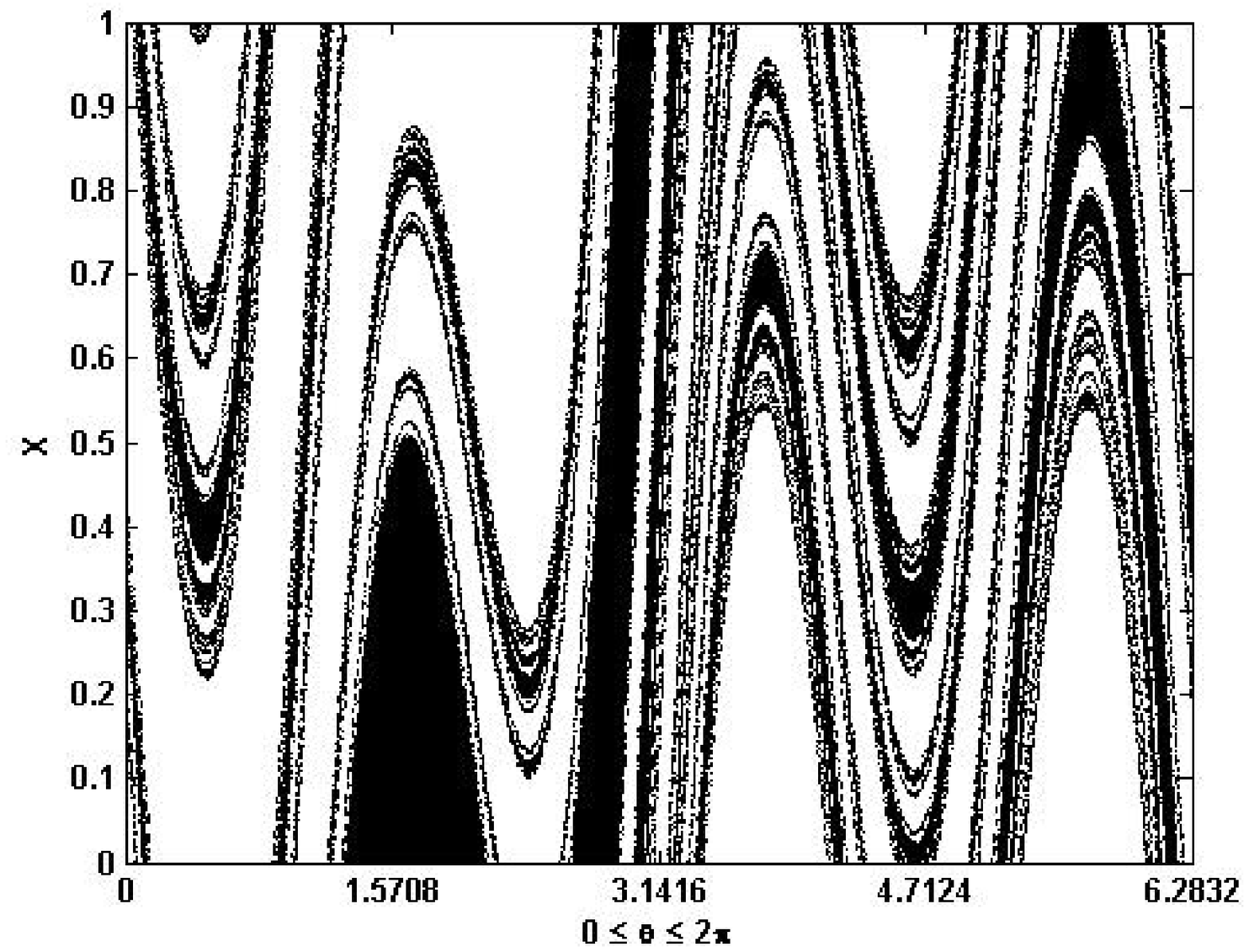,height = 5cm, width =
4.5cm} }
\end{picture}
\begin{picture}(9, 5.5)
\put(3, 0){ \psfig{figure=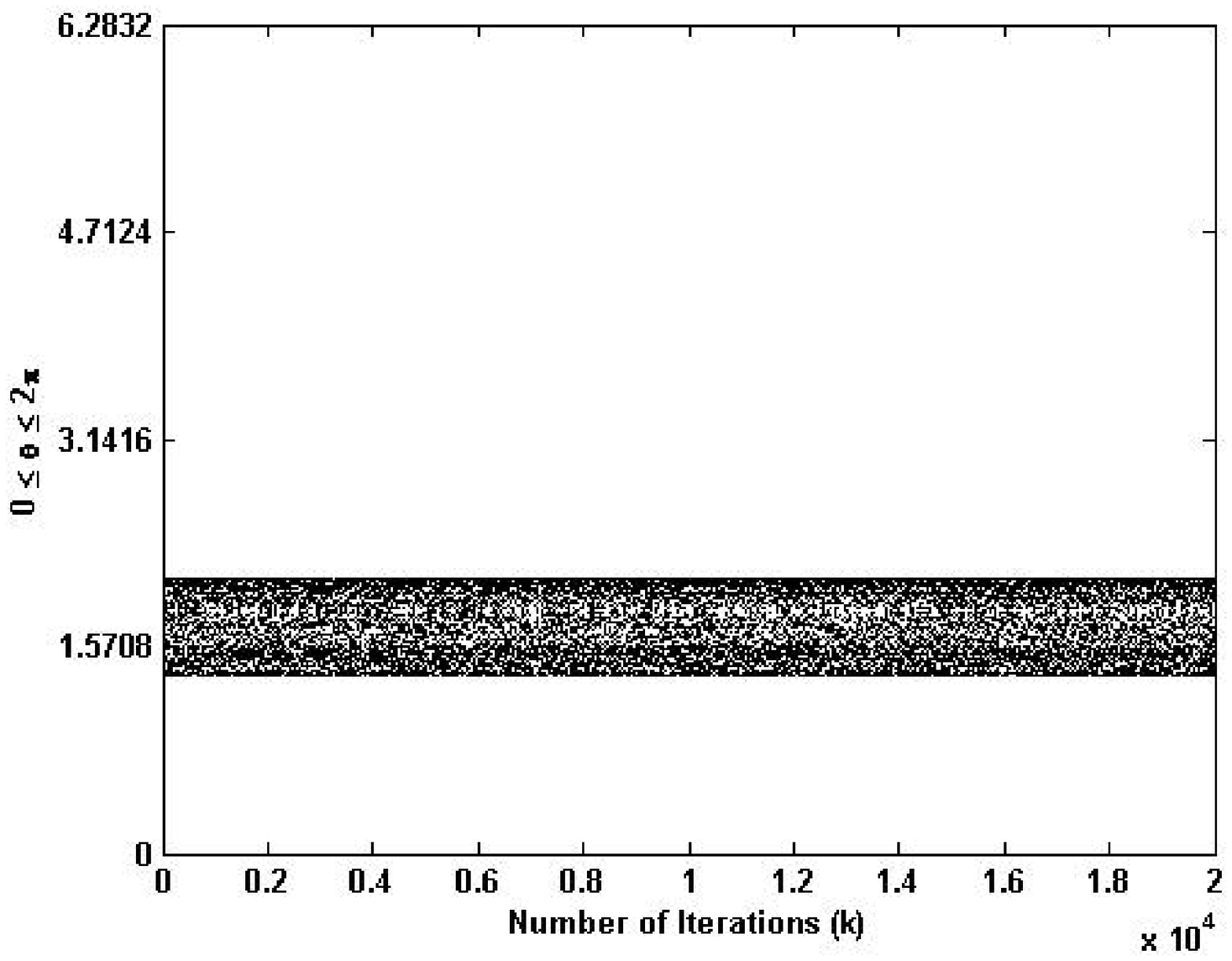,height = 5cm, width =
5cm} }
\end{picture}

\hspace{3.5cm} (a) \hspace{6cm} (b)

\centerline{Fig. 14 \ Homoclinic tangles with observable chaos.}

\centerline{($a = 0.5, b = 0.005, c = 1, d = 2$ and $\gamma =
\sqrt{2}$)}

\bigskip

\appendix

\section{Proof of Theorem \ref{th3}}

In this appendix we prove Theorem \ref{th3}.\footnote{Minus Claim
\ref{claim6-th3}(b), which we prove in Appendix B.} Let ${\bf
a}_n$ be the value of ${\bf a}$ at $\mu = \mu^{(r)}_n$ and $[{\bf
a}_n] = {\bf a}_n - {\bf a}_n mod (2 \pi)$. Let ${\bf a}(\mu)$ be
the value of ${\bf a}$ at $\mu \in [\mu^{(r)}_n,
\mu^{(l)}_{n+1}]$. We divide the proof of this theorem into the
following steps.

\medskip

\noindent {\it Step 1. Solving for hyperbolic fixed points} \ For
$\mu \in [\mu^{(r)}_n, \mu^{(l)}_{n+1}]$, let $m$ be an integer
$\geq 3\omega \beta^{-1}$ and $q_m({\bf a}) = (\theta_m, z_m)$ be
the solution of the equations
\begin{equation}\label{f1-th3}
\begin{split}
\theta + [{\bf a}_n] + 2 \pi m & = \theta + {\bf a}(\mu) -
\frac{\omega}{\beta} \ln {\mathbb
F}(\theta, z, \mu) \\
z & =  {\bf b} [{\mathbb F}(\theta, z,
\mu)]^{\frac{\alpha}{\beta}}.
\end{split}
\end{equation}
$\theta_m$ is determined by
\begin{equation}\label{f3-th3}
{\mathbb F}(\theta_m, z_m, \mu) = e^{\omega^{-1} \beta ( {\bf
a}(\mu) -[{\bf a}_n] - 2 \pi m)},
\end{equation}
and
\begin{equation}\label{f2-th3}
z_m = {\bf b} e^{\omega^{-1} \alpha ({\bf a}(\mu) - [{\bf a}_n]- 2
\pi m)}.
\end{equation}

\begin{claim}\label{claim1-th3}
$q_m({\bf a}) = (\theta_m, z_m)$ is saddle fixed point.
\end{claim}
\noindent {\it Proof of Claim \ref{claim1-th3}:}\ Recall that
\begin{equation*}
\begin{split}
Tr(D{\mathcal F}) & = \frac{\partial \theta_1}{\partial \theta} +
\alpha \beta^{-1} {\bf b} {\mathbb F}^{\alpha \beta^{-1} -1}
\frac{\partial {\mathbb
F}}{\partial z} \\
\det(D{\mathcal F}) & = \alpha \beta^{-1} {\bf b} {\mathbb
F}^{\alpha \beta^{-1} -1} \frac{\partial {\mathbb F}}{\partial z}.
\end{split}
\end{equation*}
Observe that, from (\ref{f3-th3}) and the assumption that $m \geq
3 \omega \beta^{-1}$, ${\mathbb F}(\theta_m, z_m, \mu) <
\frac{1}{100}$. It follows that $|{\bf c} \cos \theta_m| > 1$, and
$$
\left| \frac{\partial \theta_1}{\partial \theta} \right| > 101.
$$
This implies
$$
|Tr(D{\mathcal F})| > 100.
$$
Observe that we also have $\det(D{\mathcal F}) << 1$. Therefore we
have two eigenvalues, one is close to $0$ and the other is with
magnitude $> 1$. \hfill $\Diamond$

\smallskip

We also have
\begin{claim}\label{claim2-th3}
For $m \geq 3 \omega \beta^{-1}$,
$$
\left| \frac{d \theta_m}{d {\bf a}} \right| < \frac{1}{100}, \ \ \
\ \left|\frac{d z_m}{d {\bf a}} \right| < K {\bf b}.
$$
\end{claim}
\noindent {\it Proof of Claim \ref{claim2-th3}:} \ Estimate for
$\frac{d z_m}{d{\bf a}}$ follows directly from (\ref{f2-th3}). To
estimate $\frac{d \theta_m}{d {\bf a}}$ we take derivative with
respect to ${\bf a}$ on both-side of (\ref{f3-th3}) and use
${\mathbb F} < \frac{1}{100}$ to obtain $|{\bf c} \cos \theta_m|
>1$.
\hfill $\Diamond$

\medskip

\noindent {\it Step 2. The stable and the unstable manifold for
$q_m$} \  In the rest of this proof we let $m$ be the smallest
integer $> 3 \omega \beta^{-1}$. Denote $q({\bf a}) = q_m({\bf
a})$. Let
$$
\hat V = \{(\theta, z) \in V, \ {\mathbb F}
> {\mathbb F}(q_{100 m}({\bf a}), \mu) \}
$$
We obtain $\hat V$ from $V$ by taking away two thin vertical
strips at the vertical boundaries of $V$. Observe that, by
definition, $q({\bf a}) \in \hat V$. We make $\varepsilon$
sufficiently small so that (1) the distance from $q({\bf a})$ to
the vertical boundary of $\hat V$ is $>>{\bf k}$; and (2) $K_m: =
{\mathbb F}(q_{100 m}) >> {\bf k}$.

Denote the stable and the unstable manifold of $q = q({\bf a})$ as
$W^s(q)$ and $W^u(q)$ respectively. The local stable and the local
unstable manifold are denoted as $W_{loc}^s(q)$ and
$W_{loc}^u(q)$. Let $\ell^u(q)$ be the connected branch of
$W^u(q)$ in $\hat V \setminus V_f$ that contains $W_{loc}^u(q)$,
and $\ell_1^u(q) = {\mathcal F}_{\bf a}(\ell^u(q))$. $\ell_1^u(q)$
is a horizontal curve traversing $V_f$ in $\theta$-direction: it
is straight forward to verify that (1) $\ell^u(q)$ is a horizontal
curve, (2) the image of $\ell^u(q)$ is also a horizontal curve,
and (3) the length of that image at least doubles the length of
$\ell^u(q)$ so it traverses $V_f$. Let $z = w^u(\theta)$ be such
that $(\theta, w^u(\theta)) \in \ell_1^u(q)$.
\begin{claim}\label{claim3-th3}
We have on $\ell_1^u(q)$,

(a) $\left| \frac{d w^u}{d \theta} \right| < {\bf
b}^{\frac{1}{2}}$;

(b) $\left| \frac{d^2 w^u}{d \theta^2} \right| < {\bf
b}^{\frac{1}{2}}$.
\end{claim}

\noindent {\it Proof of Claim \ref{claim3-th3}:} \ Denote
${\mathcal F} = {\mathcal F}_{\bf a}$. For $(\theta, z) \in
\ell^u(q)$, let $(\theta_1, z_1) = {\mathcal F}(\theta, z)$. We
have from (\ref{f1-add-s3.2})
\begin{equation}\label{f4-th3}
\left| \frac{d w^u(\theta_1)}{d \theta_1} \right| =
\left|\frac{\alpha \beta^{-1}  {\bf b} {\mathbb F}^{\alpha
\beta^{-1}-1} \frac{\partial {\mathbb F}}{\partial \theta} +\alpha
\beta^{-1} {\bf b} {\mathbb F}^{\alpha \beta^{-1}-1}
\frac{\partial {\mathbb F}}{\partial z} \frac{d w^u(\theta)}{d
\theta} }{(1- \omega \beta^{-1} \frac{1}{{\mathbb F}}
\frac{\partial {\mathbb F}}{\partial \theta}) +  \omega \beta^{-1}
\frac{1}{{\mathbb F}} \frac{\partial {\mathbb F}}{\partial z}
\frac{d w^u(\theta)}{d \theta}}\right|.
\end{equation}
(a) holds because the magnitude of the denominator in
(\ref{f4-th3}) is $> 1$ for $(\theta, z) \in \hat V \setminus
V_f$. Remember that since $\ell^u(q)$ is horizontal we have
$|\frac{d w^u(\theta)}{d \theta}| < \frac{1}{100}$. To prove (b)
we take derivative one more time to obtain
\begin{equation}\label{f5-th3}
\frac{d^2 w^u(\theta_1)}{d^2 \theta_1} =\frac{\frac{d}{d
\theta}\left(\frac{d w^u(\theta_1)}{d \theta_1}\right)}{\frac{d
\theta_1}{d \theta}}.
\end{equation}
Observe that $\frac{d \theta_1}{d \theta}$ is the denominator in
(\ref{f4-th3}), the magnitude of which is $> 1$. Let
$$
M = \max_{(\theta, z) \in \ell_1^u}\left|\frac{d^2 w^u(\theta)}{d
\theta^2}\right|.
$$
We have from (\ref{f4-th3}) and (\ref{f5-th3})
$$
\left| \frac{d^2 z_1}{d^2 \theta_1} \right| < K_1 {\bf b} + K_2
{\bf b} M.
$$
So
$$
M < K_1 {\bf b} + K_2 {\bf b} M,
$$
and $ M < K {\bf b} < {\bf b}^{\frac{1}{2}}$. \hfill $\Diamond$

\smallskip

Let $\ell^s_{-1}(q)$ be the segment of $W^s(q)$ in $\hat V$ that
contains $W_{loc}^s(q)$, and $\ell^s(q) = {\mathcal
F}(\ell^s_{-1}(q))$. $\ell_{-1}^s(q)$ is a fully extended vertical
curve in $\hat V$, which we represent by a function $\theta =
w^s(z)$.

\begin{claim}\label{claim4-th3}
We have on $\ell^s_{-1}(q)$,

(a) $\left| \frac{d w^s(z)}{d z} \right| < {\bf k}^{\frac{1}{2}}$;
and

(b) $\left| \frac{d^2 w^s(z)}{d z^2} \right| < {\bf
b}^{\frac{1}{2}}$.
\end{claim}
\noindent {\it Proof of Claim \ref{claim4-th3}:} \ Let $(\theta,
z) \in \ell^s(q)$ and denote $(\theta_{-1}, z_{-1}) = {\mathcal
F}^{-1}(\theta, z)$. We have from (\ref{f4-add-s3.2})
\begin{equation}\label{f6-th3}
\left( \begin{array}{c} \frac{d \theta_{-1}}{d z} \\ \frac{d
z_{-1}}{d z} \end{array} \right) = \frac{1}{\alpha \beta^{-1} {\bf
b} {\mathbb F}^{\alpha \beta^{-1} -1} \frac{\partial {\mathbb
F}}{\partial z}} \left(
\begin{array}{cc} \alpha \beta^{-1}
{\bf b} {\mathbb F}^{\alpha \beta^{-1} -1} \frac{\partial {\mathbb
F}}{\partial z} &  - \omega \beta^{-1} \frac{1}{{\mathbb F}}
\frac{\partial {\mathbb F}}{\partial z} \\
- \alpha \beta^{-1}  {\bf b} {\mathbb F}^{\alpha \beta^{-1}-1}
\frac{\partial {\mathbb F}}{\partial \theta} & 1 - \omega
\beta^{-1} \frac{1}{{\mathbb F}} \frac{\partial {\mathbb
F}}{\partial \theta}
\end{array} \right) \left(
\begin{array}{c} \frac{d w^s(z)}{d z} \\ 1
\end{array} \right),
\end{equation}
and it follows that
\begin{equation}\label{f7-th3}
\left|\frac{d w^s(z_{-1})}{d z_{-1}} \right| = \left|\frac{\alpha
\beta^{-1} {\bf b} {\mathbb F}^{\alpha \beta^{-1} -1}
\frac{\partial {\mathbb F}}{\partial z} \frac{d w^s(z)}{d z}-
\omega \beta^{-1} \frac{1}{{\mathbb F}} \frac{\partial {\mathbb
F}}{\partial z}}{- \alpha \beta^{-1} {\bf b} {\mathbb F}^{\alpha
\beta^{-1}-1} \frac{\partial {\mathbb F}}{\partial \theta} \frac{d
w^s(z)}{d z} + (1 - \omega \beta^{-1} \frac{1}{{\mathbb F}}
\frac{\partial {\mathbb F}}{\partial \theta})} \right|.
\end{equation}
Let
$$
M_1 = \max_{(\theta, z) \in \ell_{-1}^s(q)} \left| \frac{d
w^s(z)}{d z} \right|.
$$
We have from (\ref{f7-th3}),
$$
\left|\frac{d w^s(z_{-1})}{d z_{-1}} \right| < K_1 {\bf k} + K_2 b
M_1
$$
because ${\mathbb F} > K_m$ and
$$
(1 - \omega \beta^{-1} \frac{1}{{\mathbb F}} \frac{\partial
{\mathbb F}}{\partial \theta}) > 2.
$$
It then follows that $M_1 < K {\bf k} < {\bf k}^{\frac{1}{2}}$.

To prove (b) we write
\begin{equation}\label{f8-th3}
\frac{d^2 w^s(z_{-1})}{d z_{-1}^2} = \frac{\frac{d}{d
z}\left(\frac{d w^s(z_{-1})}{d z_{-1}} \right)}{\frac{d
z_{-1}}{dz}},
\end{equation}
where $ \frac{d w^s(z_{-1})}{d z_{-1}}$ is as in (\ref{f7-th3})
and
\begin{equation}\label{f9-th3}
\frac{d z_{-1}}{d z} =  \frac{1}{\alpha \beta^{-1} {\bf b}
{\mathbb F}^{\alpha \beta^{-1} -1} \frac{\partial {\mathbb
F}}{\partial z}} \left(- \alpha \beta^{-1}  {\bf b} {\mathbb
F}^{\alpha \beta^{-1}-1} \frac{\partial {\mathbb F}}{\partial
\theta} \frac{d w^s(z)}{dz} + 1 - \omega \beta^{-1}
\frac{1}{{\mathbb F}} \frac{\partial {\mathbb F}}{\partial \theta}
\right).
\end{equation}
Let
$$
M_2 = \max_{(\theta, z) \in \ell_{-1}^s(q)} \left| \frac{d^2
w^s(z)}{d z^2} \right|.
$$
We have from (\ref{f7-th3}), (\ref{f8-th3}) and (\ref{f9-th3})
that
$$
\left|\frac{d^2 w^s(z_{-1})}{d z_{-1}^2}\right| < K {\bf b} (K_1
{\bf b} M_2 + K_2),
$$
from which we obtain $M_2 < K {\bf b} < {\bf b}^{\frac{1}{2}}$.
\hfill $\Diamond$

\medskip

\noindent {\it Step 3. \ Non-degenerate, transversal tangency} \
Let $\ell_{\bf a}^u$ be a connected segment of $\ell^u_1(q) \cap
V_f$, and $\ell^s_{\bf a}$ be the vertical curve $\ell_{-1}^s(q)$
where $\ell_1^u(q), \ \ell^s_{-1}(q)$ are as in Step 2. We use $z
= w^u(\theta)$ to represent $\ell_{\bf a}^u$ and $\theta = w^s(z)$
to represent $\ell_{\bf a}^s$. ${\mathcal F}_{\bf a}(\ell_{\bf
a}^u)$ traverses $\hat V$ in horizontal direction as $\mu$ runs
through $[\mu_n^{(r)}, \mu_{n+1}^{(l)}]$. Consequently there
exists $\hat \mu \in [\mu_n^{(r)}, \mu_{n+1}^{(l)}]$, the
corresponding value for ${\bf a}$ we denote as $\hat {\bf a}$, so
that $\ell_{\hat {\bf a}}^s$ and ${\mathcal F}_{\hat {\bf
a}}(\ell^u_{\hat {\bf a}})$ intersect tangentially at a point we
denote as $\tilde q= (\tilde \theta, \tilde z)$. Let $(\theta_0,
z_0) \in \ell_{\hat {\bf a}}^u$ be such that $(\tilde \theta,
\tilde z) = {\mathcal F}_{\hat {\bf a}}(\theta_0, z_0)$. Our next
claim implies that the tangential intersection of $\ell^s_{\hat
{\bf a}}$ and ${\mathcal F}_{\hat {\bf a}}(\ell^u_{\hat {\bf a}})$
at $\tilde q$ is not degenerate.
\begin{claim}\label{claim5-th3}
For $(\theta, z) \in \ell_{\hat {\bf a}}^u$, let $(\theta_1, z_1)
= {\mathcal F}_{\hat {\bf a}}(\theta, z)$. Then at $(\theta, z) =
(\theta_0, z_0)$, we have
$$
\left|\frac{d^2 \theta_1}{d z_1^2}\right|
>> 1.
$$
\end{claim}
\noindent {\it Proof of Claim \ref{claim5-th3}:} \ From
(\ref{f1-add-s3.2}) we have
\begin{equation}\label{f10-th3}
\frac{d \theta_1}{d z_1} = \frac{(1- \omega \beta^{-1}
\frac{1}{{\mathbb F}} \frac{\partial {\mathbb F}}{\partial
\theta}) +  \omega \beta^{-1} \frac{1}{{\mathbb F}} \frac{\partial
{\mathbb F}}{\partial z} \frac{d w^u(\theta)}{d \theta}}{\alpha
\beta^{-1}  {\bf b} {\mathbb F}^{\alpha \beta^{-1}-1}
\frac{\partial {\mathbb F}}{\partial \theta} +\alpha \beta^{-1}
{\bf b} {\mathbb F}^{\alpha \beta^{-1}-1} \frac{\partial {\mathbb
F}}{\partial z} \frac{d w^u(\theta)}{d \theta}}.
\end{equation}
At the point of tangential intersection, we have $\left|\frac{d
\theta_1}{d z_1}\right| < {\bf k}^{\frac{1}{2}}$, which is not
possible unless
\begin{equation}\label{f11-th3}
\left|1- \omega \beta^{-1} \frac{1}{{\mathbb F}} \frac{\partial
{\mathbb F}}{\partial \theta}\right| < {\bf b}^{\frac{1}{4}}
\end{equation}
from (\ref{f10-th3}). This is because $\frac{d w^u(\theta)}{d
\theta} < {\bf b}^{\frac{1}{2}}$ from Claim \ref{claim3-th3}(a).
The effect of ${\bf b}$ in the denominator can not be possibly
balanced if (\ref{f11-th3}) is false.

For the estimate on second derivative we start from
\begin{equation}\label{f12-th3}
\frac{d^2 \theta_1}{d z_1^2} = \frac{\frac{d}{d \theta}
\left(\frac{d \theta_1}{d z_1} \right)}{\frac{d z_1}{d \theta}}
\end{equation}
where $\frac{d z_1}{d \theta}$ is the denominator in
(\ref{f10-th3}). To compute $\frac{d}{d \theta} \left(\frac{d
\theta_1}{d z_1} \right)$, we take derivative of the function on
the right hand side of (\ref{f10-th3}) with respect to $\theta$.
Applying the quotient rule we obtain a fraction, the bottom of
which has a factor ${\bf b}^2$. On the top, we have a collection
of finitely many terms, each of which is $< K {\bf b}^{1 +
\frac{1}{4}}$ in magnitude except one in the form of
\begin{equation}\label{f13-th3}
\left(\alpha \beta^{-1}  {\bf b} {\mathbb F}^{\alpha \beta^{-1}-1}
\frac{\partial {\mathbb F}}{\partial \theta} \right)\frac{d}{d
\theta} \left(1- \omega \beta^{-1} \frac{1}{{\mathbb F}}
\frac{\partial {\mathbb F}}{\partial \theta}\right).
\end{equation}
Remember that we have ${\mathbb F} > K_m$ on $\hat V$, and
$\frac{d {\mathbb F}}{d \theta} > K^{-1}$ from (\ref{f11-th3}). We
also have
$$
\left|\frac{d}{d \theta} \left(1- \omega \beta^{-1}
\frac{1}{{\mathbb F}} \frac{\partial {\mathbb F}}{\partial
\theta}\right)\right| \approx \frac{\omega \beta^{-1}}{{\mathbb
F}^2}({\bf c}^2 + {\bf c} \sin \theta) > 1.
$$
Therefore, (\ref{f13-th3}) is the dominating term on top and we
obtain
$$
\left|\frac{d^2 \theta_1}{d z_1^2}\right| > K {\bf b}^{-2}
$$
at $\tilde q$. \hfill $\Diamond$

\smallskip

To finish our proof of Theorem \ref{th3}, we also need to prove
that, as ${\bf a}$ varies, $\ell_{\bf a}^s$ and ${\mathcal F}_{\bf
a}(\ell_{\bf a}^u)$ move with different speed at the point of
tangency. To make the dependency on parameter ${\bf a}$ explicit,
we write $w^u = w^u(\theta, {\bf a})$, $w^s = w^s(z, {\bf a})$.
Claim \ref{claim3-th3} applies to $w^u(\theta, {\bf a})$ and Claim
\ref{claim4-th3} applies to $w^s(z, {\bf a})$.

\begin{claim}\label{claim6-th3}
Let $(\theta_1(\theta_0, {\bf a}), z_1(\theta_0, {\bf a})) =
{\mathcal F}_{{\bf a}}(\theta_0, w^u_{\bf a}(\theta_0, {\bf a}))$.
Then at ${\bf a} = \hat {\bf a}$ we have

(a) $\left|\frac{\partial}{\partial {\bf a}} \theta_1(\theta_0,
{\bf a}) \right| > \frac{2}{3}$; and

(b) $\left| \frac{\partial}{\partial {\bf a}} w^s(\tilde z, {\bf
a})\right| < \frac{1}{25}$.

\noindent Recall that $\tilde q = (\tilde \theta, \tilde z)$ is
the point of tangential intersection and $(\theta_0, z_0)$ is such
that ${\mathcal F}_{\hat {\bf a}}(\theta_0, z_0)= \tilde q$.
\end{claim}
\noindent {\it Proof of Claim \ref{claim6-th3}:} \ In this prove
we use $\partial_z, \partial_{\theta}$ and $\partial_{\bf a}$ to
denote partial derivative with respect to $z, \theta$ and ${\bf
a}$ respectively.

To prove (a) we let $q = (\theta_m, z_m)$ be the saddle fixed
point and $\ell^u(q)$ and $\ell^u_1(q) = {\mathcal F}_{\bf
a}(\ell^u(q))$ be as in Claim \ref{claim3-th3}. For $(\theta, z)
\in \ell^u(q)$ and $(\theta_0, z_0) = {\mathcal F}_{\bf a}(\theta,
z)$. We have from (\ref{f1-s3.1a}),
\begin{equation}\label{f15-th3}
\begin{split}
\theta_0 & = f(\theta, z, {\bf a})= \theta + {\bf a} -
\frac{\omega}{\beta} \ln {\mathbb
F}(\theta, z, \mu) \\
z_0 & = g(\theta, z, {\bf a}) =  {\bf b} [{\mathbb F}(\theta, z,
\mu)]^{\frac{\alpha}{\beta}},
\end{split}
\end{equation}
and in (\ref{f15-th3}), $z_0 = w^u(\theta_0, {\bf a}), z =
w^u(\theta, {\bf a})$ because both $(\theta_0, z_0)$ and $(\theta,
z)$ are on $\ell^u_1(q)$. We first invert the first equality in
(\ref{f15-th3}), obtaining $\theta = \theta(\theta_0, {\bf a})$;
then we put it into the second equality in (\ref{f15-th3}) to
obtain $z_0 = w^u(\theta_0, {\bf a})$. To estimate $\partial_{\bf
a} w^u(\theta_0, {\bf a})$, we first let
$$
M_{\bf a} = \max_{(\theta, z) \in \ell^u_1(q)} |\partial_{\bf a}
w^u(\theta, {\bf a})|
$$
and obtain from the first equality in (\ref{f15-th3}),
$$
\left|\partial_{\bf a} \theta(\theta_0, {\bf a})\right| = \left|
\frac{\partial_{\bf a} f + \partial_{z} f \cdot \partial_{\bf a}
w^u}{\partial_{\theta} f + \partial_z f \cdot \partial_{\theta}
w^u}\right| < K_1 + K_2 M_{\bf a}
$$
because $|\partial_{\theta} f| > 1$ for $(\theta, z) \in \hat V
\setminus V_f$ and $|\partial_{\theta} w^u| < {\bf
b}^{\frac{1}{2}}$ from Claim \ref{claim3-th3}(a). From the second
equality in (\ref{f15-th3}) we have
\begin{equation*}
\begin{split}
|\partial_{\bf a} w^u(\theta_0, {\bf a})| & = |\partial_{\theta} g
\cdot
\partial_{\bf a} \theta +
\partial_z g \cdot (\partial_{\theta} w^u \cdot \partial_{\bf a}
\theta +
\partial_{\bf a} w^u) + \partial_{\bf a} g| \\
& \leq K_3{\bf b}(K_1 + K_2 M_{\bf a}) + K_4 {\bf b},
\end{split}
\end{equation*}
from which it follows that
\begin{equation}\label{f14-th3}
M_{\bf a} < {\bf b}^{\frac{1}{2}}.
\end{equation}
(a) now follows by taking $\partial_{\bf a}$ on
$$
\theta_1(\theta_0, {\bf a}) = \theta_0 + {\bf a} - \omega
\beta^{-1} \ln {\mathbb F}(\theta_0, w^u(\theta_0, {\bf a}), \mu)
$$
using (\ref{f14-th3}).

\medskip

Proof of (b) is more sophisticated than that of (a). We need to
study the stable manifold through the field of most contracted
directions, a method originally introduced in \cite{BC} and fully
developed in \cite{WY1} and \cite{WY2}. A detailed proof is
included in Appendix B. \hfill $\Diamond$

\smallskip

With Claim \ref{claim6-th3} we know that, as ${\bf a}$ varies
passing $\hat {\bf a}$, ${\mathcal F}_{\bf a}(\ell^u_{\bf a})$
crosses $\ell_{\bf a}^s$ transversally.

This finishes our proof of Theorem \ref{th3} owing that of Claim
\ref{claim6-th3}(b). \hfill $\square$

\section{Proof of Claim \ref{claim6-th3}(b)}
In order to produce the desired estimates in Claim
\ref{claim6-th3}(b), we need more precise controls on the stable
manifold of the saddle fixed point $q_m$. The main idea of our
proof, that is, to approximate the stable manifold by using the
integral curves of vector field defined by the most contracted
directions of the Jacobi matrix, was originated from \cite{BC},
and was fully developed in \cite{WY1} and \cite{WY2}. Here we only
need a specific version of the contents developed in the beginning
part of Section 3 in \cite{WY2}.

\subsection{Most contracted directions}

In what follows $u_1 \wedge u_2$ is the wedge product and $\langle
u_1, u_2 \rangle$ is the inner product for $u_1, u_2 \in {\mathbb
R}^2$.

Let $M$ be a $2 \times 2$ matrix and assume $M \neq cO$ where $O$
is orthogonal and $c \in \mathbb R$. Then there is a unit vector
$e$, uniquely defined up to a sign, that represents {\it the most
contracted direction} of $M$, i.e. $|Me| \leq |Mu|$ for all unit
vectors $u$. From standard linear algebra, we know $f = e^\perp$
is the most expanded direction, meaning $|Me^\perp | \geq |Mu|$
for all unit vectors $u$, and $Me \perp Me^\perp$. The numbers
$|Me|$ and $|Me^\perp|$ are the {\it singular values} of $M$.

Let $u \perp v $ be two unit vectors in ${\mathbb R}^2$. The
following formulas are results of elementary computations. First,
we write down the squares of the singular values of $M$:
\begin{equation}\label{f1-appb}
|Me|^2  =  \frac{1}{2 } (B - \sqrt{B^2 - 4C}) := \lambda , \ \
\quad |Mf|^2  =  \frac{1}{2 } (B + \sqrt{B^2 - 4C})
\end{equation}
where
\begin{equation}\label{f2-appb}
B = |M u |^2 + |M v |^2, \ \ \ \ C = |M u \wedge M v|^2.
\end{equation}
We write $e=\alpha_0 u+ \beta_0 v$, and solve for $|Me|= \sqrt
\lambda$ subject to $\alpha_0^2+ \beta_0^2=1$. There are two
solutions (a vector and its negative): either $e=\pm v$, or the
solution with a positive $u$-component is given by
\begin{equation}\label{f3-appb}
e = \frac{1}{Z}(\alpha u + \beta v)
\end{equation}
with
\begin{equation}\label{f4-appb}
 \alpha= |M v|^2 - \lambda, \ \ \ \ \beta = - \langle M v, M u \rangle
\end{equation}
and
\begin{equation}\label{f5-appb}
Z= \sqrt {\alpha^2+ \beta^2}.
\end{equation}
From this we deduce that a solution for $f$ is
\begin{equation}\label{f6-appb}
f = \frac{1}{Z}(- \beta u + \alpha v).
\end{equation}

\subsection{Stability of most contracted directions}
In what follows we let $q_i = {\mathcal F}^i_{\bf a}(q_0)$, $M_i =
D{\mathcal F}_{\bf a}({q_{i-1}})$;
$$
M_i  =  \left( \begin{array}{cc} A_i &  B_i \\
C_i & D_i \end{array} \right) = \left(
\begin{array}{cc} 1
- \omega \beta^{-1} \frac{1}{{\mathbb F}} \frac{\partial {\mathbb
F}}{\partial \theta} &  \omega \beta^{-1} \frac{1}{{\mathbb F}}
\frac{\partial {\mathbb F}}{\partial z} \\
\alpha \beta^{-1}  {\bf b} {\mathbb F}^{\alpha \beta^{-1}-1}
\frac{\partial {\mathbb F}}{\partial \theta} & \alpha \beta^{-1}
{\bf b} {\mathbb F}^{\alpha \beta^{-1} -1} \frac{\partial {\mathbb
F}}{\partial z}
\end{array} \right).
$$
We have for $q_{i-1} \in \hat V \setminus V_f$,
\begin{equation}\label{f7-appb}
2 < |A_i| < K, \ \ |B_i|< K, \ \ |C_i|, \ |D_i| < K {\bf b}.
\end{equation}

Let $M^{(n)} = D{\mathcal F}_{\bf a}^n(q_0)$. $M^{(n)} = M_n \cdot
M_{n-1} \cdots M_1$. Let the most contracted direction for
$M^{(n)}$ be $e_n$ and the most expanded direction be $f_n$.
Denote the values of $\alpha, \beta$ and $Z$ in (\ref{f4-appb})
and (\ref{f5-appb}) for $M^{(n)}$ as $\alpha_n, \beta_n$ and
$Z_n$. Observe that, assuming $q_i \in \hat V \setminus V_f, \ i <
n$,
\begin{equation}\label{f8-appb}
|M^{(n)}f_n| > 1.
\end{equation}

We have
\begin{lemma}\label{lemma-a1}
Let $q_0$ be such that $q_0, \cdots, q_n \in \hat V \setminus
V_f$. Then for all $1 \leq i \leq n$,
\begin{itemize}
\item[(a)] $|e_{i+1} - e_{i}| < (K{\bf b})^i,  \ \ |M^{(i)} e_n| <
(K{\bf b})^i$; \item[(b)] $|\partial_{\bf a} (e_{i+1} - e_{i})| <
(K{\bf b})^i, \ \ |\partial_{\bf a} M^{(i)} e_n| < (K{\bf
b})^{i}$.
\end{itemize}
\end{lemma}
\noindent {\bf Proof:} \ Let $\Delta_i:= |M^{(i)} u \wedge M^{(i)}
v|$. We have
\begin{equation}\label{f9-appb}
\Delta_i = |\det(M^{(i)})| < (K{\bf b})^i.
\end{equation}
It then follows from $|M^{(i)}e_i| |M^{(i)}f_i|= \Delta_i$ and
(\ref{f8-appb}),
\begin{equation}\label{f10-appb}
|M^{(i)} e_i| < (K {\bf b})^i.
\end{equation}

We substitute $u = e_i, v = f_i$ and $M=M^{(i+1)}$ into
(\ref{f3-appb}) for $e_{i+1}$ and (\ref{f5-appb}) for $f_{i+1}$.
By using (\ref{f1-appb}) for $M^{(i+1)} f_{i+1}$, we have
\begin{equation}\label{f11-appb}
|M^{(i+1)} f_{i}| = |M^{(i+1)} f_{i+1}| \pm {\mathcal O}((K {\bf
b})^{i}).
\end{equation}
from (\ref{f2-appb}), (\ref{f9-appb}) and (\ref{f10-appb}). We
also have
\begin{equation}\label{f12-appb}
Z_{i+1} \approx |\alpha_{i+1}| \approx |M^{(i+1)} f_{i}|^2.
\end{equation}

We now prove Lemma \ref{lemma-a1}(a). Using $u=e_i$ and $v=f_i$,
we have, from (\ref{f3-appb}),
\begin{equation}\label{f13-appb}
e_{i+1}- e_{i} =
\frac{1}{Z_{i+1}}\left(\frac{-\beta_{i+1}^2}{\alpha_{i+1} +
Z_{i+1}} e_{i} + \beta_{i+1} f_{i} \right).
\end{equation}
To estimate $|e_{i+1} - e_{i}|$, we need to obtain a suitable
upper bound for $|\beta_{i+1}|$ and lower bounds for
$|\alpha_{i+1}|$ and $Z_{i+1}$. We have from (\ref{f4-appb}),
(\ref{f10-appb}) and (\ref{f12-appb}),
\begin{equation}\label{f14-appb}
|\beta_{i+1}| \leq |M^{(i+1)} e_{i}| |M^{(i+1)} f_{i}| < K{\bf
b}^{i} \sqrt {Z_{i+1}}
\end{equation}
and $|\alpha_{i+1}| \approx Z_{i+1}$. These estimates together
with $Z_{i+1} > 1$ tell us
$$
|e_{i+1} - e_{i}| \approx \ \frac{|\beta_{i+1}|}{Z_{i+1}} < (K
{\bf b})^{i}.
$$
The second assertion follows easily from
$$
|M^{(i)} e_n| \leq
|M^{(i)}(e_n - e_{n-1})| + \cdots + |M^{(i)}(e_{i+1} - e_i)| +
|M^{(i)} e_i| < (K {\bf b})^i.
$$
This finished our proof for Lemma \ref{lemma-a1}(a).

\bigskip

To prove Lemma \ref{lemma-a1}(b) we start with
\begin{sublemma}
\label{sublem-app5.2} $|\partial_{\bf a} e_1|, |\partial_{\bf a}
f_1| < K_1$ for some $K_1$.
\end{sublemma}
\noindent {\it Proof:} \ Let $u = (0, 1)^T$, $v = (1, 0)^T$ and
use (\ref{f3-appb}) for $e_1$ and (\ref{f6-appb}) for $f_1$. We
have $Z_1 > \alpha \geq |M_1 v|^2 - K{\bf b} > 1$. Differentiating
(\ref{f3-appb}) and (\ref{f6-appb}) gives the desired result.

\hfill $\Diamond$

\medskip

\noindent In the rest of this proof, $\partial = \partial_{\bf
a}$. Our plan of proof for Lemma \ref{lemma-a1}(b) is as follows:
For $k=1, 2, \cdots$, we assume for all $i \leq k$

\medskip
(*) $|\partial e_i|, |\partial f_i| < 2K_1$ where $K_1$ is as in
Sublemma \ref{sublem-app5.2},

\medskip
\noindent and prove for all $i \leq k$:

\medskip
(A) $|\partial (M^{(i)} f_i)| < K^i$, $|\partial (M^{(i)} e_i)| <
(K{\bf b})^i$;

\medskip
(B) $|\partial(e_{i+1}- e_{i})|, \ \  |\partial (f_{i+1}- f_{i})|
< (K{\bf b})^{i}$.

\bigskip
Observe that for $i=1$, (*) is given by Sublemma
\ref{sublem-app5.2}. It is easy to see that (B) above implies (*)
with $i=k+1$, namely $ |\partial f_{k+1}| \leq |\partial (f_{k+1}
- f_{k})| + \cdots + |\partial (f_2 - f_1)| + |\partial f_1|$.
From (B), we have $|\partial (f_{i+1} - f_{i})| < (K{\bf b})^{i}$,
and from Sublemma \ref{sublem-app5.2}, we have $|\partial f_1|<
K_1$. Hence $|\partial f_{k+1}| < K{\bf b} + K_1$, which, for
${\bf b}$ sufficiently small, is $< 2 K_1$. The computation for
$e_{k+1}$ is identical.

\bigskip
\noindent {\it Proof that (*)$\implies$(A):} \  First we prove the
estimate for $\partial (M^{(i)} f_{i})$. Writing
$$
\partial (M^{(i)} f_{i})  =  \sum_{j=1}^i M_i \cdots
(\partial M_j) \cdots M_1 f_i + M^{(i)} \partial f_i,
$$
we obtain easily
$$
|\partial (M^{(i)} f_{i})| \leq \sum_{j=1}^i | M_i \cdots
(\partial M_j) \cdots M_1 f_i| + \|M^{(i)}\| |\partial f_i| \leq i
K^i + K^i (2K_1).
$$

This estimate is used to estimate $\partial (M^{(i)} e_{i})$.
Write $\partial (M^{(i)} e_{i}) = (I)+(II)$ where $(I)$ is its
component in the direction of $M^{(i)} f_{i}$ and $(II)$ is its
component orthogonal to $M^{(i)} f_{i}$. Recall that $\partial
\langle M^{(i)} e_{i}, M^{(i)} f_{i}\rangle=0 $. We have
$$
|(I)| = \left| \langle \partial (M^{(i)} e_{i}), \frac{M^{(i)}
f_{i}}{|M^{(i)} f_{i}|} \rangle \right|
 =  \frac{1}{|M^{(i)} f_i|} \ | \langle M^{(i)} e_{i},
\partial (M^{(i)} f_{i})\rangle |
<  (K {\bf b})^i K^i\ ;
$$
$$
|(II)| \ |M^{(i)} f_{i}| = |\partial (M^{(i)} e_{i}) \wedge
M^{(i)} f_{i}| \leq |\partial (M^{(i)} e_{i} \wedge M^{(i)}
f_{i})| + |M^{(i)} e_{i} \wedge \partial (M^{(i)} f_{i})|.
$$
The first term in the last line is $<(K{\bf b})^i$, noting that we
have established $|\partial e_i|, |\partial f_i|< 2K_1$; the
second term is $< (K{\bf b})^i \cdot K^i$. This completes the
proof of (A). \hfill $\Diamond$

\bigskip
To prove (B), we first compute some quantities associated with the
next iterate. Substitute $u=e_i,v=f_i, M=M^{(i+1)}$ in
(\ref{f1-appb})-(\ref{f6-appb}). The following is a
straightforward computation.

\begin{sublemma}
\label{sublem-app5.3} Assume (*) and (A). Then for all $i \leq k$:

(a) $|\partial \lambda_{i+1}| <  (K {\bf b})^{2(i+1)}$;

(b) $|\partial \beta_{i+1}| < (K {\bf b})^i \sqrt {Z_{i+1}}$;

(c) $|\partial \alpha_{i+1}|, |\partial Z_{i+1}| < K^i \sqrt
{Z_{i+1}}$.
\end{sublemma}

\noindent {\it Proof that (*), (A)$\implies$(B):} \ We work with
$e_i$; the computation for $f_i$ is similar. From (23) we have $
\partial (e_{i+1} - e_{i}) = (III) + (IV) + (V)
$ where
\begin{eqnarray*}
|(III)| & = & |\frac{1}{Z_{i+1}} (e_{i+1} - e_{i}) \partial
Z_{i+1}| \ <  \ \frac{K^i \sqrt {Z_{i+1}}}{Z_{i+1}} \cdot (K {\bf
b})^i
<(K {\bf b})^{i};\\
|(IV)| & = & |\frac{1}{Z_{i+1}} \partial (\beta_{i+1} f_{i})| \ <
\ \frac{1}{Z_{i+1}}(|\partial \beta_{i+1}| +
|\beta_{i+1}||\partial f_{i}|) <  (K{\bf b})^{i};\\
|(V)| & = & |\frac{1}{Z_{i+1}}
\partial \left(\frac{\beta_{i+1}^2}{\alpha_{i+1} + Z_{i+1}} e_{i}\right)|
\ << \ (K {\bf b})^{i}.
\end{eqnarray*}
To estimate (III), we have used Sublemma \ref{sublem-app5.3}(c)
and part (a) of Lemma \ref{lemma-a1}. To estimate (IV), we have
used Sublemma \ref{sublem-app5.3}(b), (*) and
$|\beta_{i+1}|<(\frac{Kb}{\kappa})^i$. The estimate for (V) is
easy. \hfill $\Diamond$

\medskip

This completes the proof of Lemma \ref{lemma-a1}(b).   \hfill
$\square$

\bigskip

We also need to control the speed of change for the most
contracted directions in $\hat V \setminus V_f$. Let $q_0(s, {\bf
a})$ be a curve in $\hat V \setminus V_f$ parameterized by a
parameter $s$ and assume that
$$
\|q_0(s, {\bf a})\|_{C^2} < K.
$$
Let $M^{(n)}(s) =D {\mathcal F}^n_{\bf a}(q_0(s, {\bf a}))$, and
$e_n(s)$ be the most contracted direction for $M^{(n)}(s)$.
\begin{lemma}\label{lemma-a2}
Let $q_0$ be such that $q_0, \cdots, q_n \in \hat V \setminus
V_f$. Then for all $1 \leq i \leq n$,

(a) $ |\partial_s (e_{i+1}(s) - e_{i}(s))| < (K{\bf b})^i, \ \
|\partial_s M^{(i)}(s) e_n(s)| < (K{\bf b})^{i}$; and

(b) $|\partial_s \partial_{\bf a}(e_{i+1}(s) - e_{i}(s))| < (K{\bf
b})^i, \ \ |\partial_s \partial_{\bf a} M^{(i)}(s) e_n(s)| <
(K{\bf b})^{i}$.
\end{lemma}
\noindent {\bf Proof:} \ The proof for Lemma \ref{lemma-a2}(a) is
identical to that of Lemma \ref{lemma-a1}(b). It suffices to
regard all $\partial$ as $\partial_{s}$ instead of $\partial_{\bf
a}$. The estimate for the second derivatives is proved by a
similar argument. Here we skip the details. \hfill $\square$

\subsection{Temporary stable curves and the stable manifold}
In the rest of this proof we let $\eta = {\bf b}^{\frac{1}{10}}$
and denote ${\mathcal A}_{\eta} = \{(\theta, z) \in {\mathcal A}:
\ |z| < \eta \}$. We view $e_n$ as a vector field, defined where
it makes sense, and let $\gamma_n(s)$ be the integral curve to
$e_n$ with $\gamma_n(0)=q_0$.

\begin{lemma}\label{lemma-a3}
Let $q_0 = q_m$ be the saddle fixed point of Theorem \ref{th3},
and $\gamma_n(s)$ be the integral curve to $e_n$ satisfying
$\gamma_n(0) = q_0$ in ${\mathcal A}_{\eta}$. Then, for all $n >
0$,

(a) $|{\mathcal F}_{\bf a}^i(q) - {\mathcal F}_{\bf a}^i(q_0)| <
(Kb)^i|s|$ for all $q = \gamma_n(s)$ and all $i \leq n$;

(b) $\gamma_n(s)$ is a fully extended vertical curve in $(\hat V
\setminus V_f) \cap {\mathcal A}_{\eta}$;

(c) $|\gamma_{n+1}(s) - \gamma_n(s)|, \ |\partial_{\bf a}
\gamma_{n+1}(s) - \partial_{\bf a} \gamma_n(s) | < {\bf
b}^{\frac{n}{10}}$.
\end{lemma}

\noindent {\bf Proof:}  \ Lemma \ref{lemma-a3}(a) follows directly
from $|M^{(i)} e_{n}| < (K{\bf b})^i$ for all $i \leq n$ (Lemma
\ref{lemma-a1}(a)). Denote ${\mathcal F} = {\mathcal F}_{\bf a}$.
Let $B_0$ be the ball of radius $2 \eta$ centered at $q_0$. $e_1$
is well defined on $B_0$, and substituting $u = (0, 1)^T$, $v =
(1, 0)^T$ into (\ref{f3-appb}) we obtain $s(e_1)
> K {\bf b}^{-1}$. Let $\gamma_1 = \gamma_1(s)$ be the integral
curve to $e_1$ defined for $s \in (-2\eta, 2\eta)$ with
$\gamma_1(0)=q_0$.

To construct $\gamma_2$, let $B_1$ be the $\eta^2$-neighborhood of
$\gamma_1$. For $\xi \in B_1$, let $\xi'$ be a point in $\gamma_1$
with $|\xi-\xi'|< \eta^2 $. Then $|{\mathcal F}(\xi)-{\mathcal
F}(q_0)| \leq |{\mathcal F} (\xi)-{\mathcal F}(\xi')| +|{\mathcal
F}(\xi')-{\mathcal F}(q_0)| \leq K \eta^2 + K{\bf b} \eta < K
\eta^2$. This ensures that $e_2$ is defined on all of $B_1$. Let
$\gamma_2$ be the integral curve to $e_2$ with $\gamma_2(0)=q_0$.
We verify that $\gamma_2$ is defined on $(-2\eta, 2\eta)$ and runs
alongside $\gamma_1$. More precisely, let $t \in [0, 1]$ and
$$
q(t, s) = \gamma_1(s) + t(\gamma_2(s) - \gamma_1(s)).
$$
We have
\begin{eqnarray*}
|\frac{d}{ds} (\gamma_2(s)-\gamma_1(s)) | & \leq &
|e_2(\gamma_2(s)) - e_1(\gamma_2(s))| +
|e_1(\gamma_2(s)) - e_1(\gamma_1(s))| \\
& \leq & |e_2-e_1| + |\partial_t e_1| |\gamma_2(s) - \gamma_1(s)|
\\
& \leq & K {\bf b} + K |\gamma_2(s) - \gamma_1(s)|.
\end{eqnarray*}
Here we use $|\partial_t e_1| < K$. By Gronwall's inequality, $
|\gamma_2(s) - \gamma_1(s)| \leq K{\bf b} |s| e^{K |s|}, $ which
is $<< \eta^2 $ for $|s| < 2 \eta$. This ensures that $\gamma_2$
remains in $B_1$ and hence is well defined for all $s \in (-2
\eta, 2\eta)$.

In general, we inductively construct $\gamma_i$ by letting
$B_{i-1}$ be the $\eta^{i}$-neighborhood of $\gamma_{i-1}$ in $S$.
Then for all $\xi \in B_{i-1}$, $|{\mathcal F}^j (\xi) -{\mathcal
F}^j (q_0)| = |{\mathcal F}^j (\xi) - q_0|< K \eta^{j}$ for $k<i$.
Thus $e_i$ is well defined. Integrating and arguing as above, we
obtain $\gamma_i$ with $|\gamma_i(s)-\gamma_{i-1}(s)| < (K{\bf
b})^{i-1} |s| << \eta^{i}$ for all $s$ with $|s|< 2\eta$.

\medskip

To estimate the derivative with respect to ${\bf  a}$, we let
$$
q(t, s) = \gamma_n(s) + t(\gamma_{n+1}(s) - \gamma_n(s)).
$$
We have
\begin{eqnarray*}
|\frac{d}{ds} \partial_{\bf a} (\gamma_{n+1}(s)-\gamma_n(s)) | &
\leq & |\partial_{\bf a} (e_{n+1}(\gamma_{n+1}(s)) -
e_n(\gamma_{n+1}(s)))| \\
& & \ \ +
|\partial_{\bf a}(e_{n}(\gamma_{n+1}(s)) - e_n(\gamma_n(s)))| \\
& \leq & |e_{n+1}-e_n| + |\partial_t \partial_{\bf a} e_n|
|\gamma_{n+1}(s) - \gamma_n(s)| \leq K \eta^n.
\end{eqnarray*}
From this the second item of Lemma \ref{lemma-a3}(b) follows.
\hfill $\square$

\subsection{The proof of Claim \ref{claim6-th3}(b)}
We are ready to prove Claim \ref{claim6-th3}(b). First we note
that at the point of tangency, $|z| << {\bf b}^{\frac{1}{10}}$ so
it suffices for us to consider ${\mathcal A}_{\eta}$ in the place
of ${\mathcal A}$ with $\eta = {\bf b}^{\frac{1}{10}}$.

From Lemma \ref{lemma-a3}(c), we know that $\gamma_n \to
\gamma_{\infty}$ uniformly as $n \to \infty$, and from Lemma
\ref{lemma-a3}(a) we know that $\gamma_{\infty}$ is the stable
manifold of $q_0=q_m$, which we write as
\begin{equation}\label{f20-appb}
\theta = \theta_{\infty}(s, {\bf a}), \ \ \ z = z_{\infty}(s, {\bf
a}).
\end{equation}
We also have
\begin{equation}\label{f21-appb}
|\partial_{\bf a} \theta_{\infty}(s, {\bf a}) - \partial_{\bf a}
\theta_1(s, {\bf a})|, \ \ |\partial_{\bf a} \theta_{\infty}(s,
{\bf a}) -
\partial_{\bf a} \theta_1(s, {\bf a})| < 2 {\bf b}^{\frac{1}{10}}
\end{equation}
from Lemma \ref{lemma-a2}(c).

Write $\gamma_{\infty}$ using $\theta = w^s(z, {\bf a})$. We first
solve $s = s_{\infty}(z, {\bf a})$ from the second item in
(\ref{f20-appb}), then substitute to the first item in
(\ref{f20-appb}) to obtain
$$
w^s(z, {\bf a}) = \theta_{\infty}(s_{\infty}(z, {\bf a}), a).
$$
Differentiating on both side, we obtain
\begin{eqnarray*}
\partial_{\bf a} w^s(z, {\bf a}) & = & \partial_s \theta_{\infty} (s, {\bf a})
\partial_{\bf a} s_{\infty}(z, {\bf a}) + \partial_{\bf a}
\theta_{\infty}(s, {\bf a}) \\
& = & - \partial_s \theta_{\infty} (s, {\bf a})
\frac{\partial_{\bf a} z_{\infty}(s, {\bf a})}{\partial_s
z_{\infty}(s, {\bf a})} +
\partial_{\bf a} \theta_{\infty}(s, {\bf a}) \\
& = & \partial_{\bf a} w^s_1(z, {\bf a}) + {\mathcal O}({\bf
b}^{\frac{1}{10}})
\end{eqnarray*}
where $\theta = w^s_1(z, {\bf a})$ is the equation for $\gamma_1$.
To obtain the last estimate we use (\ref{f21-appb}) and the fact
that $\partial_s z_1(s, {\bf a}) > \frac{1}{2}$.

To prove Claim \ref{claim6-th3}(b), it now suffices for us to
confirm that
$$
\partial_{\bf a} w^s_1(z, {\bf a}) < \frac{1}{50}.
$$
To verify this estimate we observe
\begin{equation}\label{f22-appb}
|\partial_{\bf a}\frac{d}{dz}w^s_1(z, {\bf a})| = |\partial_{\bf
a} s^{-1}(e_1)| < K {\bf b},
\end{equation}
where $s(e_1)$ is the slope for $e_1$. This follows from a direct
computation using (\ref{f3-appb}). From (\ref{f22-appb}),
$$
|\partial_{\bf a} w^s_1(z, {\bf a})| < |\partial_{\bf a}
w^s_1(z_m, {\bf a})| + K {\bf b} \eta
$$
where $\partial_{\bf a} w^s_1(z_m, {\bf a})$ is the value of
$\partial_{\bf a} w^s_1(z, {\bf a})$ at $q_0 = q_m$. We now use
Claim \ref{claim2-th3} for $\partial_{\bf a} w^s_1(z_m, {\bf a})$.
\hfill $\square$

\end{document}